\documentclass[12pt]{amsart}

\def\Cal{\mathcal}

\def\b1{\text{\bf 1}}

\def\BB{{\Bbb B}}
\def\BC{{\Bbb C}}
\def\BP{{\Bbb P}}

\def\BP{{\Bbb P}}
\def\BQ{{\Bbb Q}}
\def\BZ{{\Bbb Z}}

\def\CA{{\Cal A}}

\def\CJ{{\Cal J}}
\def\CB{{\Cal B}}
\def\CC{{\Cal C}}

\def\CH{{\Cal H}}
\def\CL{{\Cal L}}

\def\CO{{\Cal O}}

\def\CT{{\Cal T}}

\def\CU{{\Cal U}}

\def\fh{{\frak h}}

\def\fa{{\frak a}}

\def\fF{{\frak F}}

\def\fh{{\frak h}}

\def\Ell{\text{Ell}}

\def\pr{\text{par}}
\def\iso{ \buildrel\sim\over\longrightarrow }
\def\eqd{\buildrel\text{def}\over =}

\begin{document}

\title[Vertex Algebras and  LG/CY Correspondence]
{\bf Vertex Algebras and the Landau-Ginzburg/Calabi-Yau Correspondence}
\thanks{partially supported by the National Science Foundation}

\author{Vassily Gorbounov}

\author{Fyodor Malikov}

\maketitle

\centerline{{\it To Borya Feigin on his 50th birthday}}

%\centerline{\bf Abstract}
\begin{abstract}

We construct a spectral sequence that converges to the cohomology
of the chiral de Rham complex over a Calabi-Yau hypersurface and
whose first term is a vertex algebra  closely related to the
Landau-Ginzburg orbifold. As an application, we prove an explicit
orbifold formula for the elliptic genus of Calabi-Yau
hypersurfaces.
\end{abstract}
%%%

\centerline {\bf Introduction}

\bigskip

Introduced in [MSV] for any smooth (algebraic, analytic, etc.)
manifold $X$ there is a sheaf of vertex algebras
$\Omega^{ch}_{X}$. For example, the vertex algebra of global
sections over $\BC^{N}$, $\Omega^{ch}_{\BC^{N}}(\BC^{N})$ or
simply $\Omega^{ch}(\BC^{N})$, is well known as
``$bc-\beta\gamma$-system''; it is an apparently unsophisticated
object. Despite various important contributions [B, BD, BL, KV1],
however, very little is known about the cohomology vertex algebra
$H^{*}(X,\Omega^{ch}_{X})$ even if $X$ is toric -- except perhaps
the case of $\BP^{2n}$ where at least  the character of the  space
of global sections, $H^{0}(\BP^{2n},\Omega^{ch}_{\BP^{2n}})$, has
been computed: it was shown to be equal to the elliptic genus of
$\BP^{2n}$ in  [MS].

Let
$$
\fF=\{f=0\}\subset\BP^{N-1}.
\eqno{(1)}
$$
be a Calabi-Yau hypersurface.  The present paper is devoted to an interplay
between $\Omega^{ch}(\BC^{N})$ and closely related algebras
on the one hand, and
$H^{*}(\fF,\Omega^{ch}_{\fF})$ on the other.
Let us now  formulate the main result.

{\it Preparations.}
Set
$$
\Lambda=\BZ^{N}\oplus(\BZ^{N})^{*}.
$$
Associated to this lattice in the standard manner
 there are the lattice vertex algebra
$V_{\Lambda}$ and the fermionic vertex algebra ( $bc$-system), $F_{\Lambda}$,
which is none other than the vacuum representation of the
infinite-dimensional Clifford algebra $Cl(\BC\otimes_{\BZ^{N}}\Lambda)$.
In the important paper [B] Borisov introduces the vertex algebra
$\BB_{\Lambda}=V_{\Lambda}\otimes F_{\Lambda}$ and the vertex algebra embedding
$$
\Omega^{ch}(\BC^{N})\hookrightarrow \BB_{\Lambda}.
$$
Define the $\BZ_{N}$-action
$$
\BZ_{N}\times \BC^{N}\rightarrow \BC^{N},\;
(m,\vec{x})\mapsto \exp{(2\pi i\frac{m}{N})}\vec{x}.
$$
There arises the vertex subalgebra of $\BZ_{N}$-invariants
$$
(\Omega^{ch}(\BC^{N}))^{\BZ_{N}}\subset \Omega^{ch}(\BC^{N}).
$$
Let us now warp the lattice $\Lambda$: denote by
 $\BZ^{N}_{orb}$  the sublattice
of $\BZ^{N}$ consisting $(m_{0},m_{1},..., m_{N-1})$ such that
$\sum_{j}m_{j}$ is divisible by $N$ and define
$$
\Lambda_{orb}= \BZ^{N}_{orb}\oplus (\BZ^{N}_{orb})^{*}.
$$
Just as above, there arises the vertex algebra $\BB_{\Lambda_{orb}}$.
Note that $\BB_{\Lambda}$ and $\BB_{\Lambda_{orb}}$ have a
non-empty intersection, which contains $(\Omega^{ch}(\BC^{N}))^{\BZ_{N}}$;
thus
$$
(\Omega^{ch}(\BC^{N}))^{\BZ_{N}}\hookrightarrow \BB_{\Lambda_{orb}}.
$$
Now we extend $(\Omega^{ch}(\BC^{N}))^{\BZ_{N}}$ inside
$\BB_{\Lambda_{orb}}$ to a bi-differential vertex algebra.

Let $\{X_{i}\}\subset\BZ^{N}$, $\{X_{i}^{*}\}\subset(\BZ^{N})^{*}$
be the standard dual bases. Associated to them inside
$V_{\Lambda}$ there are fields, such as $X_{i}(z)$,
$X_{i}^{*}(z)$, $e^{X_{i}^{*}}(z)$. Let the corresponding tilded
letters denote their superpartners inside $F_{\Lambda}$, e.g.,
$\tilde{X}_{i}(z)$, $\tilde{X}_{i}^{*}(z)$.

Denote
$$
X^{*}_{orb}=\frac{1}{N}(X^{*}_{0}+X^{*}_{1}+\cdots X_{N-1}^{*})
\in (\BZ^{N}_{orb})^{*}.
$$
Form
$$
\widetilde{\text{LG}}=\bigoplus_{n=0}^{\infty}
\widetilde{\text{LG}}^{(n)},\;
\widetilde{\text{LG}}^{(n)}=(\Omega^{ch}(\BC^{N}))^{\BZ_{N}}e^{nX^{*}_{orb}}.
$$
This is clearly a $\BZ_{+}$-graded subalgebra of
$\BB_{\Lambda_{orb}}$ (but not of $\BB_{\Lambda}$).

Now define two operators
$$
D_{orb}=
(\sum_{j=0}^{N-1}e^{X^{*}_{orb}}(z)\tilde{X}^{*}_{j}(z))_{(0)},
d_{LG}=df(z)_{(0)}\in \text{End}(\widetilde{\text{LG}}),
$$
where  $f$ is the polynomial appearing in (1) and $df(z)$ is computed
by using the definition
$$
d(x_{0}^{m_{0}}x_{1}^{m_{1}}\cdots x_{N-1}^{m_{N-1}})(z)=
e^{\sum_{j}m_{j}X_{j}}(z)\sum_{j}m_{j}\tilde{X}_{j}(z).
$$

These are commuting, square zero derivations of $\widetilde{\text{LG}}$;
thus we have obtained the bi-differential vertex algebra
 $(\widetilde{\text{LG}}; D_{orb}, d_{LG})$. Note that it is filtered
by the bi-differential vertex ideals
$$
\widetilde{\text{LG}}^{\geq n}=\bigoplus_{m=n}^{\infty}
\widetilde{\text{LG}}^{(m)}.
$$
Hence there arises the projective system of bi-differential vertex
algebras
$$
\widetilde{\text{LG}}^{< n}=\widetilde{\text{LG}}/
\widetilde{\text{LG}}^{\geq n}.
$$

{\bf Theorem 1.} (cf. Theorems 4.7, 4.9) {\it There is a spectral sequence }
$$
(E^{*,*}_{*}, d_{*})\Rightarrow H^{*}(\fF,\Omega^{ch}_{\fF})
$$
{\it that satisfies:}

(i)
$$
 (E^{*,i-*}_{1}, d_{1})\sim
(H^{i}_{D_{orb}}(\widetilde{\text{LG}}^{< N}), d_{LG});
$$
{\it (ii) at the conformal weight
zero component this spectral system degenerates
in the 2nd term so that}
$$
  H^{*}(\fF,\Lambda^{*}\CT_{\fF}) \iso
H^{*}(\fF,\Omega^{ch}_{\fF})_{0}
\sim (E^{*,*}_{2})_{0}=H_{d_{LG}}(\widetilde{\text{LG}}^{< N})_{0},
$$
{\it where $\Lambda^{*}\CT_{\fF}$ is the algebra of polyvector
fields over $\fF$. Further,}
$$
H_{d_{LG}}(\widetilde{\text{LG}}^{(i)})_{0}=
\left\{\aligned
\BC &\text{ if } 1\leq i\leq N-1,\\
 M_{f}^{\BZ_{N}} &\text{ if } i=0,
\endaligned\right.
\eqno{(2)}
$$
{\it where $M_{f}^{\BZ_{N}}=(\BC[x_{0},...,x_{N-1}]/<df>)^{\BZ_{N}}$
is the $\BZ_{N}$-invariant part of the Milnor ring.
$\qed$}

\bigskip

{\it Remarks.}

(i) The sign $\sim$ in item (i) means that rather than being
isomorphic the complexes are filtered, and the corresponding
graded complexes are naturally isomorphic. This is not too serious
a complication; in fact, $\sim$ is a genuine isomorphism if
$i<N-1$, and there is a 1-step filtration if $i=N-1$.

(ii) The reader familiar with previous work might
expect \newline $ H^{*}(\fF,\Omega^{*}_{\fF}) \iso
H^{*}(\fF,\Omega^{ch}_{\fF})_{0}$ instead of
 the first isomorphism
in the item (ii) of the theorem. We have indeed changed the conformal
grading and find this important; so much so that in the main body of the text
(see especially 2.3.3)
we change the terminology and notation: we write $\Lambda^{ch}\CT_{\fF}$
for $\Omega^{ch}_{\fF}$ with the changed grading and call it the
{\it the algebra of chiral polyvector fields.} $\qed$

\bigskip

Now we would like to make two points. First, let us demonstrate how this
result works.

{\it Application: an elliptic genus formula.} Let
$\text{Ell}_{\fF}(\tau,s)$ be the 2-variable elliptic genus of $\fF$
as defined, for example, in [BL] or [KYY].

Introduce
$$
E(\tau,s)=\prod_{n=0}^{\infty}\frac{\left(1-e^{2\pi
i\left(\left(n+1\right)\tau+\left(1-1/N\right)s\right)}\right)^{N}\left(1-e^{2\pi
i\left(n\tau+\left(-1+1/N\right)s\right)}\right)^{N}}{\left(1-e^{2\pi
i\left(\left(n+1\right)\tau+s/N\right)}\right)^{N}\left(1-e^{2\pi
i\left(n\tau-s/N\right)}\right)^{N}}.
$$
It follows easily from Theorem 1 (see 4.11- 4.12) that
$$
\Ell_{\fF}(\tau,s)= \frac{1}{N}\sum_{l=0}^{N-1}\sum_{j=0}^{N-1}
e^{\pi
i\left(N-2\right)\left\{-js+\left(j^{2}-j\right)\tau+j^{2}\right\}}E(\tau,s-j\tau-l).
\eqno{(3)}
$$
The structure of this formula is rather clear: the infinite
product $E(\tau,s)$ reflects the polynomial nature of the space
$\Omega^{ch}(\BC^{N})$ of which it is indeed the Euler character,
 the summation with respect to $l$ extracts the
$\BZ_{N}$-invariants, and the summation w.r.t. $j$ reminds of the
summation over ``twisted sectors'' because the change of variable
$s\mapsto s-j\tau$ is reminiscent of the spectral flow. All of
this smacks of an orbifold,  and indeed formula (3) was proposed
in [KYY] as the elliptic genus of the Landau-Ginzburg orbifold.
Furthermore, it was shown in [KYY] that the specialization of the
r.h.s. of (3) to $\tau=i\infty$ or $s=0$ gives Vafa's orbifold
formulas  for the Poincare polynomial and the Euler characteristic
of the Fermat hypersurface respectively. These formulas have been
proved and further discussed in [OR, R]. Of course, both easily
follow from (3).

 This brings about the 2nd point we would like to
make.

{\it Landau-Ginzburg orbifold interpretation.} The Landau-Ginzburg model
is associated to an affine manifold and a function over it known
as superpotential. In the case
where the manifold is $\BC^{N}$ and the function
is a homogeneous polynomial $f$ with a unique singularity at 0, Witten's discovery [W2] can
 perhaps be formulated
 in the language accessible to us as follows:

{\it 1) There is an action of the $N=2$
superconformal algebra on $\Omega^{ch}(\BC^{N})$ such that
it commutes with the differential $df(z)_{(0)}$.

2) The cohomology  vertex algebra $H_{df(z)_{(0)}}(\Omega^{ch}(\BC^{N}))$
with thus defined action of  the $N=2$
superconformal algebra is
the chiral algebra attached to the  Landau-Ginzburg model with
superpotential $f$. }

We find it convenient
not to pass to the cohomology but to declare the Landau-Ginzburg
model to be the differential vertex algebra $(\Omega^{ch}(\BC^{N}), df(z)_{(0)})$
with  the above fixed $N=2$ superconformal algebra action.

\bigskip
{\it Remark.}
We shall
argue in 5.1.5 that alternatively one can
 think of\newline $(\Omega^{ch}(\BC^{N}), df(z)_{(0)})$ as
 the ``right definition'' of the chiral de Rham complex
$\Omega^{ch}_{\text{Spec}M_{f}}$
over the spectrum of the Milnor ring. $\qed$
\bigskip

Next consider the space $\Omega^{ch}(\BC^{N})e^{nX^{*}_{orb}}$. It
does not belong to $\BB_{\Lambda_{orb}}$ but carries a canonical
structure of a {\it twisted}  $\Omega^{ch}(\BC^{N})$-module, and
this is synonymous to being a twisted sector. Therefore, taking
the direct sum of these, $0\leq n\leq N-1$, and then extracting
$\BZ_{N}$-invariants corresponds accurately with what physicists
call the Landau-Ginzburg orbifold, see e.g. [V]. In the notation
we have adopted, the formula describing the outcome of this
process is $(\widetilde{\text{LG}}^{< N}; d_{LG})$; see sect. 5
for details. Thus item (i) of Theorem 1 can be interpreted as
follows:

{\it there is a spectral sequence abutting to $ H^{*}(\fF,\Omega^{ch}_{\fF})$
whose 1st term is isomorphic to the $D_{orb}$-cohomology of the
 Landau-Ginzburg orbifold
$H_{D_{orb}}(\widetilde{\text{LG}}^{< N})$.}

One can say that the Landau-Ginzburg orbifold $(\widetilde{\text{LG}}^{< N}; d_{LG})$
 approximates the vertex algebra $ H^{*}(\fF,\Omega^{ch}_{\fF})$.
This approximation is consistent with the $N=2$ supersymmetry.
Indeed, on the one hand,
the $N=2$ superconformal algebra action on the Landau-Ginzburg model
commutes with $D_{orb}$ and all higher differentials and thus descends to
an action on $ H^{*}(\fF,\Omega^{ch}_{\fF})$. On the other hand, since
$\fF$ is Calabi-Yau, $ H^{*}(\fF,\Omega^{ch}_{\fF})$ carries a canonical
$N=2$ superconformal algebra action  [MSV]. A direct computation
(Lemmas 4.10.1, 5.1.1, 5.2.14) shows that

{\it both the $N=2$ superconformal algebra actions coincide.}

Furthermore, not only  one of the spaces involved in assertion (i)
of Theorem 1, but both assertions (i, ii) themselves are
reminiscent of some of the important developments in string
theory. In order to explain this we shall have to pluck courage
and discuss a little more of physics.

It seems that except for the torus case, see e.g. a mathematical
 exposition in [KO], the  space of states  of the model describing the
 string propagation on a manifold is unknown even as a vector space to say
nothing about its algebraic, vertex or otherwise, structure. One
 striking result towards understanding what this might be is Gepner's
model proposed in [G].
Gepner's paper, a combination of guesswork
and computational {\it tour de force}, is not an easy reading. More conceptual
approach emerged soon afterwards, e.g. [V,VW], proclaiming
that the Landau-Ginzburg orbifold is equivalent (in this or that sense)
 to the string theory
on Calabi-Yau hypersurfaces in weighted projective spaces.
(So far as we can tell, apart from both the
theories carrying an $N=2$ superconformal algebra
action with the same central charge
and integral $U(1)$-charges, most of the supporting evidence amounted to the
isomorphism of chiral rings -- exactly as in Theorem 1 (ii).)
This activity
seems to have been crowned by Witten's paper [W1] where it is asserted, and
we cite, ``that rather than Landau-Ginzburg being "equivalent" to
Calabi-Yau, they are two different phases of the same system.''

Now, if one is allowed to think of a phase transition as a family where
at certain values of the parameter something happens, then it seems that
 spectral sequences might be relevant.
For example, the spectral sequence of Theorem 1 comes from a
double complex equipped with two differentials. Let us denote them
for the purposes of introduction by $d_{+}$ and $d_{-}$. The
vertex algebra $H^{*}(\fF,\Omega^{ch}_{\fF})$ is the cohomology of
the total differential $d_{+}+d_{-}$, and
$H_{D_{orb}}(\widetilde{\text{LG}}^{< N})$ arises  as the
$d_{-}$-cohomology. Introduce a parameter, $t$, and form the
differential $td_{+}+d_{-}$. This defines a family of complexes
over a line such that  at $t=0$
 the cohomology is $H_{D_{orb}}(\widetilde{\text{LG}}^{< N})$, and elsewhere
 it is $H^{*}(\fF,\Omega^{ch}_{\fF})$.

Furthermore, the geometric background as explained in sect. 4 of [W1]
is very similar to that we use in the proof of Theorem 1.

{\it Sketch of proof.} The proof is based on the computation of
two spectral sequences, both due to [B]. The first allows, in a sense,
to replace $\fF$ with the canonical line bundle, $\CL^{*}$, over $\BP^{N-1}$.
The 1st term of this sequence
is $H^{*}(\CL^{*},\Omega^{ch}_{\CL^{*}})$ and the corresponding
differential $d_{1}$ has the meaning of the chiral Koszul differential
\footnote{this correspondence between $H^{*}(\CL^{*},\Omega^{ch}_{\CL^{*}})$
and $H^{*}(\fF,\Omega^{ch}_{\fF})$ reminds one of the main result in [S]}.
Most of Theorem 1 is about identification of the complex
$(H^{*}(\CL^{*},\Omega^{ch}_{\CL^{*}}), d_{1})$.

We identify this space in the following way. Along with $\CL^{*}$ consider
$\CL^{*}-0$ obtained by deleting the zero section. The \v Cech complex
that computes the desired $H^{*}(\CL^{*},\Omega^{ch}_{\CL^{*}})$
naturally embeds into the analogous \v Cech complex over $\CL^{*}-0$.
Write this down schematically as
$$
\check{C}(\CL^{*})\hookrightarrow \check{C}(\CL^{*}-0).
$$
Proposed in [B] there is a vertex algebra resolution that allows to extend the latter
embedding to a resolution of complexes:
$$
\check{C}(\CL^{*})\hookrightarrow \check{C}(\CL^{*}-0)
\rightarrow \check{C}(\CL^{*}-0)^{(1)}
\rightarrow \check{C}(\CL^{*}-0)^{(2)}\rightarrow\cdots.
$$
This resolution serves the same purpose as the Cousin resolution
but has a  different flavor: its terms
$\check{C}(\CL^{*}-0)^{(j)}$ are identified with each other as
vector spaces, but are spectral flow transforms of each other as
vertex modules. The latter property is responsible for the
occurrence of the orbifold twisted sectors.

 There arises then the
bi-complex
$$
0\rightarrow \check{C}(\CL^{*}-0)
\rightarrow \check{C}(\CL^{*}-0)^{(1)}
\rightarrow \check{C}(\CL^{*}-0)^{(2)}\rightarrow\cdots.
\eqno{(4)}
$$
 whose total cohomology is the desired
$H^{*}(\CL^{*},\Omega^{ch}_{\CL^{*}})$. Finally, and this is the geometry bit
 reminiscent of [W1], there is an isomorphism
$$
(\BC^{N}-0)/\BZ_{N}\iso \CL^{*}-0.
$$
Pulling bi-complex (4) back onto $(\BC^{N}-0)/\BZ_{N}$ and further
writing it down in terms of the coordinates on the universal
covering space $\BC^{N}-0$
allows to compute all the terms of the corresponding spectral
sequence and thus identify $(H^{*}(\CL^{*},\Omega^{ch}_{\CL^{*}}), d_{1})$
with $(H_{D_{orb}}(\widetilde{\text{LG}}^{< N}), d_{LG})$
 as asserted in Theorem (i). (Note that these two pull-backs
are made possible by the {\it naturality}  property of $\Omega^{ch}_{X}$, see 2.1.) $\qed$

The ``conformal weight zero component of this argument''
gives a self-contained computation of the  cohomology algebra of
polyvector fields $H^{*}(\fF,\Lambda^{*}\CT_{\fF})$ along
with explicit ``vertex'' formulas for the cocycles representing the cohomology
classes,
see 4.13. Classically, the computation uses  $H^{*}(\fF,\BC)$
obtained in [Gr], the Serre
duality, and variations of the Hodge structure [Don, Theorem 2.2] -- and even then powers of
the class representing  hyperplane sections require special, although
 simple, treatment. In our approach, it is
the other way around: the fact that the Milnor ring is realized
inside $H^{*}(\fF,\Lambda^{*}\CT_{\fF})$ is almost obvious, even
its Koszul resolution arises naturally at the 1st term of the
spectral sequence. On the other hand,
 the Serre duals of powers of the hyperplane section are ``interesting'' because
 they are produced
by the twisted sectors as  indicated in Theorem 1 (ii),
 1st line in (2). We restore the multiplicative structure of
 $H^{*}(\fF,\Lambda^{*}\CT_{\fF})$ by vertex algebra methods  twice:
first, in 4.13 in the context of the proof of Theorem 1; second,
in 5.2.19 by way of testing the vertex algebra structure on the
space $H_{d_{LG}}(\widetilde{\text{LG}}^{< N})$  we propose in
Theorem 5.2.18 in the case of the ``diagonal'' $f$. We believe
that this vertex algebra structure is exactly Witten's chiral
algebra of the Landau-Ginzburg orbifold. Note that, in general,
orbifoldizing multiplicative structures is  a problem. The vertex
algebra structure of Theorem 5.2.18 owes its existence to a
remarkable periodicity property of the category of unitary modules
over the $N=2$ superconformal algebra [FS].

\begin{sloppypar}
{\bf Acknowledgements.} The indebtedness of this work to Borisov's
constructions [B]
should be clear to any reader. We gratefully acknowledge many
illuminating conversations with V.Batyrev, A.Gerasimov, A.Givental,
C.Hertling, M.Kapranov, R.Kaufmann, Yu.I.Manin, A.Semikhatov,
 A.Vaintrob. Special
thanks go to V.Schechtman who participated at the early stages
of this research and was the first enthusiast of the algebra
of chiral polyvector fields. This work was started in 2002 when we were visiting
the Max-Planck-Institut f\"ur Mathematik in Bonn and finished
a year later at the same place. We are grateful to the institute for excellent working conditions.
\end{sloppypar}

\bigskip
\bigskip
\centerline{{\bf 1. Vertex algebras }}

\bigskip
This section is only a collection of well-known facts and examples that
will be needed in the sequel
and the reader may want to consult either [K]
or [FB-Z] for more detail. We would like to single out sect. 1.12 on the spectral
 flow, the notion that seems to be not too popular in mathematics
literature but is essential for understanding orbifolds and will reappear
several times
in sect. 2.3.5, 3.10, 4.6, 5.2.15. Our treatment of the spectral flow is greatly
influenced by [LVW].

{\bf 1.1.} A vector space $V$  is called a supervector space if it
is $\BZ_{2}$-graded, that is, $V=V^{(0)}\oplus V^{(1)}$. We define
the parity $\pr(a)$ of $a\in V$ so that $\pr(a)=\epsilon$ if and
only if $a\in V^{(\epsilon)}$. If $V$ and $W$ are supervector
spaces, then $V\otimes W$ is also with $\pr(a\otimes
b)=\pr(a)+\pr(b)$, and so is $\text{Hom}_{\BC}(V,W)$.

Given a supervector space $V$, let $\text{Field}(V)$ be the  subspace
of $\text{End}(V)[[z,z^{\pm 1}]]$ consisting of such formal series
$x(z)=\sum_{n\in\BZ}x_{(n)}z^{-n-1}$ that for any  $a\in V$
$$
x_{(n)}a=0\text{ if } n>>0.
\eqno{(1.1.1)}
$$

{\bf 1.2. Definition.} A vertex algebra is a supervector space $V$
with a distinguished element $\b1\in V$ called {\it vacuum} and a parity
preserving map
$Y(.,z):\; V\rightarrow \text{Field}(V)$,
$Y(a,z)=\sum_{n\in\BZ}a_{(n)}z^{-n-1}$, such that the following axioms hold:

(i) {\it vacuum}:
$$
Y(\b1,z)=Id_{V},\; a_{(-1)}\b1=a;
\eqno{(1.2.1)}
$$
(ii) {\it Borcherds identity:}
 for any $a,b\in V$ and any rational function $F(z,w)$ in $z$, $w$
with poles only at $z=0$, $w=0$, $z-w=0$
$$
\aligned
&\text{Res}_{z-w}Y(Y(a,z-w)b,w)i_{w,z-w}F(z,w)\\
&=
\text{Res}_{z}\left(Y(a,z)Y(b,w)i_{z,w}F(z,w)-
(-1)^{\pr(a)\pr(b)}Y(b,w)Y(a,z)i_{w,z}F(z,w)\right).
\endaligned
\eqno{(1.2.2)}
$$
$\qed$

In the  latter formula the standard notation is used: $\text{Res}_{t}$ means
the coefficient of $t^{-1}$ in the indicated formal Laurent expansion;
$i_{\bullet,\bullet}$ specifies exactly which Laurent expansion is to be used, e.g.
$i_{z,w}$ stands for the expansion in the
domain  $|w|<|z|$, $i_{w,z-w}$ for that in the domain $|z-w|<|w|$, etc.

One thinks of $a_{(n)}$ as the `` n-th multiplication by $a$'', so
there arises a family of multiplications
$$
_{(n)}:\; V\otimes V\rightarrow V.
\eqno{(1.2.3)}
$$

{\bf 1.3. W- and cohomology vertex algebras.}

If $F=1$, then  (1.2.2) gives
$$
(a_{(0)}b)_{(n)}=[a_{(0)},b_{(n)}],\; n\in\BZ.
\eqno{(1.3.1)}
$$
In other words, for any $a\in V$, $a_{(0)}\in\text{End}V$ is a derivation
of all the multiplications. Hence,
$$
\text{Ker}a_{(0)}\subset V
\eqno{(1.3.2)}
$$
is a vertex subalgebra known as a $W$-algebra.

Furthermore, suppose $a\in V$ is odd and $a_{(0)}a=0$. Then
(1.3.1) implies that  $a_{(0)}^{2}=0$ and $\text{Im}a_{(0)}\subset\text{Ker}a_{(0)}$
is an ideal. Therefore, the cohomology
$$
H_{a_{(0)}}(V) \eqd \text{Ker}a_{(0)}/\text{Im}a_{(0)}
\eqno{(1.3.3)}
$$
carries a canonical vertex algebra structure. The vertex algebras to be used in this
text will mostly be either W- or cohomology vertex algebras.

{\bf 1.4. Chiral rings.}
 Suppose that $V$ is graded so that
$$
V=\oplus_{n=0}^{\infty} V_{n}\text{ and }
V_{n\;(r)}V_{m}\subset V_{n+m-r-1}.
\eqno{(1.4.1)}
$$
The grading satisfying this condition will be called {\it conformal}.

One can show that if (1.4.1) is valid, then
$$
_{(-1)}: V_{0}\otimes V_{0}\rightarrow V_{0}
$$
is associative and supercommutative. In the context of the unitary
 $N=2$ supersymmetry, supercommutative associative
algebras attached to graded vertex algebras in this way are often called
{\it chiral rings} [LVW] -- not to be confused
with chiral algebras although these rings are indeed algebras.
 We shall take the liberty to call these rings
chiral in any case.

{\bf 1.5. Remarks.}

(i) Let $\delta(z-w)=\sum_{n\in\BZ}z^{n}w^{-n-1}$.
It follows from (1.2.2) that
$[Y(a,z), Y(b,w)]$ is local, that is, equals a linear combination
of the delta-function derivatives, $\partial_{w}^{n}\delta(z-w)$,
 over fields in $w$.

(ii) A vertex algebra $V$ is said to be generated by a collection of fields
$Y(a_{\alpha},z)$, $\{a_{\alpha}\}\subset V$ if $V$ is the linear span
of (non-commutative) monomials in $(a_{\alpha})_{(j)}$ applied to vacuum $\b1\in V$.
The important {\it reconstruction theorem}, e.g. [K, Theorem 4.5], says, and we are omitting
some details, that this can be reversed: if there is a collection of
mutually local fields $v_{\alpha}(z)$ which generate
$V$ from a fixed vector, then $V$ carries a unique vertex algebra structure
such that $v_{\alpha}(z)=Y(v_{\alpha},z)$ for some $v_{\alpha}\in V$.
Because of this we will allow ourselves in our list of well-known examples,
which we are about to begin, to fix only  a space $V$ and a collection of
mutually local fields that generate this space. Typically, we shall have a Lie
algebra, a collection of fields with values in this algebra, and a representation
of this algebra such that  nilpotency condition (1.1.1) is satisfied by
the fields.

(iii) It should be clear what a vertex algebra  homomorphism is. As in (ii), speaking
of homomorphisms we shall often specify only images of generating fields.

Now to some basic excamples.

{\bf 1.6. $bc$-system.}
 Let $Cl$  be
 the Lie superalgebra with basis $b_{(i)}, c_{(i)}$,
$i\in\BZ$ (all odd), and $C$ (even) and commutation relations
$$
[b_{(i)},c_{(j)}]=\delta_{i,-j-1}C,\; [C,c_{i}]=[C,b_{i}]=0,\; i,j\in\BZ.
\eqno{(1.6.1)}
$$
Let
$$
F=\text{Ind}_{Cl}^{Cl_{+}}\BC,
\eqno{(1.6.2)}
$$
where $Cl_{+}$ is the Lie subalgebra spanned by $x_{(i)}$, $C$, $i\geq 0$,
and $\BC$ is an   $Cl_{+}$-module where $x_{(i)}$'s act by 0 and $C$ as
multiplication by 1.
In terms of fields $x(z)\sum_{i\in\BZ}x_{(i)}z^{-i-1}$, $x=b$ or $c$, (1.6.1)
becomes:
$$
[b(z),c(z)]=\delta(z-w).
\eqno{(1.6.3)}
$$
Hence the vertex algebra structure on $F$, see 1.5.

A little more generally, to any purely odd $\BC$-vector $W$ with
a non-degenerate symmetric form $(.,.)$ one can attach
the Lie superalgebra $Cl(W)=W\otimes\BC[t,t^{-1}]\oplus\BC C$ with
the following bracket: if
$x(z)=\sum_{i\in\BZ}(x\otimes t^{i})z^{-i-1}$, then
$$
[x(z),y(z)]=\delta(z-w)(x,y)C,\; [x(z),C]=0.
\eqno{(1.6.4)}
$$
The corresponding vertex algebra is
$$
F_{W}=\text{Ind}_{Cl(W)_{+}}^{Cl(W)}\BC,
\eqno{(1.6.5)}
$$
where $Cl(W)_{+}=W\otimes \BC[t]\oplus\BC C$,
$W\otimes \BC[t]$ operates on $\BC$ by 0, and $C$ by 1.

{\bf 1.7. $\beta\gamma$- and $bc-\beta\gamma$-system.}
 The $\beta\gamma$-system  is obtained by the ``parity change''
functor applied to the beginning of 1.6: the even Lie algebra $\fa$  is spanned
by $\beta_{(i)}$, $\gamma_{(i)}$, $C$, the bracket is
$$
[\beta(z),\gamma(w)]=-[\gamma(z),\beta(w)]=\delta(z-w)C,\; [C, x(z)]=0.
\eqno{(1.7.1)}
$$
The vertex algebra, $B$, is defined in the same way as $F$, see (1.6.2).

If $V$ and $W$ are vertex algebras, then $V\otimes W$ carries the standard
vertex algebra structure. Denote
$$
FB=F\otimes B.
\eqno{(1.7.2)}
$$
This algebra and its modifications  are local models for the chiral de Rham
complex.

{\bf 1.8. Heisenberg algebra.} Let now $\fh$ be a purely even
vector space with a non-degenerate symmetric form $(.,.)$. There
arises then the Heisenberg Lie algebra
$$
\hat{\fh}=\fh\otimes\BC[t,t^{-1}]\oplus\BC\cdot C
$$
with bracket defined by letting fields be
$x(z)=\sum_{i\in\BZ}(x\otimes t^{i})z^{-i-1}$ and then setting
$$
[a(z),b(w)]=(a,b)\partial_{w}\delta(z-w),\; [C,a(z)]=0.
\eqno{(1.8.1)}
$$
The vertex algebra attached to this Lie algebra is
$$
V(\fh)=\text{Ind}_{\hat{\fh}_{+}}^{\hat{\fh}}\BC,
\eqno{(1.8.2)}
$$
where $\hat{\fh}_{+}=\fh\otimes \BC[t]\oplus\BC\cdot C$,
$\hat{\fh}_{+}=\fh\otimes \BC[t]$ perates on $\BC$ by zero, $C$ by 1.

{\bf 1.9. Lattice vertex algebras.} We shall neeed a lattice
$L$, that is, a free abelian group with integral bilinear
form $(.,.)$ and a 2-cocycle
$$
\epsilon:\; L\times L\rightarrow \BC^{*}.
\eqno{(1.9.1)}
$$
There arise the group algebra $\BC[L]$ with multiplication
$e^{\alpha}\cdot e^{\beta}=e^{\alpha+\beta}$, $\alpha,\beta\in L$, and the twisted  group algebra, $\BC_{\epsilon}[L]$, equal
to $\BC[L]$ as a vector space but with twisted multiplication:
$$
e^{\alpha}\cdot_{\epsilon}
e^{\beta}=\epsilon(\alpha,\beta)e^{\alpha+\beta}.
\eqno{(1.9.2)}
$$
Let $\fh_{L}=\BC\otimes_{\BZ} L$. There arises the Heisenberg vertex
algebra $V(\fh_{L})$, see (1.8.2).

As a vector space, the lattice vertex algebra is defined by
$$
V_{L}=V(\fh_{L})\otimes \BC_{\epsilon}[L],\;
\pr (V(\fh_{L})\otimes e^{\alpha})\equiv (\alpha,\alpha)\text{ mod }2,
\eqno{(1.9.3)}
$$
where $\pr$ means the parity, see 1.1.

This vertex algebra is generated by the familiar fields
$x(z)=\sum_{i\in\BZ}(x\otimes t^{i})z^{-i-1}$ attached to $x_{(-1)}\otimes 1$
and the celebrated vertex operators
$$
e^{\alpha}(z)=e^{\alpha}
\exp{(\sum_{j<0}\frac{\alpha_{(j)}}{-j}z^{-j})}
\exp{(\sum_{j>0}\frac{\alpha_{(j)}}{-j}z^{-j})}z^{\alpha_{(0)}},
\eqno{(1.9.4)}
$$
attached to $1\otimes e^{\alpha}$.

The action of $x_{(i)}$, $i\neq 0$, ignores $\BC_{\epsilon}(L)$,
the action of $x_{(0)}$ is uniquely determined by
$$
x_{(0)}(1\otimes e^{\alpha})=(\alpha,x)\otimes e^{\alpha}.
\eqno{(1.9.5)}
$$
The following commutator and OPE formulas are valid:
$$
[x(z),e^{\alpha}(w)]=\delta(z-w)(x,\alpha)e^{\alpha}(w),
\eqno{(1.9.6)}
$$
$$
e^{\alpha}(z)e^{\beta}(w)=(z-w)^{(\alpha,\beta)}:e^{\alpha}(z)e^{\beta}(w):,
\eqno{(1.9.7)}
$$
and the reader is advised to consult [K, (5.4.5b)] for the meaning of the
two-variable field $:e^{\alpha}(z)e^{\beta}(w):$. Note that (1.9.7) allows
to compute all operations $(e^{\alpha})_{(n)}e^{\beta}$.
We shall need the following particular cases:
$$
\aligned
(e^{\alpha})_{(-1)}e^{\beta}=&\lim_{z\rightarrow w}
(z-w)^{(\alpha,\beta)}:e^{\alpha}(z)e^{\beta}(w):=\\
&\left\{\aligned
0&\text{ if }(\alpha,\beta)>0\\
\epsilon(\alpha,\beta)e^{\alpha+\beta}&\text{ if }(\alpha,\beta)=0,
\endaligned\right.
\endaligned
\eqno{(1.9.8)}
$$
$$
[e^{\alpha}(z),e^{\beta}(w)]=0\text{ if } (\alpha,\beta)\geq 0.
\eqno{(1.9.9)}
$$

For the future use let us mention that for any sub-semigroup
$M\subset L$ there arises the vertex subalgebra
$$
V_{M,L} \buildrel\text{def}\over = V(\fh_{L})\otimes \BC_{\epsilon}[M]\subset V_{L}.
\eqno{(1.9.10)}
$$
naturally graded by $M$.

{\bf 1.10. $N=2$ super-Virasoro algebra.}
The celebrated $N=2$ super-Virasoro algebra, to be denoted $N2$ following
[K], is a supervector space with basis  $G_{(n)}$, $Q_{(n)}$, $n\in\BZ$ (all odd), $L_{(n)}$,
$J_{(n)}$, $n\in\BZ$, $C$ (all even), and bracket
$$
\aligned
[L(z),L(w)]&=2\partial_{w}\delta(z-w)L(w)+ \delta(z-w)L(w)'\\
[J(z),J(w)]&=\partial_{w}\delta(z-w)C/3,
\endaligned
\eqno{(1.10.1a)}
$$
$$
\aligned
[L(z),G(w)]&=2\partial_{w}\delta(z-w)G(w)+
\delta(z-w)G(w)',\\
 [J(z),G(w)]&=\delta(z-w)G(w),
\endaligned
\eqno{(1.10.1b)}
$$
$$
\aligned
[L(z),Q(w)]&= \partial_{w}\delta(z-w)Q(w)+
\delta(z-w)Q(w)',\\
 [J(z),Q(w)]&=-\delta(z-w)Q(w),
\endaligned
\eqno{(1.10.1c)}
$$
$$
[L(z),J(w)]= \partial_{w}^{2}\delta(z-w)\frac{C}{6}+
\partial_{w}\delta(z-w)J(w)+
\delta(z-w)J(w)',
\eqno{(1.10.1d)}
$$
$$
[Q(z),G(w)]=  \partial_{w}^{2}\delta(z-w)\frac{C}{6}
-\partial_{w}\delta(z-w)J(w)+
\delta(z-w)L(w).
\eqno{(1.10.1e)}
$$
The vertex algebra structure is carried by the following $N2$-module:
$$
V(N2)_{c}=\text{Ind}_{N2_{\geq}}^{N2}\BC_{c}.
\eqno{(1.10.2)}
$$
where $N2_{\geq}$ is the subalgebra linearly spanned by
$G_{(n)},L_{(n)},Q_{(n)},J_{(n)},C$, $n\geq 0$, and on $\BC_{c}$
$G_{(n)},L_{(n)},Q_{(n)},J_{(n)}$ operate by 0, and $C$ as multiplication
by $c$.

{\bf 1.10.1.} {\it Definition.} An $N2$-structure on a vertex algebra $W$
is a vertex algebra homomorphism $V(N2)_{c}\rightarrow W$. $\qed$

Note that the field $L(z)$ generates the Virasoro algebra; thus an $N2$-structure on a vertex algebra induces a conformal structure, and the
grading by eigenvalues of $L_{(1)}$ is conformal, cf. (1.4.1).
\bigskip

{\bf 1.11. Automorphisms.}
It is obvious that if $W$ is a vector space with a symmetric
 non-degenerate bilinear
form, then there is an embedding
$$
\text{O}(W)\hookrightarrow \text{Aut}F_{W},
\eqno{(1.11.1)}
$$
where $\text{O}(W)$ is the orthogonal group and $F_{W}$ a vertex
algebra defined in (1.6.5); this comes from the standard action of
$\text{O}(W)$ on the Clifford Lie algebra $Cl(W)$.

The analogous construction with $\text{O}(W)$ replaced with
$\text{Aut}(L)$ and $F_{W}$  with $V_{L}$ does not quite work because
of the cocycle (1.9.1). Here is one trivial observation: if we let
$\text{Aut}_{\epsilon}(L)$ be the subgroup of $\text{Aut}(L)$ stabilizing
$\epsilon$, then there is an (obvious) embedding:
$$
\text{Aut}_{\epsilon}(L)\hookrightarrow \text{Aut}V_{L}.
\eqno{(1.11.2)}
$$

The Lie algebra $N2$ affords an exceptional, mirror symmetry automorphism
$$
Q(z)\mapsto G(z), G(z)\mapsto Q(z), J(z)\mapsto -J(z),
L(z)\mapsto L(z)+J(z)'.
\eqno{(1.11.3)}
$$

{\bf 1.12. Spectral flow.} In all our examples except for $V_{L}$
vertex algebras have come from infinite dimensional Lie algebras.
One feature these Lie algebras have in common is that they admit a
{\it spectral flow}.

Define for any $n\in\BZ$ a linear transformation:
$$
S_{n}:\;\;
\aligned
Cl\rightarrow Cl,\; \text{s.t. }&b(z)\rightarrow b(z)z^{-n},
 c(z)\rightarrow c(z)z^{n},\\
\fa\rightarrow\fa,\; \text{s.t. }&\beta(z)\rightarrow \beta(z)z^{-n},
 \gamma(z)\rightarrow \gamma(z)z^{n},\\
N2\rightarrow N2,\; \text{s.t. }&Q(z)\rightarrow Q(z)z^{n},
 G(z)\rightarrow G(z)z^{-n},\\
                    & J(z)\mapsto J(z)-\frac{1}{z}\frac{nC}{3},\\
        & L(z)\mapsto L(z)-\frac{1}{z} nJ(z)+\frac{1}{z^{2}}n(n-1)\frac{C}{6},
\endaligned
\eqno{(1.12.1)}
$$
where we abused the notation by letting the same letter stand
for the maps of different spaces, $Cl$, $\fa$, $N2$, defined in 1.6,7,10 resp.
We hope this will not lead to confusion.
An untiring reader will check that in each of the cases,
$S_{n}$ is an automorphism of the Lie algebra in question.

Maps (1.12.1) generate a $\BZ$-action on each of the algebras known
as the spectral flow. In each of the cases, therefore, there arises
a family of functors on the category of modules
$$
S_{n}: \text{Mod}\rightarrow\text{Mod},\;
M\mapsto S_{n}(M),
\eqno{(1.12.2)}
$$
action on $S_{n}(M)$ being defined by precomposing that on $M$ with
$S_{n}$ of (1.12.1).

The origin of  spectral flows (1.12.1) belongs to lattice vertex
algebras. In order to see this, let $V_{L}$ be a lattice vertex
algebra and $\text{Lie}V_{L}$ the linear span inside
$\text{End}V_{L}$ of the coefficients of the fields $v(z)$, $v\in
V_{L}$. It is well known [K, F-BZ] that
 $\text{Lie}V_{L}$ is a Lie subalgebra of $\text{End}V_{L}$.

Furthermore, if $M\subset L$ is a sub-semigroup, then we have, cf.
(1.9.10),
$$
V_{M,L}\subset V_{L},\; \text{Lie}V_{M,L}\subset \text{Lie}V_{L},
\eqno{(1.12.3)}
$$
Let
$$
e^{\alpha}: V_{L}\rightarrow V_{L}
\eqno{(1.12.4)}
$$
be multiplication by $e^{\alpha}$. If the restriction of
the cocycle $\epsilon(.,.)$ to $M\subset L$ is trivial, that is,
$$
\epsilon(M,.)=1,
\eqno{(1.12.5)}
$$
then the conjugation by
map (1.12.4)
defines an automorphism
$$
S_{\alpha}: \text{Lie}V_{M,L}\rightarrow \text{Lie}V_{M,L},\;
X\mapsto (e^{\alpha})^{-1}\circ X\circ e^{\alpha};\; \alpha\in L.
\eqno{(1.12.6)}
$$
For example, under this map
$$
e^{\beta}(z)\mapsto  e^{\beta}(z)z^{(\alpha,\beta)},\; \beta\in M
\eqno{(1.12.7)}
$$
cf. (1.12.1); the desired power of $z$ owes its appearance to the
factor $z^{\alpha_{(0)}}$ in (1.9.4). Likewise,
$$
x(z)\mapsto x(z)+\frac{(\alpha,x)}{z}.
\eqno{(1.12.8)}
$$
All spectral flows (1.12.1) are  obtained as follows: embed the
corresponding vertex algebra into an appropriate $V_{L}$, thus
obtain an embedding of the corresponding Lie algebra
$\text{Lie}(\bullet)\rightarrow \text{Lie}V_{L}$, and then
restrict  ``spectral flow in the direction $\alpha$'' (1.12.6) to
the image. This operation will be of importance for us in
 3.10,  5.2.15. Here is another such example.

{\bf 1.13. Boson-fermion correspondence.} Let $\BZ$ be the standard 1-dimensional lattice:
this means that if we let $\chi$ be the generator,
then $(\chi,\chi)=1$. There arise $V_{\BZ}$,
where $\epsilon(.,.)=1$, and the famous vertex algebra isomorphism:
$$
\aligned
&F\iso V_{\BZ},\; b(z)\mapsto e^{\chi}(z),\\
 &c(z)\mapsto e^{-\chi}(z), :b(z)c(z):\mapsto \chi(z)
\endaligned
\eqno{(1.13.1)}
$$
The interested reader will check that $S_{N}|_{Cl}$ of (1.12.1) is indeed
implemented by conjugation with $e^{-n\chi}$.

Likewise
$$
\aligned
&F^{\otimes n}\iso V_{\BZ^{n}},\; b_{i}(z)\mapsto e^{\chi_{i}}(z),\\
 &c_{i}(z)\mapsto e^{-\chi_{i}}(z), :b_{i}(z)c_{i}(z):
\mapsto \chi_{i}(z),
\endaligned
\eqno{(1.13.2)}
$$
where the cocycle on $V_{\BZ^{n}}$ is chosen so as to ensure that
$V_{\BZ^{n}}\iso V_{\BZ}^{\otimes n}$.

\bigskip

\centerline{{\bf 2. The algebra of chiral polyvector fields }}
\bigskip
This section is a reminder on the chiral de Rham complex. Our exposition
is close to [MSV] but has been influenced by [GMS]. We would like to single
out 2.3.3, where the title is clarified, and 2.3.5, where a simple cohomology computation is carried out; the results of this computation will
play an important role in the proof of Theorem 4.7. Sect. 2.4 is an exposition
of a result of [B].

{\bf 2.1.} Suppose we have a family of sheaves of vector
spaces $\CA_{X}$, one for each
smooth algebraic manifold $X$. We shall call $\CA_{X}$ {\it natural} if
for any \'etale morphism $\phi: Y\rightarrow X$ there is a sheaf embedding
$\CA(\phi): \phi^{-1}\CA_{X}\hookrightarrow\CA_{Y}$ such that
given a diagram:
$$
X\buildrel \phi\over\rightarrow Y\buildrel \psi\over\rightarrow Z
$$
 the following
associativity condition holds:
$$
\CA(\psi\circ\phi)=\CA(\phi)\circ \phi^{-1}(\CA(\psi)),
$$
where $\phi^{-1}$ is understood as the inverse
image functor on the category of sheaves of vector spaces.

It should  be
clear what a natural sheaf morphism $\CA_{X}\rightarrow\CB_{X}$ means.

Here are some obvious examples: $\CO_{X}$, $\Omega^{1}_{X}$,
$\CT_{X}$, and sheaves obtained as a result of all sorts of tensor
operation performed on these. Note that all these sheaves are
sheaves of $\CO_{X}$-modules, and our definition ignores this
extra structure. However, one   talks about natural sheaves
$\CA_{X}$ of different classes of algebras, such as commutative,
associative, Lie, vertex, etc., by requiring that $\CA_{\phi}$
preserve this structure.

Constructed in [MSV] for any smooth algebraic manifold $X$ there is a sheaf
of vertex algebras, $\Omega^{ch}_{X}$. It satisfies the following conditions.

(i) $\Omega^{ch}_{X}$ is natural as a sheaf of vertex algebras and it carries
a bi-grading $\Omega^{ch}_{X}=\oplus_{m,n}\Omega^{ch,m}_{X, n}$
  such that each homogeneous component $\Omega^{ch,m}_{X, n}$ is a natural
sheaf of vector spaces.

(ii) There are natural morphisms:
$$
\Omega^{*}_{X}\hookrightarrow \Omega^{ch}_{X}\hookleftarrow
\Lambda^{*}\CT_{X};
\eqno{(2.1.1)}
$$

(iii)  $\Omega^{ch}_{X}$ is not a sheaf of $\CO_{X}$-modules, but it carries a filtration
such that there is the following family of   natural sheaf isomorphisms
 $$
 \text{Gr}\Omega^{ch}_{X}\iso \bigotimes_{n\geq 0} \left(S^{*}_{q^{n+1}}(\CT_{X})\otimes
S^{*}_{q^{n}}(\CT_{X}^{*})\otimes
\Lambda^{*}_{q^{n+1}y}(\CT_{X})\otimes\Lambda^{*}_{q^{n}y^{-1}}
(\CT_{X}^{*})\right).
\eqno{(2.1.2)}
$$
where we habitually
 use the following `` generating functions of families of sheaves'':
$$
\text{Gr}\Omega^{ch}_{X}=\bigoplus_{m,n}
q^{n}y^{m}\text{Gr}\Omega^{ch,m}_{X,n},
S^{*}_{t}(\CA)=\bigoplus_{n=0}^{\infty}t^{n}S^{n}(\CA),
\Lambda^{*}_{t}(\CA)=\bigoplus_{n=0}^{\infty}t^{n}\Lambda^{n}(\CA).
$$
(iv) it follows from (i) that for any  $X$ there is a canonical group
embedding
$$
\rho_{X}: \text{Aut}X\rightarrow\text{Aut}\Omega^{ch}_{X}(X).
\eqno{(2.1.3)}
$$
Explicit formulas for the latter appeared  in [MSV, (3.1.6)]
as a result of guesswork and were used
 to define $\Omega^{ch}_{X}$ satisfying (i-iii).

In 2.2 we shall look at some examples that serve as a local model
and are needed later;
 in 2.3 we shall very briefly discuss how these local models
are glued together
 and what effect the gluing has on $N2$-structures and chiral rings.

{\bf 2.2. A local model.}

It is easiest to begin with a local situation in the presence of a coordinate
system.

{\bf 2.2.1.} Let $U$ be a smooth affine manifold with a coordinate
system $\vec{x}$ by which we mean a collection of functions
$x_{1},...,x_{n}\in\CO(U)$, $n=\text{dim }U$, such that the
differential forms $dx_{1},...,dx_{n}$ form a basis of the
space of 1-forms $\CT^{*}(U)$
over $\CO(U)$. A coordinate system determines a collection of
  vector fields
$\partial_{x_{1}},...,\partial_{x_{n}}$ such that
$$
\partial_{x_{i}}x_{j}=<\partial_{x_{i}},dx_{j}>=\delta_{ij}.
$$
It follows that $[\partial_{x_{i}},\partial_{x_{j}}]=0$ for all
$i,j$, and $\partial_{x_{1}},...,\partial_{x_{n}}$  form a basis of the
space of vector fields $\CT^{*}(U)$.

 Let
$\Omega^{ch}(U,\vec{x})$ be the following superpolynomial ring
over $\CO(U)$.
$$
\Omega^{ch}(U,\vec{x})=\CO(U)[x_{i,(-j-1)},\;\partial_{x_{i},(-j)};\;
dx_{i,(-j)}\partial_{dx_{i},(-j)},\;1\leq i\leq n,\, j\geq 1],
\eqno{(2.2.1)}
$$
the generators $x_{i,(-j-1)},\partial_{x_{i},(-j)}$ being even,
$dx_{i,(-j)}\partial_{dx_{i},(-j)}$ odd. Identifying $x_{i,(-j-1)}$, $\partial_{x_{i},(-j)}$ with different even
copies of $dx_{i}$ and $\partial_{x_{i}}$ resp., and
$dx_{i,(-j)}$, $\partial_{dx_{i},(-j)}$ with different odd copies
thereof, we
obtain an identification of superalgebras
$$
\Omega^{ch}(U,\vec{x})\iso \bigotimes_{n\geq 0}(S^{*}(\CT(U))\otimes
S^{*}(\CT^{*}(U))\otimes\Lambda^{*}(\CT(U))\otimes\Lambda^{*}(\CT^{*}(U))).
\eqno{(2.2.2)}
$$
This is a local version of (2.1.2), and the images of embeddings
(2.1.1) are generated over $\CO(U)$ by $dx_{i,(-1)}$,
$\partial_{dx_{i}, (-1)}$.

The ring $\Omega^{ch}(U,\vec{x})$  carries
a canonical vertex algebra structure [MSV]. To formulate the result introduce an
 even derivation $T\in\text{End}\Omega^{ch}(U,\vec{x})$
determined by the conditions
$$
T(f)=\sum_{i}x_{i,(-2)}\partial_{x_{i}}f,\;
T(a_{(-n)})=na_{(-n-1)},\eqno{(2.2.3)}
$$
where $f\in\CO(U)$, $a=x_{i}$, $\partial_{x_{i}}$, or
$\partial_{dx_{i}}$, and $n\geq 1$. Note that under identification (2.2.2)
the first of these conditions says that
$T(\CO(U))\subset\CT^{*}(U)$ and the restriction $T|_{\CO(U)}$ equals
the de Rham differential.

{\bf 2.2.2. Lemma.} {\it There is a unique vertex algebra structure on
$\Omega^{ch}(U,\vec{x})$}
$$
Y:\;\Omega^{ch}(U,\vec{x})\rightarrow\text{Field}(\Omega^{ch}(U,\vec{x}))
$$
{\it determined by the conditions:

(i) $\Omega^{ch}(U,\vec{x})$  is generated by the fields
$Y(f,z)$, $f\in\CO(U)$, $Y(\partial_{x_{i},(-1)},z)$, \newline
$Y(\partial_{dx_{i},(-1)},z)$, $Y(dx_{j,(-1)},w)$, the list of
non-zero brackets amongst them being as follows:
$$
[Y(\partial_{x_{i},(-1)},z),Y(f,z)]=\delta(z-w)Y(\partial_{x_{i}}f,w),
\eqno{(2.2.4a)}
$$
$$
[Y(\partial_{dx_{i},(-1)},z),Y(dx_{j,(-1)},w)]=\delta_{ij}\delta(z-w);
\eqno{(2.2.4b)}
$$
(ii)  $T$-covariance:
$$
\aligned
 & \; [T,a(z)]=a(z)',\\
 & a(z)= Y(f,z), Y(\partial_{x_{i},(-1)},z),
Y(\partial_{dx_{i},(-1)},z), Y(dx_{j,(-2)},w); \endaligned
\eqno{(2.2.5)}
$$
(iii) vacuum:}
$$
 \; Y(1,z)=\text{Id},\; Y(f,z)g|_{z=0}=fg, Y(a,z)1|_{z=0}=a,\eqno{(2.2.6)}
$$
{\it where $1,f,g\in\CO(U)$, $a= \partial_{x_{i},(-1)}$, $dx_{i,(-1)}$
or $\partial_{dx_{i}, (-1)}$.}
\bigskip

The uniqueness assertion of this lemma is an immediate consequence
 of Theorem 4.5 in [K]. While proving the existence
assertion in general is something of a problem, in many examples,
sufficient for our present purposes, this is easy. Before we begin
discussing these examples, let us unburden the notation by setting:
$$
\aligned
f(z)=Y(f,z),\;
&dx_{i}(z)=Y(dx_{i,(-1)},z),\\
\partial_{x_{i}}(z)=Y(\partial_{x_{i},(-1)},z),\;
&\partial_{dx_{i}}(z)=Y(\partial_{dx_{i},(-1)},z).
\endaligned
\eqno{(2.2.7)}
$$

\bigskip

{\bf 2.2.3.} {\it Example: an affine space.} If $U=\BC$ with the
canonical coordinate $x=\vec{x}$, then $\Omega^{ch}(\BC,x)$ is
nothing but $FB$ of (1.7.2). Indeed, since
in this case $\CO(U)=\BC[x]$,  $\Omega^{ch}(\BC,x)$ is generated
by the fields $x(z)$, $\partial_{x}(z)$, $dx(z)$, $\partial_{dx}(z)$,
which according to (2.2.4a,b) satisfy
$$
[\partial_{x}(z), x(w)]=\delta(z-w),\;
[\partial_{dx}(z), dx(w)]=\delta(z-w).
$$
A quick glance at (1.6.3, 1.7.1) shows that
$$
b(z)\mapsto \partial_{dx}(z), c(z)\mapsto dx(z),
\gamma(z)\mapsto x(z), \beta(z)\mapsto\partial_{x}(z)
$$
identifies $\Omega^{ch}(\BC,x)$ with $FB$.

Likewise,
$$
\Omega^{ch}(\BC^{N},\vec{x})=FB^{\otimes n}.
\eqno{(2.2.8)}
$$
Incidentally, the same formulas define a vertex algebra morphism
$$
FB^{\otimes N}\hookrightarrow \Omega^{ch}(U,\vec{x}).
\eqno{(2.2.9)}
$$
The nature of this morphism is this: a coordinate system $\vec{x}$
determines an \'etale map  $U\rightarrow\BC^{N}$;
hence (2.2.9) is  a manifestation of the naturality
of $\Omega^{ch}_{X}$, see 2.1.

{\bf 2.2.4.} {\it Example: localization of an affine space.} Let
$f\in\BC[x_{1},...,x_{n}]$, $U_{f}=\BC^{N}-
\{\vec{x}:f(\vec{x})=0\}$, and $\BC[x_{1},...,x_{n}]_{f}$ the
corresponding localization. To extend the vertex algebra structure
from $\Omega^{ch}(\BC^{N},\vec{x})$ to
$$
\Omega^{ch}(U_{f},\vec{x})=\BC[x_{1},...,x_{n}]_{f}\otimes_{\BC[x_{1},...,x_{n}]}
\Omega^{ch}(\BC^{N},\vec{x})
$$
it suffices to define the field $f^{-1}(z)$. In [MSV] an explicit
formula for this field was written down using Feigin's insight.
Lemma 2.2.2 is a convenient alternative tool to compute the action of this
(and similar) fields. Indeed, in view of the commutation relations
(2.2.4a-b) it suffices to know $f^{-1}(z)_{(n)}g$, $g\in\CO(U)$.
Due to (2.2.6) we have
$$
f^{-1}(z)_{(n)}g=\left\{\aligned 0&\text{ if }n\geq 0\\
\frac{g}{f}&\text{ if }n=-1.\endaligned\right.
$$
The values $f^{-1}(z)_{(n)}g$, $n\leq -2$, are determined by using
(2.2.5). The case where $g=1$ suffices and the repeated application
of (2.2.5) gives
$$
f^{-1}(z)1=e^{zT}f^{-1}.
$$
For example,
$$
f^{-1}(z)_{(-2)}1=T(\frac{1}{f})=
-\sum_{i}x_{i,(-2)}\frac{\partial_{x_{i}}f}{f^{2}}.
$$

We shall mostly need localization to the complements of hyperplanes. The
corresponding vertex algebras can be realized, thanks to [B], inside
lattice vertex algebras; this will be reviewed in
some detail in sect. 3.

\bigskip
{\bf 2.2.5.} {\it Two $N2$-structures and two chiral rings.} We
shall need two morphisms of vertex algebras
$$
\rho_{1},\rho_{2}:\;V(N2)_{3n}\rightarrow
\Omega^{ch}(U,\vec{x}),\; n=\text{dim}U.
$$
The first was used in [MSV] and in terms of fields is defined by
$$
\aligned
Q^{(1)}(z)&=\sum_{i}dx_{i}(z)\partial_{x_{i}}(z),\;
G^{(1)}(z)=\sum_{i}:x_{i}(z)'\partial_{dx_{i}}(z),\\
J^{(1)}(z)&=-\sum_{i}:dx_{i}(z)\partial_{dx_{i}}(z):,
L^{(1)}(z)=\sum_{i}:x_{i}(z)'\partial_{x_{i}}(z):+
:dx_{i}(z)'\partial_{dx_{i}}(z):,\endaligned
\eqno{(2.2.10)}
$$
where we let $A^{(1)}=\rho_{1}(A)$, $A$=$Q$, $G$, $J$, or $L$.

The second is obtained by composing the first with
automorphism (1.11.3); the result is this:
$$
\aligned
Q^{(2)}(z)&=\sum_{i}:x_{i}(z)'\partial_{dx_{i}}(z),\;
G^{(2)}(z)=\sum_{i}dx_{i}(z)\partial_{x_{i}}(z),\\
J^{(2)}(z)&=\sum_{i}:dx_{i}(z)\partial_{dx_{i}}(z):,
L^{(2)}(z)=\sum_{i}:x_{i}(z)'\partial_{x_{i}}(z):-
:dx_{i}(z)\partial_{dx_{i}}(z)':.\endaligned \eqno{(2.2.11)}
$$
As was noted in 1.10.1, the operators $L^{(i)}_{(1)}$ give two conformal
gradings and a
 simple computation shows that the corresponding chiral rings, 1.4, are as
follows:
$$
\aligned
\text{Ker}L^{(1)}_{(1)} = \CO(U)[dx_{1,(-1)},...,dx_{n,(-1)}],
\text{Ker}L^{(2)}_{(1)} = \CO(U)[\partial_{dx_{1},(-1)},...,\partial_{dx_{n},(-1)}].
\endaligned
\eqno{(2.2.12)}
$$

\bigskip

{\bf 2.3.} {\it Gluing the local models.}

{\bf 2.3.1.}
Localisation procedure explained in 2.2.4 carries over to any
$\Omega^{ch}(U,\vec{x})$, see
2.2.1, and defines,
in the presence of a coordinate system, a sheaf of vertex algebras
$$
U\supset V\mapsto \Omega^{ch}_{U,\vec{x}}(V)\buildrel \text{def}\over =
\Omega^{ch}(V,\vec{x})
$$
over $U$. By using the action
of the group of coordinate changes [MSV,  (3.1.6)] one obtains canonical identifications
$$
\Omega^{ch}_{U,\vec{x}}\iso \Omega^{ch}_{U,\vec{y}}
$$
for any two coordinate systems $\vec{x}$, $\vec{y}$. This defines
a family of sheaves $U\mapsto \Omega^{ch}_{U}$, where $U$ is \'etale over
$\BC^{N}$, natural w.r.t. to \'etale morphisms.

Finally covering a smooth manifold $X$ by charts $\{U_{\alpha}\}$
\'etale over
$\BC^{N}$ one defines a sheaf $\Omega^{ch}_{X}$ by gluing over intersections
according to the diagram
$$
\Omega^{ch}_{U_{\alpha}}\hookrightarrow
\Omega^{ch}_{U_{\alpha}\cap U_{\beta}}\hookleftarrow
\Omega^{ch}_{U_{\beta}}.
$$
Let us recall, briefly but in some more detail, the effect of this procedure
on the $N2$-structure.

{\bf 2.3.2.} {\it Two $N2$-structures,  the chiral de Rham complex
and algebra of chiral polyvector fields.} It was computed in [MSV]
that a coordinate change $\vec{x}\mapsto\vec{y}$, via (2.1.3),
induces the following transformation of  fields
(2.2.10):
$$
\aligned
Q^{(1)}(z)&\mapsto Q^{(1)}(z)+(d_{DR}(\text{Tr}\log\{(\partial_{x_{i}}y_{j})\}))(z)',\\
G^{(1)}(z)&\mapsto G^{(1)}(z),\\
J^{(1)}(z)&\mapsto J^{(1)}(z)
+(\text{Tr}\log\{(\partial_{x_{i}}y_{j})\})(z)',\\
L^{(1)}(z)&\mapsto L^{(1)}(z),\\
\endaligned
\eqno{(2.3.1)}
$$
and of course similar transformation formulas can be written for
fields (2.2.11). It follows that
$$
L^{(i)}(z)_{(1)}, J^{(i)}(z)_{(0)},\; i=1,2,
$$
is a well-defined quadruple of operators acting on
$\Omega^{ch}_{X}$. Since $[L^{(i)}(z)_{(1)}, J^{(i)}(z)_{(0)}]=0$,
there arise two competing bi-gradings by ``conformal weight,
fermionic charge'':
$$
\aligned &\Omega^{ch}_{X}=\oplus_{n\geq
0,m\in\BZ}\, ^{(i)}\!\Omega^{ch,m}_{X,n},\; i=1,2,\\
&^{(i)}\!\Omega^{ch,m}_{X,n}=\text{Ker}(L^{(i)}(z)_{(1)}-n\text{Id})\cap
\text{Ker}(J^{(i)}(z)_{(0)}-m\text{Id}).
\endaligned\eqno{(2.3.2)}
$$

Formula (2.2.12) shows that
the chiral ring, 1.4, now technically a sheaf of chiral rings,
associated to the first is the  algebra of differential forms:
$$
\CC^{(1)}=\Omega^{*}_{X}:\; U\mapsto \CO(U)[dx_{1,(-1)},...,dx_{n,(-1)}];
\eqno{(2.3.3a)}
$$
the chiral ring associated to
the second is the algebra of polyvector fields
$$
\CC^{(2)}=\Lambda^{*}\CT_{X},\; U\mapsto \CO(U)[\partial_{dx_{1},(-1)},...,\partial_{dx_{n},(-1)}].
\eqno{(2.3.3b)}
$$
This is how morphisms (2.1.1) come about.

{\bf 2.3.3.} {\it Terminology.} From now on we shall call
$\Omega^{ch}_{X}$ equipped with grading (2.3.2) where $i=2$ the
algebra of chiral polyvector fields and re-denote it by
$\Lambda^{ch}\CT_{X}$. The sheaf $\Omega^{ch}_{X}$ equipped with
grading (2.3.2) where $i=1$ will retain the name of the chiral de
Rham complex.

Assertions (2.3.3a,b) are one justification of this terminology.
Note that (2.1.2) uses bi-grading (2.3.2,i=1); the i=2 analogue is
as follows:
 $$
 \text{Gr}\Lambda^{ch}\CT_{X}\iso
\bigotimes_{n\geq 0}\left( S^{*}_{q^{n}}(\CT_{X})\otimes
S^{*}_{q^{n+1}}(\CT_{X}^{*})\otimes\Lambda^{*}_{q^{n}y^{-1}}
(\CT_{X})\otimes\Lambda^{*}_{q^{n+1}y}
(\CT_{X}^{*})\right).
\eqno{(2.3.5)}
$$
Transformation formulas (2.3.1) imply that $L^{(1)}(z)$ is
preserved; hence $\Omega^{ch}_{X}$ always carries a conformal
structure, see 1.10.1 for the definition. The situation is
different with $\Lambda^{ch}\CT_{X}$: it does not carry a
conformal structure compatible with its conformal grading unless
$X$ is Calabi-Yau. Indeed, as follows from the last of formulas
(2.2.11), $L^{(2)}(z)=L^{(1)}(z)+J^{(1)}(z)'$, and the latter
picks the 1st Chern class as a result of transformation (2.3.1).

 If, however, $X$ is a projective Calabi-Yau
manifold, then it can be derived from (2.3.1), [MSV], that
the quadruple of fields $Q^{(i)}(z)$, $G^{(i)}(z)$, $J^{(i)}(z)$,
$Q^{(i)}(z)$, $i=1,2$, can be made sense of globally, and both
$\Omega^{ch}_{X}$, $\Lambda^{ch}\CT_{X}$ acquire an
$N2$-structure, see 1.10.1 for the definition.
 What is especially clear is

{\bf 2.3.4. Lemma.} {\it If $\omega$ is a non-vanishing
holomorphic form over $X$, and $X$ admits an atlas consisting of
charts $\{(U,\vec{x})\}$ such that locally
$\omega=dx_{1}\wedge\cdots\wedge dx_{n}$, then formulas (2.2.10,11)
define an $N2$-structure on $\Omega^{ch}_{X}$ and
$\Lambda^{ch}\CT_{X}$ resp. }

\bigskip

Indeed, in this case the jacobian
$\text{det}(\partial_{x_{i}}y_{j})$ equals 1, and the correction
terms in (2.3.1) vanish.

{\bf 2.3.5.} {\it An example: $X=\BC^{N}-0$.}
As an illsutration, let us compute the cohomology vertex algebra
$H^{*}(\BC^{N}-0,\Lambda^{ch}\CT_{\BC^{N}-0})$, an example that will
prove important later on.

The manifold $\BC^{N}-0$ is quasiaffine and, therefore, it has the
standard global coordinate system $x_{i},\partial_{x_{i}}$, $0\leq
i\leq N-1$, inherited from $\BC^{N}$. This places us in the
situation of Lemma 2.2.2, and we obtain a morphism of bi-graded
sheaves
 $$
 \Omega^{ch}_{\BC^{N}-0}\iso
\bigoplus_{n\geq 0}(S^{*}_{q^{n}}(\CT_{\BC^{N}-0})\otimes
S^{*}_{q^{n+1}}(\CT_{\BC^{N}-0}^{*})\otimes\Lambda^{*}_{q^{n}y^{-1}}
(\CT_{\BC^{N}-0})\otimes\Lambda^{*}_{q^{n}y}
(\CT_{\BC^{N}-0}^{*})),
\eqno{(2.3.6)}
$$
cf. (2.3.5). For the same reason, the sheaf on the R.H.S. of
(2.3.6) is free, hence it suffices to compute
$H^{*}(\BC^{N}-0,\CO_{\BC^{N}-0})$. It is a pleasing excersise in
\v Cech cohomology to prove that
$$
H^{n}(\BC^{N}-0,\CO_{\BC^{N}-0})=
\left\{\aligned
\BC[x_{0},...,x_{N-1}]&\text{ if }n=0\\
\bigotimes_{i=0}^{N-1} \BC[x_{i}^{\pm 1}]/\BC[x_{i}]&\text{ if }n=N-1\\
0&\text{ otherwise}.
\endaligned
\right.
\eqno{(2.3.7)}
$$
The first line of (2.3.7) says that all the global sections of
$\Lambda^{ch}\CT_{\BC^{N}-0}$ are  restrictions from
$\BC^{N}$. Therefore,
$$
H^{0}(\BC^{N}-0,\Lambda^{ch}\CT_{\BC^{N}-0})= \Lambda^{ch}\CT(\BC^{N})=
FB^{\otimes N},
\eqno{(2.3.8)}
$$
cf. (2.2.8).

Similarly, it follows from the 2nd line of (2.3.7) that
$$
H^{N-1}(\BC^{N}-0,\Lambda^{ch}\CT_{\BC^{N}-0})
= \left(\bigotimes_{i=0}^{N-1} \BC[x_{i}^{\pm 1}]/\BC[x_{i}]\right)\otimes_{\BC[\vec{x}]}
\Lambda^{ch}\CT(\BC^{N}),
\eqno{(2.3.9)}
$$
where we use the notation of (2.2.1) with $(U,\vec{x})=(\BC^{N},\vec{x})$.

A moment's thought shows that, as an
$\Lambda^{ch}\CT(\BC^{N})$-module,
$H^{N-1}(\BC^{N}-0,\Lambda^{ch}\CT_{\BC^{N}-0})$ is obtained from
 $\Lambda^{ch}\CT(\BC^{N})$ by  spectral flow (1.12.2):
$$
H^{N-1}(\BC^{N}-0,\Lambda^{ch}\CT_{\BC^{N}-0})=
S_{1}(\Lambda^{ch}\CT(\BC^{N}))=
S_{1}(FB^{\otimes N}).
\eqno{(2.3.10)}
$$
Indeed, by definition 1.7,  $\Lambda^{ch}\CT(\BC^{N})=
FB^{\otimes N}$ is generated by a vector $\b1$ annihilated by
$x_{i,(j)}$, $\partial_{x_{i},(j)}$, $j\geq 0$, (and we identify $x_{i}=x_{i,(-1)}$);
according to (2.3.9), $H^{N-1}(\BC^{N}-0,\Lambda^{ch}\CT_{\BC^{N}-0})$
is generated by a vector annihilated by
$x_{i,(j-1)}$, $\partial_{x_{i},(j+1)}$, $j\geq 0$; the latter annihilating
subalgebra is mapped onto the former by $S_{1}$ of (1.12.1). The odd variables
are treated similarly; but notice also that the Clifford algebra
has only one irreducible module and so the spectral flow on it is inessential.

Of course, the 3rd line of (2.3.7) implies
$$
H^{i}(\BC^{N}-0,\Lambda^{ch}\CT_{\BC^{N}-0})= 0
\text{ if } i\neq 0, N-1.
\eqno{(2.3.11)}
$$

\bigskip

{\bf 2.4. The algebra of chiral polyvector fields over
hypersurfaces.}

This is an exposition of a result of [B].

{\bf 2.4.1.} Let
$$
\CL\rightarrow X\eqno{(2.4.1)}
$$
be a line bundle,
$$
\CL^{*}\rightarrow X\eqno{(2.4.2)}
$$
its dual,
$$
t: X\rightarrow\CL\eqno{(2.4.3)}
$$
its section with smooth zero locus $Z(t)$. Following [B] we shall relate
$\Lambda^{ch}\CT_{\CL^{*}}$ and $\Lambda^{ch}\CT_{Z(t)}$ as
follows.

Identify $t$ with a fiberwise linear function on $\CL^{*}$. We have the
de Rham differential of $t$,
$dt\in H^{0}(\CL^{*},\Omega^{1}_{\CL^{*}})$; via (2.1.1),
$dt\in H^{0}(\CL^{*},\Lambda^{ch}\CT_{\CL^{*}})$. It is clear
that
$$
dt_{(0)}: \Lambda^{ch}\CT_{\CL^{*}}\rightarrow
\Lambda^{ch}\CT_{\CL^{*}} \eqno{(2.4.4)}
$$
is a derivation with zero square, cf. 1.3.

$\CL^{*}$ carries an action of $\BC^{*}$ defined by fiberwise
multiplication.  By the naturality, 2.1 (i),
this action lifts to an action on
$\Lambda^{ch}\CT_{\CL^{*}}$. Hence there  arises  the grading
$$
\Lambda^{ch}\CT_{\CL^{*}}=\bigoplus_{n\in\BZ}R^{n}(\Lambda^{ch}\CT_{\CL^{*}}).
$$
The operator $dt_{(0)}$ has degree 1 with respect to this
grading and, therefore, (2.4.4) is actually a complex of sheaves
$R^{*}\Lambda^{ch}\CT_{\CL^{*}}=
(\Lambda^{ch}\CT_{\CL^{*}},dt_{(0)})$ such that
$$
\cdots\buildrel dt_{(0)}\over\longrightarrow
R^{-1}(\Lambda^{ch}\CT_{\CL^{*}})\buildrel dt_{(0)}\over\longrightarrow
R^{0}(\Lambda^{ch}\CT_{\CL^{*}})\buildrel dt_{(0)}\over\longrightarrow
R^{1}(\Lambda^{ch}\CT_{\CL^{*}})\buildrel dt_{(0)}\over\longrightarrow\cdots
\eqno{(2.4.5)}
$$

{\bf 2.4.2. Lemma} ([B]) {\it The cohomology sheaf
$\CH^{n}_{dt_{(0)}}(\Lambda^{ch}\CT_{\CL^{*}})$ of complex
(2.4.5) is zero unless $n=0$.
The sheaf $\CH^{0}_{dt_{(0)}}(\Lambda^{ch}\CT_{\CL^{*}})$  is supported on $Z(t)$ and naturally
isomorphic to $\Lambda^{ch}\CT_{Z(t)}$.}

{\bf 2.4.3.} Observe that
$dt\in H^{0}(\CL^{*},\Lambda^{ch}\CT_{\CL^{*}})$ is of
conformal weight 1, as follows e.g. from (2.3.5). Therefore,
the differential of complex (2.4.5) preserves
conformal weight, and
the conformal weight 0 component of (2.4.5) is the following classical
complex:
$$
\cdots\buildrel dt\over\longrightarrow
 \Lambda^{i+1}\CT_{\CL^{*}}\buildrel dt\over\longrightarrow
 \Lambda^{i}\CT_{\CL^{*}}\buildrel dt\over\longrightarrow
 \Lambda^{i-1}\CT_{\CL^{*}}\buildrel dt\over\longrightarrow
 \cdots,
\eqno{(2.4.6)}
$$
with differential equal to the contraction with the 1-form $d_{DR}t$.
Lemma 2.4.2 says, in particular, that this complex computes the
algebra of polyvector fields on $Z(t)$, a well-known result
perhaps. Note that in the chiral de Rham complex setting, cf. [B],
this classical construction  is somewhat harder to discern because there
 conformal weight is not preserved by $dt_{(0)}$.

{\bf 2.4.4.} {\it The $N2$-structure.} Let $\CL^{*}$ be the canonical line bundle.
Then both $\CL^{*}$ and $Z(t)$ are Calabi-Yau -- both have a  nowhere zero global
holomorphic volume form --
 and both $\Lambda^{ch}\CT_{\CL^{*}}$ and $\Lambda^{ch}\CT_{Z(t)}$ carry
an $N2$-structure, see the end of 2.3.3. Let us write down some explicit formulas.

Suppose there is a nowhere zero global holomorphic volume form $\omega$
and $X$ can be covered by charts
$s,y_{1},...,y_{N-1}$, $s$ being the coordinate along  the fiber,
such that $\omega=ds\wedge y_{1}\wedge\cdots\wedge y_{N-1}$. Then,
as follows from Lemma 2.3.4,
$$
\aligned
Q(z)\mapsto s(z)'\partial_{ds}(z)+\sum_{j=1}^{N-1}y_{j}(z)'\partial_{dy_{j}}(z),\;
&G(z)\mapsto ds(z)\partial_{s}(z)+\sum_{j=1}^{N-1}dy_{j}(z)\partial_{y_{j}}(z),\\
J(z)\mapsto -ds(z)\partial_{ds}(z)-\sum_{j=1}^{N-1}dy_{j}(z)\partial_{dy_{j}}(z),\;
&L(z)\mapsto s(z)'\partial_{s}(z)+\sum_{j=1}^{N-1}y_{j}(z)'\partial_{y_{j}}(z)-\\
&-ds(z)\partial_{ds}(z)'-\sum_{j=1}^{N-1}dy_{j}(z)\partial_{dy_{j}}(z)'
\endaligned
\eqno{(2.4.7)}
$$
 defines an $N2$-structure on  $\Lambda^{ch}\CT_{\CL^{*}}$.

[B, Proposition 5.8] says that, via Lemma 2.4.2, the $N2$-structure on
$\Lambda^{ch}\CT_{Z(t)}$ is determined by
$$
\aligned
G(z)&\mapsto ds(z)\partial_{s}(z)+\sum_{j=1}^{N-1}dy_{j}(z)\partial_{y_{j}}(z),\\
Q(z)&\mapsto s(z)'\partial_{ds}(z)+\sum_{j=1}^{N-1}y_{j}(z)'\partial_{dy_{j}}(z)-
(s(z)\partial_{ds}(z))'.
\endaligned
\eqno{(2.4.8)}
$$

{\bf 3. The lattice vertex algebra realization and applications to toric
varieties.}
\bigskip

This section is an exposition of part of Borisov's free field realization
[B]. It does not contain any new results except perhaps Lemma 3.8, and
in order to construct the spectral sequence
appearing in the latter the entire section had to be written up.

{\bf 3.1.} Let $M$ be a rank $N$ free abelian group,
$M^{*}=Hom_{\BZ}(M,\BZ)$ its dual. Give $\Lambda=M\oplus M^{*}$
a lattice structure by defining the symmetric bilinear form
$$
\Lambda\times\Lambda\rightarrow \BZ,
\eqno{(3.1.1)}
$$
induced by the natural pairing $M^{*}\times M\rightarrow \BZ$,
$(X^{*},X)\mapsto X^{*}(X)$.
There arises the lattice vertex algebra $V_{\Lambda}$, 1.9, where
we fix the following cocycle, cf.(1.9.1),
$$
\epsilon(X+X^{*},Y+Y^{*})=(-1)^{X^{*}(Y)},\; X,Y\in M, X^{*},Y^{*}\in M^{*}.
\eqno{(3.1.2)}
$$
Next, consider the complexification $\Lambda_{\BC}=\BC\otimes_{\BZ}\Lambda$
onto which form (3.1.2) carries over. There arises the fermionic vertex
algebra $F_{\Lambda_{\BC}}$ which we re-denote by  $F_{\Lambda}$, see (1.6.5).

Finally, following [B] introduce Borisov's vertex algebra
$$
\BB_{\Lambda}=V_{\Lambda}\otimes F_{\Lambda}.
\eqno{(3.1.3)}
$$

{\bf Notation.} The notational problem one faces here is that the
lattice $\Lambda$  twice manifests itself inside
$\BB_{\Lambda}$: first, as an ingredient of $V_{\Lambda}$; second, as that
of $ F_{\Lambda}$. We attempt to resolve this issue by letting
capital latin letters, $X,Y,Z,..$
( $X^{*},Y^{*},Z^{*},..$ resp.), denote elements of $M$
($M^{*}$ resp.) in the context of $V_{\Lambda}$; and let the tilded
letters,  $\tilde{X},\tilde{Y},\tilde{Z},..$ or
 $\tilde{X}^{*},\tilde{Y}^{*},\tilde{Z}^{*},..$ denote their respective copies
in the context of   $ F_{\Lambda}$. Later on this will be related to geometry
and then we shall let the lowercase letters denote the respective
coordinates. $\qed$

Note that the assignment $M^{*}\mapsto \BB_{\Lambda}$ is functorial. Indeed,
if $g\in\text{Hom}(M^{*}_{1},M^{*}_{2})$ is an isomorphism of abelian groups,
then
$$
(g^{-1}, g^{*})\in \text{Hom}(M^{*}_{2},M^{*}_{1})\times
\text{Hom}(M_{2},M_{1})\hookrightarrow\text{Hom}(\Lambda_{2},\Lambda_{1})
$$
is an isomorphism of lattices preserving form (3.1.1) and cocycle (3.1.2).
According to (1.11.1,2), this isomorphism induces the following isomorphism
of vertex algebras
$$
\hat{g}: \BB_{\Lambda_{2}}\rightarrow \BB_{\Lambda_{1}}
\eqno{(3.1.4a)}
$$
$$
\aligned
&X^{*}(z)\mapsto g^{-1}X^{*}(z),\tilde{X}^{*}(z)\mapsto g^{-1}\tilde{X}^{*}(z),
e^{X^{*}}(z)\mapsto e^{g^{-1}X^{*}}(z)\\
&X(z)\mapsto g^{*}X(z),\tilde{X}(z)\mapsto
 g^{*}\tilde{X}(z),
e^{X}(z)\mapsto e^{ g^{*}X}(z).
\endaligned
\eqno{(3.1.4b)}
$$

Therefore, if we introduce the category of lattices $\Lambda$,
morphisms being the described isomorphisms, then
$$
M^{*}\mapsto\BB_{\Lambda},
g\mapsto\hat{g} \eqno{(3.1.5)}
$$
is a contravariant functor.

Later we shall have to work with $g$ such that
 $g(M^{*}_{1})\subset M^{*}_{2}$ but after the extension of scalars
to $\BQ$ the induced
$g\in\text{Hom}_{\BQ}((M^{*}_{1})_{\BQ},(M^{*}_{2})_{\BQ})$ is an
isomorphism. In this case the functorial nature of $M^{*}\mapsto
\BB_{\Lambda}$ is a little more subtle because $g^{-1}$ may have non-integer
entries.  There are two ways around.

Consider a vertex subalgebra $\BB_{M,\Lambda}\subset
\BB_{\Lambda}$ associated to the sublattice $M\subset \Lambda$,
see definition (1.9.10). (This simply means that all the fields
$e^{X^{*}}(z)$ are not allowed.) It is clear that
$$
\BB_{.,.}:\; M^{*}\mapsto\BB_{M,\Lambda},
g\mapsto\hat{g}|_{ \BB_{M,\Lambda}}
\eqno{(3.1.6)}
$$
is a contravariant functor because the indicated
restriction of (3.1.4a) makes sense for any lattice
embedding $g$.

Second, naturally associated to $g$ there is a map
$$
\hat{g}: \BB_{\Lambda_{2}}\hookrightarrow \BB_{g^{-1}\Lambda_{2}},
\eqno{(3.1.7)}
$$
still defined by (3.1.4a), where the lattice $g^{-1}\Lambda_{2}$ is defined
to be $g^{*}M_{2}\oplus g^{-1}M^{*}_{2}\subset (\Lambda_{1})_{\BQ}$.

Note that (3.1.6) is a ``subfunctor'' of (3.1.7).

{\bf 3.2.}
Let us introduce
the following terminology and notation
pertaining to toric variety theory: by a basic cone $\sigma\subset M^{*}$
we shall mean a sub-semigroup spanned over $\BZ_{+}$ by part of a basis
(over $\BZ$) of $M^{*}$. Let $<\sigma>$  denote the (uniquely
determined) spanning set
of $\sigma$ . Given a basic cone
$\sigma\subset M^{*}$, let $\check\sigma\subset M$ be its dual cone
 defined  by $\check\sigma=\{X\in M\text{ s.t. } \sigma (X)\geq 0\}$.

A smooth toric variety will always be defined by fixing a lattice
$\Lambda$ as in 3.1 and a regular fan $\Sigma$.
(Regular means  that $\Sigma$ is a collection of basic
cones in $M^{*}$.) If we define
$$
U_{\sigma}=\text{Spec}\BC[\check\sigma],\; \sigma\in\Sigma,
\eqno{(3.2.1)}
$$
where $\BC[\check\sigma]$ is the semigroup algebra of $\check\sigma$,
then there arises a canonical embedding
$$
U_{\sigma'}\subset U_{\sigma},\; \sigma'\subset\sigma.
\eqno{(3.2.2)}
$$

The toric variety $X_{\Sigma}$
 attached to $\Sigma$ is defined by declaring that
$$
\CU_{\Sigma}=\{U_{\sigma},\sigma\in\Sigma\}
\eqno{(3.2.3)}
$$
is its covering and by gluing the charts over intersections
$$
U_{\alpha}\hookleftarrow U_{\alpha\cap\beta}\hookrightarrow U_{\beta}
$$
according to (3.2.2).

Note that the assignment $(\sigma, M^{*})\mapsto U_{\sigma}$ is
functorial. Indeed, if we introduce the category whose objects are pairs
$(\sigma, M^{*})$ and morphisms
$(\sigma_{1}, M^{*}_{1})\rightarrow (\sigma_{2}, M^{*}_{2})$
 are abelian group morphisms
$g:  M^{*}_{1}\rightarrow  M^{*}_{2}$ such that
 $g(\sigma_{1})\subset\sigma_{2}$,
then $g^{*}(\check\sigma_{2})\subset\check\sigma_{1}$. Hence $g^*$ induces
a ring homomorphism $\BC[\check\sigma_{2}]\rightarrow
\BC[\check\sigma_{1}]$ and thus a morphism $\tilde{g}:
U_{\sigma_{1}}\rightarrow
U_{\sigma_{2}}$. Of course,
$$
(\sigma, M^{*})\mapsto U_{\sigma},\; g\mapsto\tilde{g}
\eqno{(3.2.4)}
$$
is a covariant functor.

This can be globalized: given $(\Sigma_{1}, M_{1}^{*})$and
$(\Sigma_{2}, M_{2}^{*})$ with a lattice morphism
$g:  M_{1}^{*}\rightarrow  M_{2}^{*}$ such that for each
$\sigma_{1}\in \Sigma_{1}$ there is $\sigma_{2}\in \Sigma_{2}$
containing $g(\sigma_{1})$, there arises a morphism
$$
\tilde{g}: X_{\Sigma_{1}}\rightarrow X_{\Sigma_{2}}.
\eqno{(3.2.5)}
$$

{\bf 3.3.}
We would like to define Borisov's realization [B], that is,
 a vertex algebra embedding
$\Lambda^{ch}\CT(U_{\sigma})\hookrightarrow\BB_{M,\Lambda}$
for each basic cone $\sigma\in M^{*}$. To write down an explicit formula
for this map, let us choose a basis of $M^{*}$, $X_{0}^{*},..., X_{N-1}^{*}$,
such that $\sigma$ is spanned by $X_{0}^{*},..., X_{m-1}^{*}$.
 This fixes the dual basis
$X_{0},..., X_{N-1}$ of $M$.

In order to conform
to the notation of sect.2, let
$x_{j}\buildrel\text{def}\over = e^{X_{j}}$, $\partial_{x_{i}}$,
$0\leq i\leq m-1$,
be a coordinate system on $U_{\sigma}$. Thus
$$
U_{\sigma}=\text{Spec}\BC[x_{0},...,x_{m-1},x_{m}^{\pm 1},...,x_{N-1}^{\pm 1}].
$$
Borisov proves that there is a vertex algebra homomorphism
$$
\CB(\sigma):\; \Lambda^{ch}\CT(U_{\sigma})
\rightarrow \BB_{M,\Lambda}
\eqno{(3.3.1)}
$$
determined by the assignment
$$
\aligned & x_{i}^{\pm 1}(z)\mapsto e^{\pm X_{i}}(z),
x_{j}(z)\mapsto e^{X_{j}}(z),\;
m\leq i\leq N-1,j\leq m-1,\\
& dx_{i}(z)\mapsto :e^{X_{i}}(z)\tilde{X}_{i}(z):.
\endaligned
\eqno{(3.3.2a)}
$$
$$
\partial_{x_{i}}(z)\mapsto
:(X_{i}^{*}(z)-:\tilde{X}_{i}(z)\tilde{X}_{i}^{*}(z):)e^{-X_{i}}(z):,\;
\partial_{dx_{i}}(z)\mapsto :e^{-X_{i}}(z)\tilde{X}_{i}^{*}(z):,
\eqno{(3.3.2b)}
$$
cf. sect. 3.1, Notation. (Note that  (3.3.2a-b) are exactly
(2.1.3) specialized to the exponential change of variables
$x_{i}\rightarrow e^{X_{i}}$; this remark is also borrowed from
[B].) Formulas (3.3.2a) are manifestly independent of the choice
of variables; it is easy then to infer that so are (3.3.2b).

This embedding naturally depends on $\sigma$. To make a precise statement,
give

{\bf 3.3.1. Definition.} Fix a number $N$.
 $\CC$ is a category whose objects are pairs
$(\sigma, M^{*})$, $\text{dim}M^{*}=N$, all  and morphisms
$(\sigma_{1}, M^{*}_{1})\rightarrow (\sigma_{2}, M^{*}_{2})$
 are abelian group embeddings
$g:  M^{*}_{1}\hookrightarrow  M^{*}_{2}$ such that
 $g(\sigma_{1})=\sigma_{2}$. $\qed$

\bigskip

Note that the conditions imposed  on morphisms in this definition
strengthen those used in  (3.2.4).
In fact, a short computation
shows that the map $\tilde{g}:U_{\sigma_{1}}\rightarrow U_{\sigma_{2}}$ associated
 with a morphism
$g\in\text{Mor}_{\CC}((\sigma_{1}, M_{1}^{*}),
(\sigma_{1}, M_{1}^{*}))$ in (3.2.4) is \'etale. Hence, by virtue of the naturality of $\Lambda^{ch}\CT_{X}$, see
2.1(i), the composition
$$
(\sigma, M^{*})\mapsto U_{\sigma}\mapsto\Lambda^{ch}\CT(U_{\sigma})
$$
defines a contravariant functor
$$
\Lambda^{ch}\CT(.):\;(\sigma, M^{*})\mapsto\Lambda^{ch}\CT(U_{\sigma}).
\eqno{(3.3.3)}
$$
By forgetting $\sigma$, we regard
contravariant functor  (3.1.6)  as defined on $\CC$.
 The main property of (3.3.1) is formulated as follows.

{\bf 3.4. Lemma.} {\it The assignment $(\sigma,M^{*})\mapsto
\CB(\sigma)$, see (3.3.1), is a functor morphism}
$$
\CB(\sigma):\; \Lambda^{ch}\CT(.)\rightarrow \BB_{.,.},
$$
{\it where $\BB_{.,.}$ is the functor defined in (3.1.6).}

\bigskip

{\bf Notational convention.} Using this fact we shall not distinguish
between $\Lambda^{ch}\CT(U_{\sigma})$ and
$\CB(\sigma)(\Lambda^{ch}\CT(U_{\sigma}))\subset\BB_{M,\Lambda}$.

{\bf 3.5.} Let $S^{*}\in<\sigma>$ be one of the generators of $\sigma$
and let $\sigma\setminus S^{*}$ denote the cone spanned by $<\sigma>\setminus
S^{*}$.
There arises then the restrtiction morphism
$$
res(\sigma, S^{*}): \Lambda^{ch}\CT(U_{\sigma})\hookrightarrow \Lambda^{ch}\CT(U_{\sigma\setminus S^{*}}),
$$
and one would like to extend it to a resolution.

Let
$$
\CJ^{*}(\sigma,S^{*})=\oplus_{n=0}^{\infty}\CJ^{n}(\sigma,S^{*}),\;
\CJ^{n}(\sigma,S^{*})=\Lambda^{ch}\CT(U_{\sigma\setminus S^{*}})e^{nS^{*}},
\eqno{(3.5.1)}
$$
where $\Lambda^{ch}\CT(U_{\sigma\setminus S^{*}})$ is thought of as
a vertex subalgebra of $\BB_{M,\Lambda}$, see 3.4, Notational convention,
and this makes sense out of
$\Lambda^{ch}\CT(U_{\sigma\setminus S^{*}})e^{nS^{*}}$ as a subspace of
 $\BB_{\Lambda}$.

$\CJ^{*}(\sigma,S^{*})$ is evidently a $\BZ_{+}$-graded vertex subalgebra of $\BB_{\Lambda}$.

Let
$$
D(\sigma,S^{*})=(e^{S^{*}}(z)\tilde{S}^{*}(z))_{(0)}.
\eqno{(3.5.2)}
$$
It is evidently a square zero derivation of $\CJ^{*}(\sigma,S^{*})$,
see 1.3. Thus
$(\CJ^{*}(\sigma,S^{*}),D(\sigma,S^{*}))$ is a differential graded
vertex algebra.

Now look upon $\Lambda^{ch}\CT(U_{\sigma})$  as  a differential graded
vertex algebra with  $\Lambda^{ch}\CT(U_{\sigma})$ placed in degree 0, zero
spaces placed everywhere else, and  zero differential.
Then, by taking the composition
$$
\Lambda^{ch}\CT(U_{\sigma})\buildrel res(\sigma, S^{*})\over\rightarrow
\Lambda^{ch}\CT(U_{\sigma\setminus S^{*}})=\CJ^{0}(\sigma,S^{*})
\hookrightarrow \CJ^{*}(\sigma,S^{*}),
$$
$res(\sigma, S^{*})$ can be interpreted as a morphism of differential
graded vertex algebras
$$
res(\sigma, S^{*}): (\Lambda^{ch}\CT(U_{\sigma}),0)\hookrightarrow (\CJ^{*}(\sigma,S^{*}),D(\sigma,S^{*})).
\eqno{(3.5.3)}
$$
This construction is natural. To explain this, let us give the following
definition, cf.  3.3.1.

{\bf 3.5.1. Definition.}  $\CC_{pnt}$ is a category whose objects are triples
$(S^{*},\sigma, M^{*})$ with $S^{*}\in<\sigma>$,
$(\sigma, M^{*})\in\text{Ob}(\CC)$, and morphisms
$(S_{1}^{*},\sigma_{1}, M^{*}_{1})\rightarrow (S_{2}^{*},\sigma_{2}, M^{*}_{2})$
 are abelian group morphisms
$g:  M^{*}_{1}\rightarrow  M^{*}_{2}$ such that
 $g(\sigma_{1})=\sigma_{2}$ and  $g(S^{*}_{1})=S^{*}_{2}$.  $\qed$

It is clear that both
$$
(\Lambda^{ch}\CT(.),0):\; (S^{*},\sigma, M^{*})\mapsto (\Lambda^{ch}\CT(U_{\sigma}),0)
\eqno{(3.5.4)}
$$
and
$$
(\CJ^{*}(.), D(.)):\; (S^{*},\sigma, M^{*})\mapsto (\CJ^{*}(\sigma,S^{*}),D(\sigma,S^{*}))
\eqno{(3.5.5)}
$$
are contravariant functors from $\CC_{pnt}$ to the
category of differential graded vertex algebras. Indeed, the former is essentially (3.3.3), as to the latter, one has to apply map (3.1.7)
restricted to $\CJ^{*}(\sigma,S^{*})$.

{\bf 3.6. Lemma.} {\it (i) The assignment $(S^{*},\sigma)\mapsto
res(\sigma, S^{*})$, see (3.5.3), is a functor morphism
$$
res(.):\; (\Lambda^{ch}\CT(.),0)\mapsto (\CJ^{*}(.),D(.)).
$$

(ii) For each $(S^{*},\sigma, M^{*})$ map (3.5.3) is a quasiisomorphism.}

{\bf 3.7. } Let us apply Lemma 3.6 to the situation where
the fan $\Sigma$ satisfies the following condition:
 there is $S^{*}\in M^{*}$ such that $S^{*}\in <\sigma>$ for all
 highest dimension cones $\sigma\in\Sigma$. Let
$\Sigma\setminus S^{*}=\{\sigma\setminus S^{*}:\; \sigma\in\Sigma\}$,
cf. the beginning of 3.5.
This means that the morphism
$$
X_{\Sigma}\rightarrow X_{\Sigma/\BZ S^{*}}
\eqno{(3.7.1)}
$$
induced by the canonical projection
$M^{*}\rightarrow M^{*}/\BZ S^{*}$
is a line bundle, and the map
$$
 X_{\Sigma\setminus S^{*}}\hookrightarrow X_{\Sigma}
\eqno{(3.7.2)}
$$
induced by the tautological inclusion
$\Sigma\setminus S^{*}\subset\Sigma$ is the embedding of the total space of
the line bundle without the zero section. We would like to relate
the cohomology groups $H^{*}(X_{\Sigma},\Lambda^{ch}\CT_{X_{\Sigma}})$
and $H^{*}(X_{\Sigma\setminus S^{*}},
\Lambda^{ch}\CT_{X_{\Sigma\setminus S^{*}}})$.

The nerve of the covering $\CU_{\Sigma}$
is a simplicial object in the category $\CC_{pnt}$. Applying to it
functor (3.5.4) one gets the complex $\Lambda^{ch}\CT(\CU_{\Sigma})$ commonly
known as the \v Cech complex.
(A complex, not a bi-complex, because we ignore the trivial differential
on $(\Lambda^{ch}\CT(.),0)$.) We  denote it by
$\check{C}^{*}(\CU_{\Sigma},\Lambda^{ch}\CT_{X_{\Sigma}};d_{\check C})$,
where
$d_{\check C}$ is the \v Cech differential.

Likewise  applying  functor
(3.5.5) to the nerve of $\CU_{\Sigma}$,
 we obtain the bi-complex
$\check{C}^{*}(\CU_{\Sigma},\CJ^{*};d_{\check C}, D(S^{*}))$. By definition
$$
\check{C}^{p}(\CU_{\Sigma},\CJ^{q};d_{\check C}, D(S^{*}))=
\check{C}^{p}(\CU_{\Sigma\setminus S^{*}},\Lambda^{ch}
\CT_{X_{\Sigma\setminus S^{*}}}e^{qS^{*}};d_{\check C}, D(S^{*})).
\eqno{(3.7.3)}
$$
According to Lemma 3.6,
$$
res(\CU_{\Sigma}):
\check{C}^{*}(\CU_{\Sigma},\Lambda^{ch}\CT_{X_{\Sigma}};d_{\check C})
\rightarrow
\check{C}^{*}(\CU_{\Sigma},\CJ^{*};d_{\check C}, D(S^{*}))
\eqno{(3.7.4)}
$$
is a quasiisomorphism. More precisely, the
 bi-complex
$\check{C}^{*}(\CU_{\Sigma},\CJ^{*};d_{\check C}, D)$ gives rise
to two spectral sequences both converging to its total cohomology.
The one where the vertex differential $D(.)$ is used first
degenerates in the first term to the \v Cech complex $\check{C}^{*}(\CU_{\Sigma},\Lambda^{ch}\CT_{X_{\Sigma}};d_{\check C})$
-- this follows at once from Lemma 3.6 (ii). Hence both the sequences
abut to  $H^{*}(X_{\Sigma},\Lambda^{ch}\CT_{X_{\Sigma}})$.
By definition,
the first and  the second terms of the second spectral sequence are as follows:
$$
\aligned
&(E^{p,q}_{1}, d_{1})=(H^{p}(X_{\Sigma\setminus S^{*}},\Lambda^{ch}
\CT_{X_{\Sigma\setminus S^{*}}}e^{qS^{*}}), D(S^{*})),\\
&D(S^{*})=(e^{S^{*}}(z)\tilde{S}^{*}(z))_{(0)}.
\endaligned
\eqno{(3.7.5)}
$$
$$
E^{p,q}_{2}=H_{D(S^{*})}^{q}(H^{p}(X_{\Sigma\setminus S^{*}},\oplus_{n=0}^{\infty}\Lambda^{ch}
\CT_{X_{\Sigma\setminus S^{*}}}e^{nS^{*}})).
\eqno{(3.7.6)}
$$

Let us summarize our discussion.

{\bf 3.8. Lemma.} {\it There is a spectral sequence
$\{E^{p,q}_{r},d_{r}\}\Rightarrow
H^{*}(X_{\Sigma},\Lambda^{ch}\CT_{X_{\Sigma}})$ that
satisfies (3.7.5,6).}

{\bf 3.9.} The formation of bi-complex (3.7.3), and hence of the
corresponding spectral sequences, is  functorial in $X_{\Sigma}$.
To make this precise, let a triple $S^{*}$, $\Sigma_{2}$, and
$M^{*}_{2}$ satisfy the conditions imposed on $S^{*}$, $\Sigma$,
and $M^{*}$ in 3.7, and let us give ourselves another pair
$\Sigma_{1}$ and $M^{*}_{1}$,
$\text{dim}M^{*}_{1}=\text{dim}M^{*}_{2}$, along with a lattice
embedding
$$
g:\; M_{1}^{*}\rightarrow M^{*}_{2}\;\text{ s.t. }
g(\Sigma_{1})=\Sigma_{2}\setminus S^{*}. \eqno{(3.9.1)}
$$
According to (3.2.5) this induces a map
$$
\tilde{g}:\; X_{\Sigma_{1}}\mapsto X_{\Sigma_{2}},
\eqno{(3.9.2)}
$$
which is  \'etale, cf. 3.3.
Map (3.9.1) gives rise to the lattice $g^{-1}\Lambda^{*}_{2}$ and
 the embedding of Borisov's algebras
$$
\hat{g}:\; \BB_{\Lambda_{2}}\hookrightarrow \BB_{g^{-1}\Lambda_{2}}
\eqno{(3.9.3)}
$$
due to (3.1.7). It is rather obvious that maps (3.9.2,3) allow to pull
the bi-complex
$$
\check{C}^{p}(\CU_{\Sigma_{2}\setminus S^{*}},\Lambda^{ch}
\CT_{X_{\Sigma_{2}\setminus S^{*}}}e^{qS^{*}};d_{\check C}, D(S^{*}))
\eqno{(3.9.4)}
$$
back onto $X_{\Sigma_{1}}$. Let us write down the relevant
formula. An element of bi-complex (3.9.4)
 is a  family of elements
$$
f_{\sigma}\in\{\Lambda^{ch}\CT_{X_{\Sigma_{2}\setminus
S^{*}}}(U_{\sigma})\}e^{qS^{*}},\; \sigma\in\Sigma_{2}.
$$

Define
$$
\check{C}^{p}(\CU_{\Sigma_{1}},\hat{g}\Lambda^{ch}
\CT_{X_{\Sigma_{2}}}e^{qg^{-1}S^{*}};d_{\check C}, g^{-1}D(S^{*}))
\eqno{(3.9.5)}
$$
 to be the set of all the families
$$
f_{\sigma}\in\hat{g}\{\Lambda^{ch}\CT_{X_{\Sigma_{2}}}(U_{g\sigma})\}e^{qg^{-1}S^{*}},\;
\sigma\in\Sigma_{1}.
$$
By construction, the map
$$
\{f_{\sigma},\; \sigma\in\Sigma_{2}\}\mapsto
\{\hat{g}f_{g\sigma},\; \sigma\in\Sigma_{1}\}
$$
delivers an isomorphism of
 (3.9.5) and (3.9.4):
$$
\aligned
&\check{C}^{p}(\CU_{\Sigma_{2}\setminus S^{*}},\Lambda^{ch}
\CT_{X_{\Sigma_{2}\setminus S^{*}}}e^{qS^{*}};d_{\check C}, D(S^{*}))
\iso\\
&\check{C}^{p}(\CU_{\Sigma_{1}}, \hat{g}\Lambda^{ch}
\CT_{X_{\Sigma_{1}}}e^{qg^{-1}S^{*}};d_{\check C}, g^{-1}D(S^{*})).
\endaligned
\eqno{(3.9.6)}
$$

{\bf 3.10.} {\it Digression: Borisov's realization and the spectral flow.}
We are now making good on our promise to show how spectral flow (1.12.1)
is realized via lattice vertex algebras in the case of the
$bc-\beta\gamma$-system.  We shall show in
5.2.15 that a simple version of this construction
 does the same for $N2$.

According to our conventions
$\Lambda^{ch}\CT(\BC^{N})\subset\BB_{M,\Lambda}$, and  by
definition (3.1.2) $M\subset L$ satisfies (1.12.5). Therefore, any
$\alpha\in M^{*}$ generates on $\text{Lie}\BB_{M,\Lambda}$
 the spectral flow in the direction of $\alpha$,
see (1.12.6). A glance at (3.3.2a,b) shows that if
$\alpha=\sum_{j}X_{j}^{*}$, then (1.2.6) reads
$$
S_{\sum_{j}X_{j}^{*}}:
\aligned
&x_{i}(z)\mapsto x_{i}(z)z,\; dx_{i}(z)\mapsto  dx_{i}(z)z,\\
&\partial_{x_{i}}(z)\mapsto \partial_{x_{i}}(z)z^{-1},
\partial_{dx_{i}}(z)\mapsto \partial_{dx_{i}}(z)z^{-1},
\endaligned
\eqno{(3.10.1)}
$$
cf. 1.12.7-8, and this  does coincide with $S_{1}$ of (1.12.1).

It follows that the map
$$
e^{-\sum_{j}X^{*}_{j}}:\;
e^{\sum_{j}X^{*}_{j}}\Lambda^{ch}\CT(\BC^{N})\rightarrow
\Lambda^{ch}\CT(\BC^{N});
v\mapsto e^{-\sum_{j}X^{*}_{j}}v
\eqno{(3.10.2)}
$$
identifies $e^{\sum_{j}X^{*}_{j}}\Lambda^{ch}\CT(\BC^{N})$ with
the spectral flow transform of
 $\Lambda^{ch}\CT(\BC^{N})$:
$$
e^{\sum_{j}X^{*}_{j}}\Lambda^{ch}\CT(\BC^{N})\iso
S_{1}(\Lambda^{ch}\CT(\BC^{N})),
\eqno{(3.10.3)}
$$
see definition of the spectral flow transform (1.12.2).

Therefore, the results of 2.3 can be rewritten as follows:
$$
H^{n}(\BC^{N}-0,\Lambda^{ch}\CT_{\BC^{N}-0})=\left\{
\aligned
\Lambda^{ch}\CT(\BC^{N})&\text{ if } n=0\\
\Lambda^{ch}\CT(\BC^{N})e^{\sum_{j}X^{*}_{j}}&\text{ if } n=N-1\\
0&\text{ otherwise}.
\endaligned
\right.
\eqno{(3.10.4)}
$$

\newpage

{\bf 4. Chiral polyvector fields over hypersurfaces in projective
spaces}

\bigskip

Having put all the preliminaries out of the way we can tackle our
main problem -- computation of $H^{*}(F,\Lambda^{ch}\CT_{F})$ for
a smooth hypersurface $F\subset \BP^{N-1}$.

\bigskip

{\bf 4.1.} Let
$\CL\rightarrow\BP^{N-1}$  be the  degree $N$ line bundle over
$\BP^{N-1}$ and
$$\pi:\CL^{*}\rightarrow\BP^{N-1}
\eqno{(4.1.0)}
$$
its dual.
 Let us give the spaces
$\BC_{-N}=\BC$ and $\BC^{N}-0$ a $\BC^{*}$-space structure as
follows:
$$
\aligned \BC_{-N}\times\BC^{*}&\rightarrow\BC_{-N},\; u\cdot t=
st^{-N},\\
\BC^{*}\times(\BC^{N}-0)&\rightarrow \BC^{N}-0,\;
t\cdot(x_{0},...,x_{N-1})=(tx_{0},...,tx_{N-1}).
\endaligned
\eqno{(4.1.1)}
$$
One has the quotient realization
$$
\CL^{*}= \BC_{-N}\times_{\BC^{*}}(\BC^{N}-0),\eqno{(4.1.2)}
$$
where we impose the relation
$$
(u;x_{0},...,x_{N-1})\sim (ut^{-N};tx_{0},...,tx_{N-1}),\; t\neq
0. \eqno{(4.1.3)}
$$
Let $\BZ_{N}$ act on $\BC^{N}$ as follows
$$
\aligned &\BZ^{N}\times (\BC^{N}-0)\rightarrow (\BC^{N}-0),\\
&\bar{m}\cdot(x_{0},...,x_{N-1})=(e^{2\pi\sqrt{-1}m/N}x_{0},...,e^{2\pi\sqrt{-1}m/N}x_{N-1}).
\endaligned\eqno{(4.1.4)}
$$
Crucial for our purposes is the following isomorphism of smooth algebraic
varieties
$$
\aligned
&p:\;(\BC^{N}-0)/\BZ_{N}\iso\CL^{*}-0,\\
&\text{class of }(x_{0},...,x_{N-1})\mapsto \text{class of
}(1;x_{0},...,x_{N-1}),\endaligned\eqno{(4.1.5)}
$$
where $\CL^{*}-0$ denotes $\CL^{*}$ without the zero section. (This is
well known and obvious: deleting the zero section
means requiring that $u\neq 0$; then one uses the $\BC^{*}$-action
to make $u=1$; this leaves only the classes of
$(1;x_{0},...,x_{N-1})$ and simultaneously breaks the
$\BC^{*}$-action down to the $\BZ_{N}$-action defined in (4.1.4).)

We wish to study Calabi-Yau hypersurfaces in $\BP^{N-1}$. To
define any such hypersurface, take $x_{0},...,x_{N-1}$ to be, in
accordance with (4.1.1), the homogeneous coordinates on $\BP^{N-1}$.
Let
$$
\fF=\{(x_{0}:\cdots :x_{N-1})\text{ s.t. }f(x_{0},...,x_{N-1})=0\},
\eqno{(4.1.6)}
$$
where $f$ is a degree $N$ homogeneous polynomial with a unique
singularity at 0. Such an $f$ can be regarded as a section of
$\CL$. The corresponding fiberwise linear function $t=t(u;x_{0},...,x_{N-1})$
 on $\CL^{*}$ is
$$
t=uf(x_{0},...,x_{N-1}),\eqno{(4.1.7)}
$$
cf. (4.1.2,3).
The pull-back of this function onto $(\BC^{N}-0)/\BZ_{N}$ under
(4.1.5)  is literally $f(x_{0},...,x_{N-1})$:
$$
p^{*}(t)=f(x_{0},...,x_{N-1})
\eqno{(4.1.8)}
$$
as follows from (4.1.5).

\bigskip

{\bf 4.2.} $\CL^{*}$ carries the standard affine covering
$\CU=\{U_{j},\;0\leq j\leq N-1\}$ defined by
$$
U_{j}=\{\text{class of }(u;x_{0},...,x_{N-1})\text{ s.t. } x_{j}\neq 0\}.
\eqno{(4.2.0)}
$$
Then $\pi\CU=\{\pi(U_{j}),\;0\leq j\leq N-1\}$ is (also standard)
affine covering of $\BP^{N-1}$ and $\pi\CU\cap\fF=\{\pi(U_{j})\cap\fF\}$
is an affine covering of $\fF$.

By definition, the cohomology of the \v Cech complex
$\check{C}^{*}(\pi\CU\cap\fF,\Lambda^{ch}\CT_{\fF})$ equals the cohomology
$H^{*}(\fF,\Lambda^{ch}\CT_{\fF})$. A more practical way to compute the latter
is provided by  Lemma 2.4.2. Indeed, being currently in the
situation of this lemma we obtain the bi-complex
$\check{C}^{*}(\CU,R^{*}\Lambda^{ch}\CT_{\CL^{*}};
d_{\check{C}},dt_{(0)})$, that is, the \v Cech complex
 with coefficients in the complex
$(\Lambda^{ch}\CT_{\CL^{*}},dt_{(0)})$ defined
in (2.4.5). Associated to this bi-complex
there are
 two standard spectral sequences, $'\!E^{p,q}_{r}$ and  $''\!E^{p,q}_{r}$,
such that
$$
\aligned
('\!E^{p,q}_{1},'\!d_{1})&=(H^{p}(\CL^{*},R^{q}( \Lambda^{ch}\CT_{\CL^{*}})),
dt_{(0)}),\\
'\!E^{p,q}_{2}&=
H^{q}_{dt_{(0)}}(H^{p}(\CL^{*},\Lambda^{ch}\CT_{\CL^{*}})),
\endaligned
\eqno{(4.2.1a)}
$$
$$
(''\!E^{p,q}_{1},''\!d_{1})=\check{C}^{p}(\CU,\CH^{q}_{dt_{(0)}}
(R^{q}( \Lambda^{ch}\CT_{\CL^{*}})).
\eqno{(4.2.1b)}
$$
In the latter $\CH^{q}_{dt_{(0)}}
(R^{*} \Lambda^{ch}\CT_{\CL^{*}})$ denotes the $q$-th cohomology
sheaf of complex (2.3.5).

{\bf 4.3. Lemma.} {\it Both  $'\!E^{p,q}_{r}$ and  $''\!E^{p,q}_{r}$ abut to
$H^{*}(\fF,\Lambda^{ch}\CT_{\fF})$.}

{\bf 4.4.} {\it Proof.} Observe that even though sheaf complex (2.4.5)
 appears to be infnite in both directions,
its differential preserves conformal weight (cf. 2.4.3) and it is easy to see
that each fixed conformal weight component of (2.4.5) is finite.
(Indeed, the differential $dt_{(0)}$ changes fermionic charge by one
and it follows from (2.3.5) that for a fixed conformal
weight fermionic charge may acquire only a finite number of values.)
This implies,
in a standard manner, that  both the spectral sequences abut
to the cohomology
of the total complex $C^{*}(\CU,R^{*}\Lambda^{ch}\CT_{\CL^{*}};d_{\check{C}}+dt_{(0)})$.
Lemma 2.4.2 implies that
$$
''\!E^{p,q}_{1}\iso\left\{\aligned 0&\text{ if }q\neq 0\\
C^{p}(\pi\CU\cap\fF,\Lambda^{ch}\CT_{\fF})&\text{ otherwise }.
\endaligned\right.
$$
Hence the second spectral sequence degenerates
to the \v Cech complex over $\fF$ and the lemma follows. $\qed$

{\bf 4.5.} We now wish to compute the 1st term of the 1st spectral sequence
recorded in
 (4.2.1a). Ignoring the double grading we rewrite (4.2.1a) as
$$
('\!E^{**}_{1},'\!d_{1})=(H^{*}(\CL^{*}, \Lambda^{ch}\CT_{\CL^{*}}),
dt_{(0)}),
\eqno{(4.5.1)}
$$
We have $\CL^{*}-0\hookrightarrow\CL^{*}$; this places us in the set-up of
3.7, see e.g. (3.7.1,2), and in order to compute
$(H^{*}(\CL^{*}, \Lambda^{ch}\CT_{\CL^{*}}))$ we employ the spectral
sequence of Lemma 3.8. -- in this particular case we shall be able to compute
all its terms. Begin with

{\bf 4.5.1.} {\it Toric description of $\CL^{*}$.}
Choose
$$
s=x_{0}^{N}, y_{j}=\frac{x_{j}}{x_{0}},\; 1\leq j\leq N
\eqno{(4.5.2)}
$$
to be coordinates of $\CL^{*}$ over the ``big cell'' $U_{0}$, see (4.2.0).
Let
$$
\aligned
&M_{\CL^{*}}=\BZ S\bigoplus\{\bigoplus_{i=1}^{N-1}\BZ Y_{i}\},\\
&M_{\CL^{*}}^{*}=\BZ S^{*}\bigoplus\{\bigoplus_{i=1}^{N-1}\BZ Y^{*}_{i}\},\\
&\Lambda_{\CL^{*}}=M_{\CL^{*}}\oplus M_{\CL^{*}}^{*},\;
\endaligned
\eqno{(4.5.3)}
$$
and the bases $\{S,Y_{1},...,Y_{N-1}\}$
and $\{S^{*},Y^{*}_{1},...,Y^{*}_{N-1}\}$ be dual to each other.

It follows from (4.5.2) and the identifications $s=e^{S}$, $y_{j}=e^{Y_{j}}$
that
a fan $\Sigma_{\CL^{*}}$ that defines $\CL^{*}$ can be chosen as follows:
the set of its 1-dimensional generators consists of
$$
S^{*}, NS^{*}-\sum_{i=1}^{N-1}Y^{*}_{i}, Y^{*}_{1}, Y^{*}_{2},..., Y^{*}_{N-1}.
\eqno{(4.5.4)}
$$
The set of the highest dimension cones consists of $N$ cones, each generated
by $S^*$ and the rest of the vectors in (4.5.4) except one of them. These
data determine $\Sigma_{\CL^{*}}$ uniquely.

{\bf 4.5.2.} {\it Computation of $(E^{**}_{1}, d_{1})$,
$(E^{**}_{2}, d_{2})$,..., $(E^{**}_{\infty}, d_{\infty})$.}

Proceeding along the lines of 3.7  we obtain
$$
(E(\CL^{*})^{**}_{*}, d_{*})\Rightarrow H^{*}(\CL^{*},
\Lambda^{ch}\CT_{\CL^{*}}).
\eqno{(4.5.4 1/2)}
$$
(3.7.5) reads (we skip $\CL^{*}$):
$$
(E^{p,q}_{1}, d_{1})=(H^{p}(\CL^{*}-0,\Lambda^{ch}\CT_{\CL^{*}-0})e^{qS^{*}},
(e^{qS^{*}}(z)\tilde{S}^{*}(z))_{(0)}).
\eqno{(4.5.5)}
$$
Thanks to (4.1.5) and the naturality of $\Lambda^{ch}\CT_{X}$, 2.1 (i), there
are canonical ismorphisms
$$
\aligned
&H^{p}(\CL^{*}-0,\Lambda^{ch}\CT_{\CL^{*}-0})\iso
H^{p}(\BC^{N}- 0/\BZ^{N},\Lambda^{ch}\CT_{\BC^{N}\setminus 0/\BZ^{N}})\\
&\iso
H^{p}(\BC^{N}- 0,\Lambda^{ch}\CT_{\BC^{N}- 0})^{\BZ_{N}}.
\endaligned
\eqno{(4.5.6)}
$$
The latter has been already computed, see 2.3.5, formulas (2.3.8,10,11).
Therefore,
$(E^{p,q}_{1}, d_{1})$ is as follows:
$$
E^{p,q}_{1}=\left\{
\aligned
\Lambda^{ch}\CT(\BC^{N})^{\BZ_{N}}e^{qS^{*}}&\text{ if }p=0\\
S_{1}(\Lambda^{ch}\CT(\BC^{N})^{\BZ_{N}})e^{qS^{*}}&\text{ if }p=N-1\\
0 &\text{ otherwise}
\endaligned
\right.
\eqno{(4.5.7a)}
$$
$$
d_{1}=(e^{S^{*}}(z)\tilde{S}^{*}(z))_{(0)}.
\eqno{(4.5.7b)}
$$
Note that $(E^{0,*}_{1}, d_{1})$ is canonically a differential graded
vertex algebra and  $(E^{N-1,*}_{1}, d_{1})$ is canonically
a differential graded  $(E^{0,*}_{1}, d_{1})$-module. Therefore,
as a vertex algebra,
$$
(E^{*,*}_{1}, d_{1})=(E^{0,*}_{1}, d_{1})\oplus (E^{N-1,*}_{1}, d_{1})
\eqno{(4.5.8)}
$$
is the abelian extension of the former by the latter. Therefore, (3.7.6)
reads as follows:
$$
E^{*,*}_{2}=H_{d_{1}}(E^{0,*}_{1}, d_{1})\oplus
 H_{d_{1}}(E^{N-1,*}_{1}, d_{1}),
\eqno{(4.5.9)}
$$
and again $H_{d_{1}}(E^{0,*}_{1}, d_{1})$ is a vertex algebra,
$H_{d_{1}}(E^{N-1,*}_{1}, d_{1})$ its module, and $E^{*,*}_{2}$
is the abelian extension of the former by the latter.

Thanks to (4.5.7a), the dimension consideration shows that
$$
\aligned
& d_{2}=d_{3}=\cdots d_{N-1}=0,\\
& E^{*,*}_{2}=E^{*,*}_{3}=\cdots =E^{*,*}_{N-3}=E^{*,*}_{N}.
\endaligned
\eqno{(4.5.10)}
$$
The same argument shows that the non-zero
components of the last non-zero differential are the following maps:
$$
d_{N}^{(i)}:\; E^{N-1,i}_{2}\rightarrow  E^{0,i+N}_{2}.
\eqno{(4.5.11)}
$$
As follows from 4.5.1, $\CL^{*}$ is covered by $N$ open affine sets.
Therefore, $H^{n}(\CL^{*},\Lambda^{ch}\CT_{\CL^{*}})=0$ if $n>N-1$.
This implies that the maps $ d_{N}^{(i)}$ are isomorphisms
if $i>0$ and $ d_{N}^{(0)}$ is an epimorphism. Hence
$$
E^{*,*}_{\infty}= E^{0,0}_{2}\oplus  E^{0,1}_{2}\oplus\cdots
\oplus  E^{0,N-2}_{2}\oplus
E^{0,N-1}_{2}\oplus \underbrace{\text{Ker}d_{N}^{(0)}}_{E^{N-1,0}_{\infty}}.
\eqno{(4.5.12)}
$$

By the spectral sequence definition, (4.5.12) implies
$$
\aligned
&H^{i}(\CL^{*},\Lambda^{ch}\CT_{\fF})\iso
H_{d_{1}}(E^{0,i}_{1}, d_{1}),\; 0\leq i\leq N-2,\\
0\rightarrow  H_{d_{1}}(E^{0,N-1}_{1}, d_{1})\rightarrow
&H^{N-1}(\CL^{*},\Lambda^{ch}\CT_{\fF})\rightarrow
\text{Ker}d_{N}^{(0)}\rightarrow 0.
\endaligned
\eqno{(4.5.13)}
$$

Now let us perform a change of coordinates that will reveal so far
invisible structure of  result (4.5.7ab, 4.5.13).

{\bf 4.6.} We have used the coordinates attached to $\CL^{*}$
by  definition. Now let us employ the map
$$
\BC^{N}-0\rightarrow \CL^{*}-0,
\eqno{(4.6.1)}
$$
that is, the composition of (4.1.5) and the natural projection
$\BC^{N}-0\rightarrow (\BC^{N}-0)/{\BZ_{N}}$. We would like to recast
the argument of 4.5 in terms inherent in $\BC^{N}-0$.

By invoking (4.6.1) we have placed ourselves in the situation of 3.9.
Let us make this explicit.

$\BC^{N}-0$ has the standard coordinate system $x_{i},\partial_{x_{i}}$
and has, therefore, the following toric description:
$$
M=\oplus_{i=0}^{N-1}\BZ X_{i},\;
M^{*}= \oplus_{i=0}^{N-1}\BZ X_{i}^{*},\; \Lambda=M\oplus M^{*},
\eqno{(4.6.2a)}
$$
so that
$$
X^{*}_{i}(X_{j})=\delta_{ij}.
\eqno{(4.6.2b)}
$$
$$
\Sigma=\{\sigma_{0},...,\sigma_{N-1}\},\; \sigma_{j}=\{x_{j}\neq 0\}.
\eqno{(4.6.2c)}
$$
Map (4.6.1) is induced by the lattice (cf. 3.2) embedding
$$
g: M^{*}\rightarrow M_{\CL^{*}}^{*}\; (=M_{\CL^{*}-0}^{*})
\eqno{(4.6.3)}
$$
dual to the lattice embedding
$$
\aligned
&g^{*}: M_{\CL^{*}}\rightarrow M,\\
&S\mapsto NX_{0},\; Y_{j}\mapsto X_{j}-X_{0},\; 1\leq j\leq N-1.
\endaligned
\eqno{(4.6.4a)}
$$
Indeed,   comparing (4.1.5) and (4.5.2) one obtains: under (4.6.1),
$s\mapsto x_{0}^{N}$, $y_{j}\mapsto x_{j}/x_{0}$ and it remains
to use $x_{i}=e^{X_{i}}$, $s=e^{S}$, $y_{j}=e^{Y_{j}}$ in order to obtain
(4.6.4a).

It is immediate to see that
$$
g(X_{0}^{*})=NS^{*}-\sum_{j=1}^{N-1}Y^{*}_{j},\;
g(X^{*}_{j})=Y^{*}_{j},\; j\geq 1;
\eqno{(4.6.4b)}
$$
thus $g(\Sigma)=\Sigma_{\CL^{*}}$. Therefore, we are indeed in the situation
of 3.9, and if we write down isomorphism (3.9.6)
explicitly, we will have all the assertions
of 4.5 recast in terms pertaining to $\BC^{N}-0$.

According to (3.1.7), there arises an isomorphism
$$
\hat{g}:\BB_{\Lambda_{\CL^{*}}} \rightarrow \BB_{g^{-1}\Lambda_{\CL^{*}}},
\eqno{(4.6.5)}
$$
where
$$
g^{-1}\Lambda_{\CL^{*}}=g^{*}M_{\Lambda_{\CL^{*}}}\oplus
g^{-1}M_{\Lambda_{\CL^{*}}}^{*}\subset\Lambda_{\BQ}.
$$
It follows from (4.6.4a) that
$$
g^{*}M_{\Lambda_{\CL^{*}}}=
\{ \sum_{i=0}^{N-1}m_{i}X_{i}\text{ s.t. }N| \sum_{i=0}^{N-1}\}\subset M.
\eqno{(4.6.6)}
$$
Inverting (4.6.4b) we obtain that
$g^{-1}M_{\Lambda_{\CL^{*}}}^{*}$ is spanned over $\BZ$ by
$$
g^{-1}(S^{*})=\frac{1}{N}(X^{*}_{0}+X^{*}_{1}+\cdots +X^{*}_{N-1}),
g^{-1}(Y^{*}_{j})=X^{*}_{j},1\leq j\leq N-1.
\eqno{(4.6.7)}
$$

Let us introduce the notation
$$
X^{*}_{orb}=\frac{1}{N}(X^{*}_{0}+X^{*}_{1}+\cdots +X^{*}_{N-1}).
\eqno{(4.6.9)}
$$
The first line of (4.5.7a) and (4.5.7b) rewrites as follows:
$$
\aligned
&(E^{0,*}_{1},d_{1})=(\oplus_{n=0}^{\infty}
\Lambda^{ch}\CT(\BC^{N})^{\BZ_{N}}e^{nX^{*}_{orb}},D_{orb})\\
&D_{orb}=\frac{1}{N}
\sum_{i=0}^{N-1}(\tilde{X}_{i}(z)e^{X^{*}_{orb}}(z))_{(0)}.
\endaligned
\eqno{(4.6.10)}
$$
(For the latter (3.1.4b) and (4.6.7) were used.)

This is a differential (w.r.t. $D_{orb}$) graded vertex algebra and because of
its importance and
its relation to Landau-Ginzburg models to be discovered later on,
we make a digression.

{\bf 4.6.1.} {\it Vertex algebra} $\widetilde{\text{LG}}_{orb}$.
Introduce the following notation:
$$
\widetilde{\text{LG}}_{orb}=\oplus_{n=0}^{\infty}\widetilde{\text{LG}}^{(n)}_{orb},
\; \widetilde{\text{LG}}^{(n)}_{orb}=\Lambda^{ch}\CT(\BC^{N})^{\BZ_{N}}e^{nX^{*}_{orb}}
\subset \BB_{g^{-1}\Lambda_{\CL^{*}}}.
\eqno{(4.6.11)}
$$
This vertex algebra is filtered by the system of differential vertex ideals
$$
\widetilde{\text{LG}}_{orb}^{\geq m}=\oplus_{n\geq m}
\widetilde{\text{LG}}_{orb}^{(n)}.
\eqno{(4.6.12)}
$$
Hence there arises a projective system of differential vertex algebras
$$
\widetilde{\text{LG}}_{orb}^{< m}=\widetilde{\text{LG}}_{orb}/
\widetilde{\text{LG}}_{orb}^{\geq m}.
\eqno{(4.6.13)}
$$
It is obvious that the natural projection
$$
\widetilde{\text{LG}}_{orb}\rightarrow
\widetilde{\text{LG}}_{orb}^{< m}
\eqno{(4.6.13)}
$$
induces the isomorphisms
$$
H^{i}_{D_{orb}}(\widetilde{\text{LG}}_{orb})\iso
H^{i}_{D_{orb}}(\widetilde{\text{LG}}_{orb}^{< m})
\text{ if } i<m-1.
\eqno{(4.6.14)}
$$
If $i=m-1$, then (4.6.13) induces the embedding
$$
H^{m-1}_{D_{orb}}(\widetilde{\text{LG}}_{orb})\hookrightarrow
H^{m-1}_{D_{orb}}(\widetilde{\text{LG}}_{orb}^{< m}),
\eqno{(4.6.15)}
$$
which is included in the following short exact sequence
$$
0\rightarrow H^{m-1}_{D_{orb}}(\widetilde{\text{LG}}_{orb})\rightarrow
H^{m-1}_{D_{orb}}(\widetilde{\text{LG}}_{orb}^{< m})
\rightarrow B^{m}\rightarrow 0,
\eqno{(4.6.16)}
$$
where $B^{m}\subset \widetilde{\text{LG}}_{orb}^{(m)}$ is the group
of $m$-couboundaries of
$\widetilde{\text{LG}}_{orb}$ w.r.t. $D_{orb}$.
 If we think of (4.6.15) as a 1-step filtration on
$H^{m-1}_{D_{orb}}(\widetilde{\text{LG}}_{orb}^{< m})$, then we get from (4.6.16)
$$
\text{Gr}H^{m-1}_{D_{orb}}(\widetilde{\text{LG}}_{orb}^{< m})=
 H^{m-1}_{D_{orb}}(\widetilde{\text{LG}}_{orb})\oplus B^{m}.
\eqno{(4.6.17)}
$$
$\qed$

\bigskip

Now, to the second line of (4.5.7a). The spectral flow transform
appearing there had been earlier interpreted in (3.10.4)
inside $\BB_{\Lambda}$, currently in use, by employing exactly
$X^{*}_{orb}$ although not its name.
 According to (3.10.4)
we have
$$
E^{N-1,q}_{1}=\Lambda^{ch}\CT(\BC^{N})^{\BZ^{N}}e^{(q+N)X^{*}_{orb}}.
$$
The complete translation of (4.5.7a) is then this:
$$
(E^{*,*}_{1}, d_{1})=(E^{0,*}_{1}\oplus E^{N-1,*}_{1},  d_{1}) =(\widetilde{\text{LG}}\oplus
\widetilde{\text{LG}}^{\geq N},D_{orb}).
\eqno{(4.6.18)}
$$
The rest of 4.5 carries over in a  straightforward manner:
(4.5.9) is translated as
$$
(E^{*,*}_{2}, d_{2})=(H^{*}_{D_{orb}}(\widetilde{\text{LG}})\oplus
H^{*}_{D_{orb}}(\widetilde{\text{LG}}^{\geq N}), 0),
\eqno{(4.6.19)}
$$
the 1st summand standing for $E^{0,*}_{2}$, the 2nd
for  $E^{N-1,*}_{2}$; (4.5.13) is translated as
$$
\aligned
&H^{i}(\CL^{*},\Lambda^{ch}\CT_{\CL^{*}})\iso H^{i}_{D_{orb}}(\widetilde{\text{LG}}),\; 0\leq i\leq N-2\\
0\rightarrow H^{N-1}_{D_{orb}}(\widetilde{\text{LG}}_{orb})\rightarrow
&H^{N-1}(\CL^{*},\Lambda^{ch}\CT_{\CL^{*}})\rightarrow B^{N}
\rightarrow 0,
\endaligned
\eqno{(4.6.20)}
$$
where we freely use the notation of digression 4.6.1. The first line
of (4.6.20) is obvious, and in the second only the identification
of $B^{N}$ with $\text{Ker}d^{(0)}_{N}$ needs an explanation. It follows from
(4.5.12) that the latter includes in the short exact sequence
$$
0\rightarrow\text{Ker}d^{(0)}_{N}\rightarrow
E^{N-1,0}_{2}\rightarrow E^{0,N}_{2}\rightarrow 0.
\eqno{(4.6.21)}
$$

According to (4.6.19), we have $E^{N-1,0}_{2}=Z^{N}$,
 $E^{0,N}_{2}=Z^{N}/B^{N}$, where $Z^{N}$, $B^{N}$ are the groups of $N$-cocycles
(coboundaries resp.) of the complex  $(\tilde{\text{LG}}_{orb}, D_{orb})$.
Therefore (4.6.21) can be identified with
$$
0\rightarrow B^{N}\rightarrow Z^{N}\rightarrow Z^{N}/B^{N}
\rightarrow 0.
\eqno{(4.6.22)}
$$
We have proved

{\bf 4.7. Theorem.} {\it The cohomology $H^{*}(\CL^{*},\Lambda^{ch}\CT_{\CL^{*}})$
satisfies:}
$$
H^{i}(\CL^{*},\Lambda^{ch}\CT_{\CL^{*}})\iso
H^{i}_{D_{orb}}(\widetilde{\text{LG}}_{orb}^{< N})\text{ if }i<N-1,
\eqno{(4.7.1)}
$$
{\it and if $i=N-1$, then there arises the short exact sequence}
$$
0\rightarrow H^{N-1}_{D_{orb}}(\widetilde{\text{LG}}_{orb})\rightarrow
H^{N-1}(\CL^{*},\Lambda^{ch}\CT_{\CL^{*}})\rightarrow B^{N}
\rightarrow 0.
\eqno{(4.7.2)}
$$
{\it In particular,}
$$
\text{Gr}H^{*}(\CL^{*},\Lambda^{ch}\CT_{\CL^{*}})\iso
\text{Gr}H^{*}_{D_{orb}}(\widetilde{\text{LG}}_{orb}^{< N}),
\eqno{(4.7.3)}
$$
{\it see (4.6.17) with $m=N$ for the definition of}
$\text{Gr}H^{*}_{D_{orb}}(\widetilde{\text{LG}}_{orb}^{< N})$.

{\bf 4.7.1. Remark.} Having put (4.16.14-16) on the table next to
(4.6.20) one observes that (4.7.1) has a good chance of being valid
 for $i=N-1$ as well.

{\bf 4.8.} Recall that our ultimate goal is the cohomology
vertex algebra $H^{*}(\fF,\Lambda^{ch}\CT_{\fF})$ and, according
to (4.2.1a), Theorem 4.7 only computes  the 1st term of the spectral
sequence converging to $H^{*}(\fF,\Lambda^{ch}\CT_{\fF})$. It remains
to write down explicitly its 2nd term, also see (4.2.1a), and for this we
need the differential $dt_{(0)}$ and the grading  $R^{q}(.)$
expressed in  terms of $\widetilde{\text{LG}}_{orb}$.

 Thanks to (4.1.8) the differential is as follows:
$$
\hat{g}(dt)= df(x_{0},...,x_{N-1}).
\eqno{(4.8.1)}
$$
A quick computation using, e.g., (1.9.9) shows that
$$
[D_{orb}, df(z)_{(0)}]=0.
\eqno{(4.8.2)}
$$

The grading $R^{q}(\widetilde{\text{LG}}_{orb})$, $q\in\BZ$,
 is nicely described as follows:
$$
R^{q}(\widetilde{\text{LG}}_{orb})=\text{Ker}(X^{*}_{orb,(0)}-qI).
\eqno{(4.8.3)}
$$
Indeed, $R^{q}(.)$ was defined in 2.4 as the eigenspace of the $\BC^{*}$-action.
In terms of $S,Y_{j}$, the infinitesimal generator of this
action is $S^{*}_{(0)}$; for example,
$S^{*}_{(0)}s=S^{*}_{(0)}e^{S}=e^{S}=s$,
 as it should
because according to (3.3.2a) $s=e^{S}$ is the coordinate along the fiber. It remains to use
(4.6.7).

Incidentally, the $\BZ$-grading built into definition (4.6.11)
of $\widetilde{\text{LG}}_{orb}$ is likewise given by the eigenvalues
of the operator $\sum_{j}X_{j,(0)}$. The two gradings are compatible
because $[X^{*}_{orb}, \sum_{j}X_{j,(0)}]=0$ as follows from (1.8.1).

Thus $\widetilde{\text{LG}}$ is  a bi-differential bi-graded
vertex algebra, to be denoted in this capacity
by  $(\widetilde{\text{LG}};D_{orb}, df(z)_{(0)})$,
and so is  $(\widetilde{\text{LG}}^{<N};D_{orb}, df(z)_{(0)})$.
Essentially it  remains
to summarize Lemma 4.3 and Theorem 4.7.

{\bf 4.9. Theorem.} {\it (i) Spectral sequence (4.2.1a) abuts to
$H^{*}(\fF,\Lambda^{ch}\CT_{\fF})$ so that}

$$
'\!E^{p,q}_{2}=H^{q}_{df(z)_{(0)}}(H^{p}_{D_{orb}}
(\widetilde{\text{LG}}_{orb}^{< N})),\; \text{ if } 0\leq p\leq N-2,
\eqno{(4.9.1)}
$$
{\it and if $p=N-1$, then $ '\!E^{N-1,q}_{2}$  is included into the long
exact sequence stemming from the short exact sequence of complexes,
cf. (4.7.2)}
$$
0\rightarrow (H^{N-1}_{D_{orb}}(\widetilde{\text{LG}}_{orb}), df(z)_{(0)})
\rightarrow ('\!E^{*,N-1-*}_{1}, d_{1})\rightarrow (B^{N},  df(z)_{(0)})
\rightarrow 0.
\eqno{(4.9.2)}
$$
{\it (ii) In the conformal weight zero component the spectral sequence degenerates
in the 2nd term and gives the isomorphism}
$$
\text{Gr}H^{*}(\fF,\Lambda^{*}\CT_{\fF})\iso
\underbrace{\BC\oplus\cdots\oplus\BC}_{N-1}\oplus
\left(\BC[x_{0},...,x_{N-1}]/<df>\right)^{\BZ_{N}}.
\eqno{(4.9.3)}
$$
{\bf 4.9.1. Remark.}  For the same reason that was indicated in Remark 4.7.1,
isomorphisms (4.9.1) have a good chance of being valid for
$p=N-1$ as well. Were this the case, (4.9.2) would be unnecessary.
$\qed$

\bigskip

{\it Beginning of the proof.} Item (i) is indeed simply (4.2.1a), Lemma 4.3,
and (4.7.1-2) of
Theorem 4.7 put together. The isomorphism in (ii) is well known classically
 and
the degeneration assertion follows very easily. We prefer, however, to
emphasize some additional  structure hidden in the spectral sequence
and then use it to give, among other things, a self-contained proof of (ii),
see 4.13.

\begin{sloppypar}
{\bf 4.10.} {\bf Addendum: $N2$-structure.}
 Condition (4.1.0) ensures
that $\CL^*$ is the canonical line bundle and this places us in the situation of
2.4.4: coordinates
$s,y_{1},...,y_{N-1}$ of
 (4.5.2) satisfy the conditions imposed in 2.4.4, formulas
(2.4.7) define an $N2$-structure on $\Lambda^{ch}\CT_{\CL^{*}}$
and (2.4.8) define that on $\Lambda^{ch}\CT_{\fF}$. In order to
compute this structure in terms of $\widetilde{\text{LG}}_{orb}$
one has to do the following: first, compute the images of (2.4.8)
in $\BB_{\Lambda_{\CL^{*}}}$ under  Borisov's map (3.3.2a,b);
second, apply  map  (4.6.5) to the result. This is straightforward
and tedious but rewarding, the reward being the coincidence of the
result with the free field realization that Witten related with
the Landau-Ginzburg model; this coincidence will be verified in
 Lemma 5.2.14.
 \end{sloppypar}

In order to record the result it is convenient to use the boson-fermion correspondence, 1.13. Let us introduce the standard lattice $\BZ^{N}$
so that the standard basis $\chi_{0},...,\chi_{N-1}$ is orthonormal.
Then, see 1.13,
one can make identifications
$$
\tilde{X}_{i}(z)= e^{\chi_{i}}(z),\;
\tilde{X}_{i}^{*}(z)= e^{-\chi_{i}}(z).
\eqno{(4.10.1)}
$$
{\bf 4.10.1. Lemma.} {\it The $N2$-structure on
$H^{*}(\fF,\Lambda^{ch}\CT_{\fF})$  comes from the following
$N2$-structure on} $\widetilde{\text{LG}}_{orb}$:
$$
\aligned
&G(z)\mapsto \sum_{j=0}^{N-1}
:\left(X_{j}^{*}(z)-\chi_{j}(z)\right)e^{\chi_{j}}(z):,\;
Q(z)\mapsto \sum_{j=0}^{N-1}
:\left(X_{j}(z)+\frac{1}{N}\chi_{j}(z)\right)e^{-\chi_{j}}(z):,\\
&J(z)\mapsto \sum_{j=0}^{N-1}\left(-\frac{1}{N}X^{*}_{j}(z)+X_{j}(z)+
\chi_{j}(z):\right),\\
&L(z)\mapsto \sum_{j=0}^{N-1}
\left( :X_{j}(z)X^{*}_{j}(z):+\frac{1}{2}:\chi_{j}(z)^{2}:
-\frac{1}{2}\chi_{j}(z)'-X_{j}(z)'\right).
\endaligned
$$
{\bf Proof.} $N2$ is generated as a Lie algebra by 2 fields,
$G(z)$ and $Q(z)$.
Let us do $Q(z)$, the field that acquires the geometrically
mysterious factor $1/N$. We have starting with (2.4.8)
$$
\aligned
&Q(z)\mapsto
s(z)'\partial_{ds}(z)+\sum_{j=1}^{N-1}y_{j}(z)'\partial_{dy_{j}}(z)-
(s(z)\partial_{ds}(z))'=\\
&\sum_{j=1}^{N-1}:y_{j}(z)'\partial_{dy_{j}}(z):-
:s(z)(\partial_{ds}(z))':=\\
&\sum_{j=1}^{N-1}:e^{Y_{j}}(z)'(:e^{-Y_{j}}(z)\tilde{Y}^{*}_{j}(z):)-
:e^{S}(z)(:e^{-S}(z)\tilde{S}^{*}_{j}(z):)':=\\
&\sum_{j=1}^{N-1}:Y_{j}(z)\tilde{Y}^{*}_{j}(z): +:S(z)\tilde{S}^{*}(z):-
\tilde{S}^{*}(z)'=\\
&\sum_{j=1}^{N-1}:(X_{j}(z)-X_{0}(z))\tilde{X}^{*}_{j}(z): +:NX_{0}(z)\frac{1}{N}\sum_{j=0}^{N-1}\tilde{X}^{*}_{j}(z):-
\frac{1}{N}\sum_{j=0}^{N-1}\tilde{X}^{*}_{j}(z)'=\\
&\sum_{j=0}^{N-1}:X_{j}(z)\tilde{X}^{*}_{j}(z):
-\frac{1}{N}\sum_{j=0}^{N-1}\tilde{X}^{*}_{j}(z)'=\\
&\sum_{j=0}^{N-1}
:\left(X_{j}(z)+\frac{1}{N}\chi_{j}(z)\right)e^{-\chi_{j}}(z):.
\endaligned
$$
A brief guide to this computation is as follows: the 3rd line is
(3.3.2a-b) applied to the 1st line, in the 5th line transformation
formulas (4.6.7) are used, and  boson-fermion correspondence
(4.10.1) is used in the 7th. In addition, the well-known
differentiation formula $e^{\alpha}(z)'=:\alpha(z)e^{\alpha}(z):$
has been repeatedly employed. $\qed$

{\bf 4.11. Character and Euler character formulas.}
Whenever one has a bi-graded vector space
$$
V=\bigoplus_{m,n} V^{m}_{n},\; \text{dim}V^{m}_{n}<\infty,
$$
one can define its character:
$$
chV(s,\tau)=\sum_{m,n} e^{2\pi i (ms+n\tau)}\text{dim}V^{m}_{n},
\eqno{(4.11.1)}
$$
and if in addition $V$ is a supervector space, one can define its
Euler character
$$
\text{Eu}V(s,\tau)=\sum_{m,n} e^{2\pi i (ms+n\tau)}
s\text{dim}V^{m}_{n},
\eqno{(4.11.2)}
$$
where the super-dimension $s\text{dim}V^{m}_{n}$
 is defined in the standard manner to be the dimension of the even
component of $V^{m}_{n}$ minus the dimension of its odd component. Note
that if $V$ carries an odd differential preserving the bi-grading, then
$$
\text{Eu}H_{d}(V)(s,\tau)=\text{Eu}V(s,\tau).
\eqno{(4.11.3)}
$$
As an example, we can consider $\widetilde{\text{LG}}_{orb}^{<N}$
bi-graded by the eigenvalues of $L_{(1)}$, $J_{(0)}$, see
Lemma 4.10.1. A straightforward computation
(which, however, we postpone until 5.2.16) shows that if we introduce
$$
E(\tau,s)=\prod_{n=0}^{\infty}\frac{\left(1-e^{2\pi
i\left(\left(n+1\right)\tau+\left(1-1/N\right)s\right)}\right)^{N}\left(1-e^{2\pi
i\left(n\tau+\left(-1+1/N\right)s\right)}\right)^{N}}{\left(1-e^{2\pi
i\left(\left(n+1\right)\tau+s/N\right)}\right)^{N}\left(1-e^{2\pi
i\left(n\tau-s/N\right)}\right)^{N}},
$$
then
$$
\text{Eu}\widetilde{\text{LG}}_{orb}^{<N}(\tau,s)
= \frac{1}{N}\sum_{l=0}^{N-1}\sum_{j=0}^{N-1}
e^{\pi
i\left(N-2\right)\left\{-2js+\left(j^{2}-j\right)\tau+j^{2}\right\}}E(\tau,s-j\tau-l).
\eqno{(4.11.4)}
$$
and if we introduce
$$
\tilde{E}(\tau,s)=\prod_{n=0}^{\infty}\frac{\left(1+e^{2\pi
i\left(\left(n+1\right)\tau+\left(1-1/N\right)s\right)}\right)^{N}\left(1+e^{2\pi
i\left(n\tau+\left(-1+1/N\right)s\right)}\right)^{N}}{\left(1-e^{2\pi
i\left(\left(n+1\right)\tau+s/N\right)}\right)^{N}\left(1-e^{2\pi
i\left(n\tau-s/N\right)}\right)^{N}},
$$
then
$$
ch\widetilde{\text{LG}}_{orb}^{<N}(\tau,s)
= \frac{1}{N}\sum_{l=0}^{N-1}\sum_{j=0}^{N-1}
e^{\pi
i\left(N-2\right)\left\{-2js+\left(j^{2}-j\right)\tau\right\}}
\tilde{E}(\tau,s-j\tau-l).
\eqno{(4.11.5)}
$$
The repeated application of (4.11.3) to (4.7.3)
 shows that result (4.11.4) is valid
for $H^{*}(\fF,\Lambda^{ch}\CT_{\fF}) $ as well:
$$
\text{Eu} H^{*}(\fF,\Lambda^{ch}\CT_{\fF}) =
\frac{1}{N}\sum_{l=0}^{N-1}\sum_{j=0}^{N-1}
e^{\pi
i\left(N-2\right)\left\{-2js+\left(j^{2}-j\right)\tau+j^{2}\right\}}
E(\tau,s-j\tau-l).
\eqno{(4.11.6)}
$$
The importance of the latter is that, as was observed in [BL], the elliptic
genus of $\fF$, $\text{Ell}_{\fF}(\tau,s)$, satisfies
$$
\text{Ell}_{\fF}(\tau,s)=e^{\pi i (N-2)s}
\text{Eu} H^{*}(\fF,\Lambda^{ch}\CT_{\fF}) .
\eqno{(4.11.7)}
$$
We have proved

{\bf 4.12. Corollary.}
$$
\text{Ell}_{\fF}(\tau,s)=\frac{1}{N}\sum_{l=0}^{N-1}\sum_{j=0}^{N-1}
e^{\pi
i\left(N-2\right)\left\{-js+\left(j^{2}-j\right)\tau+j^{2}\right\}}
E(\tau,s-j\tau-l).
\eqno{(4.11.8)}
$$

{\bf 4.13. Chiral rings and $ H^{*}(\fF,\Lambda^{*}\CT_{\fF})$.}

{\it 4.13.1. End of proof of Theorem 4.9.}
Arguments that led to Theorems 4.7, 9
consist of computations of two spectral sequences. Both the sequences
are graded by $L_{(1)}$, see e.g. Lemma 4.10.1, and all the differentials
preserve this grading.  The classical story is about the events unfolding
in the conformal weight zero component. Let us re-tell this story.

Therefore, we adopt that point of view according to which our task is
to compute the cohomology algebra $ H^{*}(\fF,\Lambda^{*}\CT_{\fF})$
of polyvector fields. Investigated
in Theorem 4.9 is spectral sequence (4.2.1a); the origin
of its conformal weight zero component is the classical Koszul
sheaf complex (2.4.6) and its 1st term is the  algebra
$ H^{*}(\CL^{*},\Lambda^{*}\CT_{\CL^{*}})$. The computation of this
alegbra is accomplished in the classical part of Theorem 4.7 and that is where
an important simplification occurs:
  (4.7.3) restricted to  the conformal weight zero component is valid
without the passage to the graded objects:
$$
H^{*}(\CL^{*},\Lambda^{ch}\CT_{\CL^{*}})_{0}\iso
H^{*}_{D_{orb}}(\widetilde{\text{LG}}_{orb}^{< N})_{0},
\eqno{(4.13.1)}
$$
In fact, even taking the $D_{orb}$-cohomology is not necessary:
$$
H^{*}(\CL^{*},\Lambda^{ch}\CT_{\CL^{*}})_{0}\iso
(\widetilde{\text{LG}}_{orb}^{< N})_{0}.
\eqno{(4.13.2)}
$$
Indeed, by definition
$$
(\widetilde{\text{LG}}_{orb}^{(0)})_{0}=\Lambda^{*}\CT(\BC^{N})^{\BZ_{N}}.
\eqno{(4.13.3)}
$$
Further,
$$
\text{dim}(\widetilde{\text{LG}}_{orb}^{(i)})_{0}=1,\; 1\leq i\leq N-1,
\eqno{(4.13.4)}
$$
as easily follows from  character formula (4.11.5). The corresponding generator
is
$$
(\widetilde{\text{LG}}_{orb}^{(i)})_{0}= \BC
e^{iX^{*}_{orb}-\sum_{j}(X_{j}+\chi_{j})},\; 1\leq i\leq N-1,
\eqno{(4.13.5)}
$$
where we again use boson-fermion correspondence (4.10.1).
Now observe that  spaces (4.13.3) and that spanned by elements
 (4.13.5) are both annihilated by
$D_{orb}$, and (4.13.2) follows. We  also see that the 1st term of
spectral sequence (4.2.1a) is
$$
\left(('\!E_{1}^{*,*})_{0},d_{1}\right)=\left(
\underbrace {\BC\oplus\cdots\oplus\BC}_{N-1}\oplus
\left(\Lambda^{*}\CT(\BC^{N})\right)^{\BZ_{N}},
df(z)_{(0)}\right).
\eqno{(4.13.6)}
$$
Note that in view of (4.13.2), (4.9.1) is valid for $p=N-1$
as well.
A quick computation shows that
the elements (4.13.5) are annihilated by $df(z)_{(0)}$.
 As to component (4.13.3), we have:

$(\Lambda^{*}\CT(\BC^{N}),df(z)_{(0)})$ {\it is the standard Koszul
resolution of the Milnor ring. }

By virtue of (4.13.6),
$$
('\!E_{2}^{*,*})_{0}=
\underbrace {\BC\oplus\cdots\oplus\BC}_{N-1}\oplus
\left(\BC[x_{0},...,x_{N-1}]/<df>\right)^{\BZ_{N}}.
\eqno{(4.13.7)}
$$
Finally, all the higher differentials vanish simply because
the elements listed in (4.13.5,7)  are genuine cocycles.
This completes the proof of Theorem 4.9. $\qed$

It is amusing to note that we have obtained ``vertex'' representatives
of all the classes of the cohomology $H^{*}(\fF,\Lambda^{*}\CT_{\fF})$.
Since, as was explained in 4.8, the eigenvalues of
$X^{*}_{orb}+\sum_{j}X_{j,(0)}$ give the cohomological grading, we have:
$$
(\text{class of}) e^{iX^{*}_{orb}-\sum_{j}(X_{j}+\chi_{j})}\in
H^{i-1}(\fF,\Lambda^{N-i-1}\CT_{\fF}),\; 1\leq i\leq N-1,
\eqno{(4.13.8)}
$$
$$
(\text{class of})\prod_{j}x_{j}^{m_{j}}=e^{\sum_{j}m_{j}X_{j}}\in
H^{m}(\fF,\Lambda^{m}\CT_{\fF}),\;
m=\frac{1}{N}\sum_{j}m_{j}, 0\leq m_{j}\leq N-2.
\eqno{(4.13.9)}
$$

{\it 4.13.2. Multiplicative structure.}
The multiplicative structure of the ring $H^{*}(\fF,\Lambda^{*}\CT_{\fF})$
is well known, of course. But let us restore it by the ``vertex''
methods.

According to (2.3.3b),
 the chiral ring
of $H^{*}(\fF,\Lambda^{ch}\CT_{\fF})$ is isomorphic to
$H^{*}(\fF,\Lambda^{*}\CT_{\fF})$. However, what we have at our disposal
is the chiral ring $(\widetilde{\text{LG}}_{orb}^{< N})_{0}$ and as a ring
it only gives $\text{Gr}H^{*}(\fF,\Lambda^{ch}\CT_{\fF})$, as we have just
proved -- this is a common problem with spectral sequences.

Nevertheless, having re-examined the way in which spectral
sequence (4.2.1a) was defined, one concludes easily that the 0-th
component, $(\widetilde{\text{LG}}_{orb}^{(0)})_{0}$,
 remains unaffected by the passage
to the graded object and thus carries the ``right'' multiplication; therefore,
the Milnor ring, $\BC[x]/<df>$, embeds into $H^{*}(\fF,\Lambda^{*}\CT_{\fF})$
as a ring, cf. (4.9.3).

Let us now look at  elements (4.13.8).
The  computation
$$
\aligned
\left(e^{iX^{*}_{orb}-\sum_{j}(X_{i}+\chi_{i})}\right)_{(-1)}
e^{\sum_{j}m_{j}X_{j}}=&\lim_{z\rightarrow w}(z-w)^{i\sum_{j}\frac{m_{j}}{N}}
e^{iX^{*}_{orb}-\sum_{j}(X_{j}+\chi_{j})}=\\
&\left\{\aligned
0&\text{ if }\sum_{j}m_{j}>0\\
e^{iX^{*}_{orb}-\sum_{j}(X_{j}+\chi_{j})}&\text{ if
}\sum_{j}m_{j}=0,
\endaligned\right.
\endaligned
\eqno{(4.13.10)}
$$
as follows from (1.9.8), is valid even in $(\widetilde{\text{LG}}_{orb}^{(0)})_{0}$;
hence in $H^{*}(\fF,\Lambda^{*}\CT_{\fF})$ as well.

Likewise, if $s+t<N$, one obtains
$$
\aligned
&\left(e^{sX^{*}_{orb}-\sum_{j}(X_{j}+\chi_{j})}\right)_{(-1)}
e^{tX^{*}_{orb}-\sum_{j}(X_{j}+\chi_{j})}=\\
&(-1)^{s}\lim_{z\rightarrow w}(z-w)^{N-s-t}
e^{(s+t)X^{*}_{orb}-2\sum_{j}(X_{j}+\chi_{j})}=0
\endaligned
\eqno{(4.13.11)}
$$
inside $(\widetilde{\text{LG}}_{orb}^{(0)})_{0}$, hence inside
$H^{*}(\fF,\Lambda^{*}\CT_{\fF})$ as well.

Finally, the same computation shows that
$$
\left(e^{sX^{*}_{orb}-\sum_{j}(X_{j}+\chi_{j})}\right)_{(-1)}
e^{(N-s)X^{*}_{orb}-\sum_{j}(X_{j}+\chi_{j})}= (-1)^{s}
e^{NX^{*}_{orb}-2\sum_{j}(X_{j}+\chi_{j})}
\eqno{(4.13.12)}
$$
In $(\widetilde{\text{LG}}_{orb}^{< N})_{0}$ the class of this element
is zero because it belongs to $\widetilde{\text{LG}}_{orb}^{(N)}$
and this component was cut off in definition (4.6.13).
 However, an amusing diagram
search shows (we skip this computation)
 that in reality this element is cohomologous to
$\pm\prod_{j}x_{j}^{N-2}$, which is a generator of
$H^{N-2}(\fF,\Lambda^{N-2}\CT_{\fF})$, see (4.13.9). Therefore,
$$
\left(\text{class of }e^{sX^{*}_{orb}-\sum_{j}(X_{j}+\chi_{j})}\right)_{(-1)}
\left(\text{class of }e^{(N-s)X^{*}_{orb}-\sum_{j}(X_{j}+\chi_{j})}
\right)=\pm
\prod_{j}x_{j}^{N-2}
\eqno{(4.13.13)}
$$
and gives a non-degenerate pairing
$$
H^{s-1}(\fF,\Lambda^{N-s-1}\CT_{\fF})\times
H^{N-s-1}(\fF,\Lambda^{s-1}\CT_{\fF})\rightarrow
H^{N-1}(\fF,\Lambda^{N-1}\CT_{\fF}),
\eqno{(4.13.14)}
$$
as it should.
This completes the description of multiplication on
$H^{*}(\fF,\Lambda^{*}\CT_{\fF})$.

\bigskip\bigskip

\centerline{{\bf 5. Landau-Ginzburg orbifolds.}}

\bigskip

In this section we shall provide the necessary definitions and
constructions so as to identify the complex
$(\widetilde{\text{LG}}_{orb}^{<N}, df(z)_{(0)})$, which played an
important role in Theorem 4.7 and 9, with the chiral algebra of
the Landau-Ginzburg orbifold.

{\bf 5.1. Landau-Ginzburg model.} Let
$$
f\in\BC[x_{0},...,x_{N-1}],\; \text{deg}f=p
\eqno{(5.1.1)}
$$
be a homogeneous polynomial such that  its partials
 $\partial_{x_{i}}f$, $0\leq i\leq N-1$, have only one common zero
occurring at $\vec{x}=0$. For the
time being this $f$ need not be identified with $f$
of (4.1.6), that is, $p$ need not be N, but
in the main application this assumption will be made.

By the (chiral algebra of the) Landau-Ginzburg model associated to
$\BC^{N}$ and $f$ as in (5.1.1) we understand the differential
vertex algebra
$$
\text{LG}_{f}=(\Lambda^{ch}\CT(\BC^{N}),df(z)_{(0)}).
\eqno{(5.1.2)}
$$
Note  that the chiral ring of $\text{LG}_{f}$
is the standard Koszul resolution:
$$
K_{f}^{*}:\; (\text{LG}_{f})_{0}=(\Lambda^{*}\CT(\BC^{N}),df),
\eqno{(5.1.3)}
$$
the differential being the contraction with the 1-form $df$, cf. (4.13.6).

The following lemma is an important ingredient in Witten's approach [W2]
to the Landau-Ginzburg model; exactly this form of the result
appears as formula (3.1.1) in [KYY].

{\bf 5.1.1. Lemma.} {\it The assignment}
$$
\rho:
\aligned
&G(z)\mapsto\sum_{i=0}^{N-1}\partial_{x_{i}}(z)dx_{i}(z),\;
Q(z)\mapsto \sum_{i=0}^{N-1}-\frac{1}{p}x_{i}(z)\partial_{dx_{i}}(z)'
-(\frac{1}{p}-1)x_{i}(z)'\partial_{dx_{i}}(z),\\
&J(z)\mapsto \sum_{i=0}^{N-1}-\frac{1}{p}:\partial_{x_{i}}(z)x_{i}(z):+
(\frac{1}{p}-1):\partial_{dx_{i}}(z)dx_{i}(z):,\\
&L(z)\mapsto \sum_{i=0}^{N-1}:\partial_{x_{i}}(z)x_{i}(z)':+
:\partial_{dx_{i}}(z)'dx_{i}(z):
\endaligned
$$
{\it determines a vertex algebra morphism}
$$
\rho:\; V(N2)_{N\frac{p-2}{p}}\rightarrow \text{LG}_{f}
\eqno{(5.1.4)}
$$
{\it such that}
$$
df(z)_{(0)}\rho\left(V(N2)_{N\frac{p-2}{p}}\right)=0.
\eqno{(5.1.5)}
$$

{\bf 5.1.2. Corollary.} {\it The vertex algebra
$H_{df(z)_{(0)}}(\text{LG}_{f})$
carries an $N2$-structure inherited from that on $\text{LG}_{f}$.}

This follows at once from Lemma 5.1.1 and (1.3.1).
Let us look at some basic examples.

{\bf 5.1.3. Theorem} ([FS]) {\it If $N=1$, $f=x^{p}$, then}
$H_{df(z)_{(0)}}(\text{LG}_{f})$ {\it is the direct sum of $p-1$
unitary $N2$-modules generated by $1,x,x^{2},...,x^{p-2}$.}

Denote by $U_{i,p}$ the unitary $N2$-module generated according to
 Theorem 5.1.3 by $x^{i}$, $0\leq i\leq p-2$;  here
$p$ keeps track of the ``central charge'', that is, the value by
which the central element  $C$ operates on the module; in this case
$C\mapsto (p-2)/p$. These modules are pairwise
non-isomorphic.

{\bf 5.1.4. Corollary.} {\it If $f=\sum_{j}x_{j}^{p}$, then}
$$
H_{df(z)_{(0)}}(\text{LG}_{f})=\bigoplus_{0\leq j_{0},...,j_{N-1}\leq p-2}
\bigotimes_{t=0}^{N-1}U_{j_{t},p}
$$
{\it which is a unitary $N2$-module w.r.t. the diagonal action of central
charge $N(p-2)/p$.} $\qed$

{\bf Proof} follows at once from Theorem 5.1.3 because in the case
of the ``diagonal'' $f$ the complex $\text{LG}_{f}$ is the tensor product
of the complexes of Theorem 5.1.3 -- hence so is its cohomology
as follows from the K\"unneth formula.

{\bf 5.1.5. Remark.} The reader will  notice that
complex (5.1.2)
is nothing but (2.4.5) computed in the purely local situation with the function
$t$ replaced with $f$. Furthermore, the cohomology of (2.4.5)
gives $\Lambda^{ch}\CT_{Z(t)}$, see Lemma 2.4.2, and $Z(t)$ is exactly
the singular locus of $t$: for any affine $U\subset \CL^{*}$
$$
Z(t)\cap U=\text{Spec}\left(\CO_{\CL^{*}}(U)/<dt>\right).
$$
Thus one is tempted to set
$$
\Lambda^{ch}\CT_{\text{Spec}M_{f}}(\text{Spec}M_{f})=
H_{df(z)_{(0)}}(\text{LG}_{f}),
$$
thereby defining the sheaf $\Lambda^{ch}\CT_{\text{Spec}M_{f}}$,
where $M_{f}=\BC[x_{0},...,x_{N-1}]/<df>$.

To put this somewhat differently, we have resolved the singularity of
$M_{f}$ by passing to the DGA $K_{f}^{*}$, see (5.1.3), and then chiralized
the latter so as to obtain (5.1.2). There seems to be a natural construction
[KV2]
allowing to chiralize in a similar manner other free DGA's thereby extending
algebras of chiral polyvector fields from smooth varieties to
spectra of Milnor rings to a wider class of schemes.

{\bf 5.2. Landau-Ginzburg orbifold.}

{\bf 5.2.1.} Familiar in  vertex algebra
 theory is the following
pattern: $V$ is a vertex algebra; $g$ is its order $N$ automorphism;
$V^{(i)}$ is a ``naturally defined'' $g^{i}$-twisted $V$-module,
$1\leq i\leq N-1$; assuming that the group $\{1,g,...,g^{N-1}\}$
also acts on each $V^{(i)}$, naturally w.r.t to the $V$-action, one forms
the vertex algebra of $g$-invariants, $V^{g}$, and its (untwisted) modules $(V^{(i)})^{g}$,
$1\leq i\leq N-1$. It sometimes so happens that the space
$$
V^{g}\oplus (V^{(1)})^{g} \oplus\cdots (V^{(N-1)})^{g}
\eqno{(5.2.1)}
$$
carries an ``interesting'' vertex algebra structure compatible
with the described $V^{g}$-module structure. Vertex algebra
(5.2.1) is often referred to as an orbifold or a $V$ orbifold.

The most famous example of an orbifold
is undoubtedly the Monster vertex algebra $ V^{\text{Mnstr}} $
[FLM]. Indeed,
$$
 V^{\text{Mnstr}} = V_{L}^{g}\oplus V_{L,1}^{g},
$$
where $L$ is the Leech lattice, $g$ its involution, and
$ V_{L,1}$ is an irreducible
 $g$-twisted $V_{L}$-module (unique by Dong's theorem [D2]).

\begin{sloppypar}
A number of physics papers, an incomplete list including [V,VW,W2]
and references therein, suggests that a realization of this
pattern in the case where $V=\text{LG}_{f}$, $\text{deg}f=N$, and
$g=\exp{(2\pi i \rho J_{0})}$, $\rho J(z)$ being defined in Lemma
5.1.1, may  be closely related to ``string vacua''. The remainder
of the paper is our attempt to understand this idea. We shall,
first, construct the candidates for $\text{LG}_{f}^{(i)}$, and in
order to do so we shall need a recollection on vertex algebra
twisted modules, the notion  introduced in [FFR]. In our
presentation we shall mostly follow [KR]. In 5.2.14 we shall
conclude that the space
$\oplus_{i=0}^{N-1}\left(\text{LG}_{f}^{(i)}\right)^{g}$ with
$\text{deg}f=N$ coincides with $\widetilde{\text{LG}}^{<N}_{orb}$
of Theorems 4.7, 4.9 as an $N2$-module. Second, in 5.2.17-5.2.20,
we shall make
$\oplus_{i=0}^{N-1}\left(\text{LG}_{f}^{(i)}\right)^{g}$ into a
vertex algebra such that its chiral ring is isomorphic to the
cohomology ring of the polyvector fields,
$H_{\fF}(\fF,\Lambda^{*}\CT_{\fF})$, on the corresponding
Calabi-Yau hypersurface  (4.1.6).
\end{sloppypar}

{\bf 5.2.2. Twisting data.}

 Let $G$ be an additive subgroup of $\BC$ containing $\BZ$.
A vertex algebra $V$ is called $G/\BZ$-graded if
$$
V=\oplus_{\bar{m}\in G/\BZ}V[\bar{m}],
\eqno{(5.2.2a)}
$$
so that
$$
V[\bar{m}]_{(n)}V[\bar{l}]\subseteq V[\bar{m}+\bar{l}].
\eqno{(5.2.2b)}
$$
It is clear that $V[0]\subset V$ is a vertex subalgebra.

 To give an example, let $g$ be an order $N$ automorphism of
a vertex algebra $V$. Let $G=\frac{1}{N}\BZ$. We have $G/\BZ\iso
\BZ_{N}$.

Then
$$
V=\oplus_{m=0}^{N-1} V[m]=\{v\in V: gv=e^{2\pi im/N}v\}
\eqno{(5.2.3)}
$$
is a $\BZ_{N}$-grading. By definition, in this case $V[0]$ is the
vertex subalgebra of $g$-invariants, $V^{g}$.

 Let $W$ be a vector space and $\bar{m}\in G/\BZ$,
where $G$ is as in 2.1.1. An $\bar{m}$-twisted $\text{End} M$-valued
 field is a series
$$
a(z)=\sum_{m\in \bar{m}}a_{(m)}z^{-m-1},
$$
where $a_{(m)}\in\text{End} M$ is such that for any $v\in W$,
$a_{(m)}v=0$ if $\text{Re } m>>0$. Let $\text{Field}_{G}(W)$ be
 the linear
space of $m$-twisted $\text{End} W$-valued fields for all $m\in G/\BZ$.

{\bf 5.2.3. Definition.} (cf. Definition 1.2 and [KR, sect.5)
 A $G$-twisted $V$-module $W$ is a parity preserving linear map
$$
\rho: V\rightarrow \text{Field}_{G}(W),\;
(\rho a)(z)=\sum_{m}\rho a_{(m)}z^{-m-1}
$$
satisfying the following axioms:

(i) if $a\in V[\bar{m}]$, then $(\rho a)(z)$ is $\bar{m}$-twisted;

(ii) $\rho (\b1)=\text{id}$;

(iii) (twisted  Borcherds identity)
for any $a\in V[\bar{m}]$, $b\in V$ and
$F(z,w)=z^{m}(z-w)^{l}$ such that $m\in\bar{m}$, $l\in\BZ$
$$
\aligned
&\text{Res}_{z-w}\rho(a(z-w)b,w)i_{w,z-w}F(z,w)\\
&=
\text{Res}_{z}\left((\rho a)(z)(\rho b)(w)i_{z,w}F(z,w)-
(-1)^{\pr(a)\pr(b)}(\rho b)(w)(\rho a)(z)i_{w,z}F(z,w)\right).
\endaligned
\eqno{(5.2.4)}
$$

{\bf 5.2.4. Remarks.}

(i) Note that the $l=0$ case of the twisted Borcherds identity
is the following commutator formula
$$
[\rho a_{(m)},\rho b_{(k)}] =\sum_{j=0}^{\infty}\binom{m}{j}
\rho (a_{(j)}b)_{(m+k-j)} .
\eqno{(5.2.5)}
$$
This shows that the coefficients of the fields $(\rho a)(z)$
form a Lie algebra.

(ii) A $\BZ$-twisted vertex algebra module is called
simply a vertex algebra module. In particular, the restriction of
a twisted $V$-module $W$ to the vertex subalgebra $V[0]\subset V$
  is a $V[0]$-module. If $G$ arises from an order $N$ automorphism $g$ as in
(5.2.3), then
a $G$-twisted module is called $g$-twisted or twisted by $g$. The
restriction of a $g$-twisted $V$-module to the vertex subalgebra
$V^{g}$ is a $V^{g}$-module.

(iii) It should be clear what the phrases ``$W$ is an irreducible twisted
$V$-module'' and `` $W$ is a twisted $V$-module generated by a collection
of fields $\{(\rho a_{\alpha})(z)$ from a given vector $m\in M$''.

(iv) A vertex algebra is canonically a module over itself. Given an
arbitrary
 $G/\BZ$-graded vertex algebra, there is no obvious way to construct
a $G$-twisted module, but let us consider some concrete examples.

{\bf 5.2.5.} {\it The twisted module
$\Lambda^{ch}\CT(\BC^{N})_{\vec{\lambda},\vec{\mu}}$.}
Given 2 n-tuples $\vec{\lambda}=(\lambda_{0},...,\lambda_{N-1})\in\BC^{N}$,
$\vec{\mu}=(\mu_{0},...,\mu_{N-1})\in\BC^{N}$, let
$G$ be the $\BZ$-span of $\lambda_{i}$, $\mu_{i}$, and 1, $0\leq i\leq N-1$.
Define the $G/\BZ$ grading on $\Lambda^{ch}\CT(\BC^{N})$ by
declaring that
$$
\aligned
&x_{i}\in \Lambda^{ch}\CT(\BC^{N})[-\bar{\lambda_{i}}],\;
\partial_{x_{i}}\in \Lambda^{ch}\CT(\BC^{N})[\bar{\lambda_{i}}],\\
&dx_{i}\in \Lambda^{ch}\CT(\BC^{N})[-\bar{\mu_{i}}],\;
\partial_{dx_{i}}\in \Lambda^{ch}\CT(\BC^{N})[\bar{-\mu_{i}}],
\endaligned
$$
cf. (5.2.2a).

 {\bf 5.2.6. Lemma.}

{\it (i) There is a unique up to isomorphism structure of a $G$-twisted
$\Lambda^{ch}\CT(\BC^{N})$-module }
$$
\rho_{\vec{\lambda},\vec{\mu}}:
\Lambda^{ch}\CT(\BC^{N})\rightarrow \text{Field}_{G}W
$$
generated by $vac\in W$ such that
$$
\aligned
&(\rho_{\vec{\lambda},\vec{\mu}} x_{i})_{(-\lambda_{i}+j)}vac=
(\rho_{\vec{\lambda},\vec{\mu}} \partial_{x_{i}})_{(\lambda_{i}+j)}vac=\\
&(\rho_{\vec{\lambda},\vec{\mu}} (dx_{i}))_{(-\mu_{i}+j)}vac=
(\rho_{\vec{\lambda},\vec{\mu}} \partial_{dx_{i}})_{(\mu_{i}+j)}vac=0.
\endaligned
\eqno{(5.2.6)}
$$
{\it (ii) This module can be constructed by letting
$W=\Lambda^{ch}\CT(\BC^{N})$ as a vector space and}
$$
\aligned
&(\rho_{\vec{\lambda},\vec{\mu}} x_{i})(z)=x_{i}(z)z^{\lambda_{i}},\;
(\rho_{\vec{\lambda},\vec{\mu}} \partial_{x_{i}})(z)=
\partial_{x_{i}}(z)z^{-\lambda_{i}},\\
&(\rho_{\vec{\lambda},\vec{\mu}} dx_{i})(z)=dx_{i}(z)z^{\mu_{i}},\;
(\rho_{\vec{\lambda},\vec{\mu}} \partial_{dx_{i}})(z)=
\partial_{dx_{i}}(z)z^{-\mu_{i}}.
\endaligned
\eqno{(5.2.7)}
$$

{\bf Proof.} The commutation relations (5.2.5) applied to the quadruple
of fields $(\rho_{\vec{\lambda},\vec{\mu}} x_{i})(z)$
$(\rho_{\vec{\lambda},\vec{\mu}} \partial_{x_{i}})(z)$,
$(\rho_{\vec{\lambda},\vec{\mu}} dx_{i})(z)$,
$(\rho_{\vec{\lambda},\vec{\mu}} \partial_{x_{i}})(z)$
imply that their coefficients span the Lie algebra isomorphic
to $Cl\oplus\fa$, see 1.6, 1.7. For example,
$$
[(\rho_{\vec{\lambda},\vec{\mu}} \partial_{x_{i}})_{(\alpha)},
(\rho_{\vec{\lambda},\vec{\mu}} x_{i})_{(\beta)}]=\delta_{\alpha, -\beta-1}.
$$
Conditions (5.2.6) become then the vacuum vector conditions for
$Cl\oplus\fa$ which of course determine $W$ uniquely even as
$Cl\oplus\fa$-module, and uniqueness follows. This fixes recipe
(5.2.7), and a standard argument lucidly explained e.g. in [KR, sect. 5]
allows one to
 extend  (5.2.7) naturally to a twisted module structure
$\rho_{\vec{\lambda},\vec{\mu}}:
\Lambda^{ch}\CT(\BC^{N})\rightarrow \text{Field}_{G}W$. $\qed$

{\bf Notation.} We shall let $\Lambda^{ch}\CT(\BC^{N})_{\vec{\lambda},\vec{\mu}}$ denote the twisted
$\Lambda^{ch}\CT(\BC^{N})$-module occurring in Lemma 5.2.6.

{\bf 5.2.7.} {\it The twisted module $V_{\Lambda+\eta}$.} While
constructing $\Lambda^{ch}\CT(\BC^{N})_{\vec{\lambda},\vec{\mu}}$
requires a little effort, a family of twisted modules over a
lattice vertex algebra seems to be built into the definition of
the latter. For simplicity we shall only consider the case where
the lattice $\Lambda$ is that defined in 3.1. (The following
should be regarded as known even though we failed to find the
needed results in the literature, but see, e.g., a similar and
more involved discussion in [D1].)

Fix $\eta\in\BC\otimes_{\BZ}M^{*}$ and consider the
abelian group $L_{\eta}=\Lambda+\BC\eta\subset \BC\otimes\Lambda$. Let
$\BC[\Lambda_{\eta}]$ be its group algebra. By analogy with
the lattice vertex algebra
$V_{\Lambda}=V(\fh_{L})\otimes\BC_{\epsilon}[\Lambda]$ introduce
$V_{\Lambda_{\eta}}=V(\fh_{L})\otimes\BC[\Lambda_{\eta}]$.

Let
$$
\rho_{\Lambda}:V_{L}\rightarrow \text{Field}V_{\Lambda}
\eqno{(5.2.8)}
$$
be the vertex algebra structure map. Note that
$V_{\Lambda_{\eta}}=V(\fh_{L})\otimes\BC[\Lambda_{\eta}]$ is
naturally a $V(\fh_{L})$-module: indeed  $x(z) 1\otimes
e^{\beta}$, $x\in\fh_{L}$, makes perfect sense if $\beta\in
L_{\eta}$, or indeed if $\beta\in \BC\otimes_{\BZ}L$, see (1.9.5).
Similarly, formula (1.9.2) may be extended without any changes to
define a $\BC_{\epsilon}[\Lambda]$-action on
$V_{\Lambda_{\eta}}=V(\fh_{L})\otimes\BC[\Lambda_{\eta}]$; indeed,
cocycle  $\epsilon(\alpha,\beta)$ defined in (3.1.2) makes perfect
sense for any $\beta\in L+\BC\eta$ because it, so to say, ignores
$\eta$.
 Since any field
$\rho_{\Lambda}a(z)$, $a\in V_{L}$, is written in terms of the
operators we have just described, cf (1.9.4-5), the very formula for $\rho_{\Lambda}a(z)$
defines it as a ``field'' operating on $V_{\Lambda_{\eta}}=V(\fh_{L})\otimes\BC[\Lambda_{\eta}]$. For example,
$$
x(z)e^{\eta}=\left(x(z)+\frac{1}{z}(x,\eta)\b1\right)e^{\eta},\;
x\in\fh_{L}
\eqno{(5.2.9a)}
$$
$$
e^{\alpha}(z)e^{\eta}=\left( e^{\alpha}(z)z^{(\alpha,\eta)}\b1\right)e^{\eta},
\eqno{(5.2.9b)}
$$
where $\b1=e^{0}$ is the vacuum vector of $V_{\Lambda}$.
(Note, by the way, that (5.2.9a-b) are analogous to
spectral flow formulas (1.12.7-8).)
In order to fix the twist, see 5.2.2, let us restrict such fields to
the subspace $V_{\Lambda+\eta}\subset V_{\Lambda_{\eta}}$ where
only $e^{\alpha+\eta}$, $\alpha\in L$, are allowed. Formulas (5.2.9a-b)
allow us to conclude that there is a tautological embedding
$$
\rho_{\Lambda}V_{\Lambda}\hookrightarrow
\text{Field}_{G_{\eta}}V_{\Lambda+\eta},
\eqno{(5.2.10)}
$$
where $G_{\eta}$ is the grading on $V_{\Lambda}$ that assigns
degree
$-(\alpha,\eta)$ to $e^{\alpha}$, $\alpha\in \Lambda$. Hence there
arises the composition
$$
\rho_{\Lambda+\eta}: V_{\Lambda}\buildrel \rho_{\Lambda}\over\rightarrow
\rho_{\Lambda}V_{\Lambda}\hookrightarrow
\text{Field}_{G_{\eta}}V_{\Lambda+\eta},
\eqno{(5.2.11)}
$$
{\bf 5.2.8. Lemma.} {\it Map (5.2.11) endows $V_{\Lambda+\eta}$
with a $G_{\eta}$-twisted $V_{\Lambda}$-module structure.}

{\bf Proof.} Recall that the intuition behind
 twisted Borcherds identity (5.2.4) -- as well
as its untwisted version (1.2.2) -- is that the 3 expressions appearing in
 it, $\rho(a(z-w),w)$, $(\rho a)(z)(\rho b)(w)$, and
$(-1)^{\pr(a)\pr(b)}(\rho b)(w)(\rho a)(z)$ are  Laurent series
expansions of the same function in 3 respective domains, see
 the short discussion after Definition 1.2. This can be made precise:
if $W$ is a vector space graded by finite dimensional subspaces,
 then one can define matrix elements
of fields and their products, such as $<v^{*},\rho(a(z-w),w)v>$,
$<v^{*},(\rho a)(z)(\rho b)(w)v>$, and
$(-1)^{\pr(a)\pr(b)}<v^{*},(\rho b)(w)(\rho a)(z)v>$. If these are
indeed Laurent series expansions of the same rational function
twisted by $z^{-m}$, $w^{k}$, $m$, $k$ not necessarily integral,
with poles on  $z=w$, $z=0$, $w=0$, then (5.2.4) holds; cf. [K,
Remark 4.9a]. In our case  a suitable grading is easy to exhibit
and a familiar argument along the lines of [K, sect.5.4] shows
that matrix elements of the products of fields from
$\rho_{\Lambda+\eta}(V_{\Lambda})$ are indeed such rational
functions. Furthermore, it is easy to see that if we identify
$$
V_{\Lambda+\eta}\iso V_{\Lambda},\; v\mapsto v e^{-\eta},
$$
then these matrix elements, with fixed $v^{*},v\in V_{\Lambda}$, become
analytic functions of $\eta$, see (5.2.9a-b). But if $\eta\in M^{*}$, then
$V_{\Lambda+\eta}=V_{\Lambda}$ by definition and the matrix
elements are indeed equal to each other rational functions.
Thanks to analyticity, this must hold for all $\eta$ and the lemma
follows.

{\bf 5.2.9.} The constructions of 5.2.5 and 5.2.6 are  related
to each other.
Invoke Borisov's algebra $\BB_{\Lambda}$, see (3.1.3). Since
 $\BB_{\Lambda}=V_{\Lambda}\otimes F_{\Lambda}$, the space
$$
\BB_{\Lambda+\eta}\buildrel \text{def}\over =
V_{\Lambda+\eta}\otimes F_{\Lambda}
\eqno{(5.2.12)}
$$
is naturally a twisted $\BB_{\Lambda}$-module. Hence its pull-back
w.r.t. Borisov's embedding
$$
\Lambda^{ch}\CT(\BC^{N})\rightarrow \BB_{\Lambda}
\eqno{(5.2.13)}
$$
is a twisted $\Lambda^{ch}\CT(\BC^{N})$-module.

{\bf 5.2.10. Lemma.} {\it If
$$
\eta=\sum_{j=0}^{N-1}\lambda_{j}X^{*}_{j},
\eqno{(5.2.13)}
$$
then the twisted  $\Lambda^{ch}\CT(\BC^{N})$-module generated by
 $\Lambda^{ch}\CT(\BC^{N})$ from
the vector
$e^{\eta}\otimes\b1\in V_{\Lambda+\eta}$ is isomorphic to
$\Lambda^{ch}\CT(\BC^{N})_{\vec{\lambda}, \vec{\lambda}}$, where
$\vec{\lambda}=(\lambda_{0},...,\lambda_{N-1})$.}

{\bf Proof.} According to Lemma 5.2.5 in order to prove the
lemma one only has to check that the vector $e^{\eta}$
satifies relations (5.2.6). But this is obvious.
For example,  (5.2.9b) gives
$$
\rho_{L+\eta}x_{i}(z)e^{\eta}=
\left(e^{X_{i}}(z)z^{(X_{i},\eta)}\b1\right)e^{\eta}=
\left(e^{X_{i}}(z)z^{\lambda_{i}}\b1\right)e^{\eta},
$$
which is exactly the first of conditions (5.2.6). The remaining
3 fields are dealt with in exactly the same manner.

{\bf 5.2.11. Notation.} It is natural to denote the twisted
module $\Lambda^{ch}\CT(\BC^{N})_{\vec{\lambda}, \vec{\lambda}}$
realized as in Lemma 5.2.10 by
$\Lambda^{ch}\CT(\BC^{N})e^{\sum_{j}\lambda_{j}X_{j}^{*}}$.

{\bf 5.2.12. Landau-Ginzburg orbifold.}
We now wish to orbifoldize the differential graded vertex algebra
$\text{LG}_{f}$, (5.1.2), with respect to the $\BZ_{p}$-action generated by
the operator $g=\exp{(\rho J(z)_{(0)})}$, where $\rho J(z)_{(0)}$ is
 defined in Lemma 5.1.1.
According to the pattern reviewed in 5.2.1, in order to do so one has
to exhibit a $g^{i}$-twisted differential $\text{LG}_{f}$-module
for each $1\leq i\leq N-1$. (Remark 5.2.4 (ii) explains what
``$g^{i}$-twisted'' means)

As follows from (5.2.3), the $i$th twisting gradation is determined by
the action of $\rho J(z)_{(0)}$ on the generators. Formulas of Lemma 5.1.1
imply
$$
\aligned
&\rho J(z)_{(0)}x_{i}=-\frac{1}{p}x_{i},\;
\rho J(z)_{(0)}\partial{x_{i}}=\frac{1}{p}\partial_{x_{i}},\\
&\rho J(z)_{(0)}dx_{i}=(1-\frac{1}{p})dx_{i},\;
\rho J(z)_{(0)}\partial{dx_{i}}=(\frac{1}{p}-1)\partial_{x_{i}}.
\endaligned
\eqno{(5.2.14)}
$$
Therefore,
$\Lambda^{ch}\CT(\BC^{N})_{i\vec{\frac{1}{p}},i\vec{\frac{1}{p}}}$, where
$\vec{\frac{1}{p}}=(1/p,...,1/p)$,
is a natural candidate for the $\exp{i(\rho J(z)_{(0)})}$-twisted
module. So we define, cf. 5.2.11,
$$
\text{LG}^{(i)}_{f}\buildrel \text{def}\over =
\Lambda^{ch}\CT(\BC^{N})e^{\frac{i}{p}\sum_{j}X_{j}^{*}}.
\eqno{(5.2.15)}
$$
As mentioned in 5.2.4 (ii), $\text{LG}^{(i)}_{f}$ is a
(untwisted) $\Lambda^{ch}\CT(\BC^{N})^{g}$-module. (Recall that
$\Lambda^{ch}\CT(\BC^{N})=\text{LG}_{f}=\text{LG}_{f}^{(0)}$
as vertex algebras.)

Since the Landau-Ginzburg differential $df(z)_{(0)}$ comes
from $df\in \Lambda^{ch}\CT(\BC^{N})^{g}$, it operates
naturally on $\text{LG}^{(i)}_{f}$ thus making it a differential
$\text{LG}_{f}^{g}$-module.

Since  (5.1.4) maps $V(N2)_{N\frac{p-2}{p}}$  into $\text{LG}_{f}^{g}$, each
$\text{LG}^{(i)}_{f}$ acquires an $N2$-structure. In particular,
$g$ operates on $\text{LG}^{(i)}_{f}$ and $\left(\text{LG}^{(i)}_{f}\right)^{g}$
is also a differential $\left(\text{LG}^{(0)}_{f}\right)^{g}$-module.

{\bf 5.2.13. Definition.} Define the Landau-Ginzburg orbifold
to be the following differential $\text{LG}_{f}^{g}$-module:
$$
\text{LG}_{f,\text{orb}}=\left(\bigoplus_{j=0}^{p-1}
\left(\text{LG}^{(j)}_{f}\right)^{g},\; df(z)_{(0)}\right).
$$
$\qed$

Note that if $\text{deg}f=p=N$, then $\text{LG}_{f,\text{orb}}$
simply coincides with
$(\widetilde{\text{LG}}_{orb}^{<N},df(z)_{(0)})$ appearing in
Theorems 4.7, 4.9. This space, however, carries two {\it a priori}
different $N2$-structures: one computed in Lemma 4.10.1 and having
purely geometric origin and another, Landau-Ginzburg structure,
recorded in Lemma 5.1.1.

{\bf 5.2.14. Lemma.} {\it The two $N2$-structures coincide with each other.}

\bigskip

{\bf Proof} is, of course, a routine computation consisting in
applying Borisov's formulas (3.3.2a-b) to the 4 fields of Lemma 5.1.1.
Let us consider  $Q(z)$ and leave the rest as an exercise for the interested
reader. We have (and recall that we are using boson-fermion
correspondence (4.10.1)):
$$
\aligned
&Q(z)\mapsto \sum_{i=0}^{N-1}-\frac{1}{N}x_{i}(z)\partial_{dx_{i}}(z)'
-(\frac{1}{N}-1)x_{i}(z)'\partial_{dx_{i}}(z)\mapsto\\
&\sum_{i=0}^{N-1}-\frac{1}{N}:e^{X_{i}}(z)
\left(:\left(-X_{i}(z)-\chi_{i}(z)\right)e^{-X_{i}-\chi_{i}}(z):\right):-\\
&\left(\frac{1}{N}-1\right):\left(:X_{i}(z)e^{X_{i}}(z):\right)
e^{-X_{i}-\chi_{i}}(z):=\\
&\sum_{i=0}^{N-1} \frac{1}{N}:\left(X_{i}\left(z\right)+\chi_{i}\left(z\right)
\right)e^{-\chi_{i}}\left(z\right):+ \left(1-\frac{1}{N}\right):X_{i}\left(z\right)
e^{-\chi_{i}}\left(z\right):=\\
&\sum_{i=0}^{N-1} \frac{1}{N}:\chi_{i}(z)e^{-\chi_{i}}(z):
+:X_{i}(z)e^{-\chi_{i}}(z):,
\endaligned
$$
as it should, cf. Lemma 4.10.1.  Note that this computation is
parallel to that we performed in the proof of Lemma 4.10.1. $\qed$

Lemma 5.2.14 represents the last step in the identification of
$\widetilde{\text{LG}}^{<N}_{orb}$ that played an important role
in Theorems 4.7, 4.9 with the Landau-Ginzburg
orbifold\footnote{the tensor product decomposition of
$H_{df(z)_{(0)}}(\widetilde{\text{LG}}^{<N}_{orb})$ arising in the
case of a diagonal equation, see Corollary 5.1.4,
 is perhaps a bridge to Gepner's
models [G] also cooked up of the tensor products of $N2$ unitary
representations}.
The next lemma
says that as $N2$-modules the twisted sectors are spectral flow
transforms of the untwisted one, see (1.12.1-2) for the definition
of the spectral flow.

{\bf 5.2.15. Lemma.} {\it The map}
$$
e^{-j\sum_{i}(\frac{1}{p}X^{*}_{i}-X_{i}-\chi_{i})}:\;
\text{LG}_{f}^{(j)}\rightarrow \text{LG}_{f};\text{ s.t. }
x\mapsto xe^{-j\sum_{i}(\frac{1}{p}X^{*}_{i}-X_{i}-\chi_{i})}.
\eqno{(5.2.16)}
$$
{\it delivers the following isomorphisms:

(i)}
$$
\text{LG}^{(j)}_{f}\iso
S_{j}\left(\text{LG}_{f}\right),
H_{df(z)_{(0)}}(\text{LG}^{(j)}_{f})\iso
S_{j}\left(H_{df(z)_{(0)}}(\text{LG}_{f})\right)\;.
$$
{\it (ii) If $p=N$, then}
$$
\left(\text{LG}^{(j)}_{f}\right)^{g}\iso
S_{j}\left(\text{LG}_{f}^{g}\right),\;
H_{df(z)_{(0)}}\left(\text{LG}^{(j)}_{f}\right)^{g}\iso
S_{j}\left(H_{df(z)_{(0)}}(\text{LG}_{f}^{g})\right)
$$
{\it if $p=N$.}

{\bf Proof} is  parallel to the discussion in 3.10, and we will be
brief: a glance at the formulas in Lemma 4.10.1 shows that under  map
(5.2.16)
$Q(z)$ gets multiplied by $z^{j}$, $G(z)$ by $z^{-j}$ thereby recovering
(1.12.1). This proves the first of ismorphisms (i); the second
follows from an equally obvious observation that the differential
$df(z)_{(0)}$ is invariant under (5.2.16).

As to (ii), this argument does not quite apply if $p\neq N$ because then
\newline
$e^{-j\sum_{i}(\frac{1}{p}X^{*}_{i}-X_{i}-\chi_{i}})$ is not $Z_{N}$-invariant
and does not give a well-defined map. If, however, $p=N$, then
(ii) is proved in the same way as (i).

{\bf 5.2.16. Application: the character formulas.} Lemma 5.2.15
is an intelligent way to obtain (4.11.4-5). Focus on the Euler
character formulas. First of all, $\text{LG}_{f}$ being bi-graded
by the eigenvalues of $\rho J_{(0)}$ and $\rho L_{(1)}$,
definition (4.11.2) is equivalent to
$$
\text{Eu}(\text{LG}_{f})(s,\tau)=\text{Tr}|_{\text{LG}_{f}}(-1)^{\pr(.)}
e^{2\pi i(sJ_{(0)}+\tau L_{(1)})},
\eqno{(5.2.17)}
$$
where $\pr(.)$ is the parity function, see 1.1 -- and we are skipping the
$\rho$.
It is quite standard to deduce that
$$
\text{Eu}(\text{LG}_{f})(s,\tau)= E(\tau,s),
\eqno{(5.2.18)}
$$
where $E(\tau,s)$ is the function appearing right before (4.11.4), after
all $\text{LG}_{f}$ is a superpolynomial ring in infinitely many variables.
 In view
of (5.2.16), to obtain the Euler character of a spectral flow transformed
module, one has to perform the linear change of variables, determined by
(1.12.1), in the Euler
character of the original module and then take
care of the parity. By virtue of Lemma 5.2.15, in the case of
 $\text{LG}_{f}^{(j)}$ one has to replace
$$
\aligned
J_{(0)} &\text{ with } J_{(0)}-(N-2)j,\\
L_{(1)} &\text{ with }L_{(1)}-jJ_{(0)}+j(j-1)(N-2)/2,\\
\pr &\text{ with }\pr +||j\sum_{i}\frac{1}{N}X^{*}_{i}-X_{i}-\chi_{i}||^{2}
\text{ mod } 2.
\endaligned
\eqno{(5.2.19)}
$$
The 1st two of these follow from (1.12.1), the last from (5.2.16)
and definition of parity (1.9.3). Plugging (5.2.19) in (5.2.17) we
obtain
$$
\aligned
\text{Eu}(\text{LG}_{f}^{(j)})(s,\tau)=&\text{Tr}|_{\text{LG}_{f}}
(-1)^{\pr(.)+j^{2}(N-2)}\times\\
&e^{2\pi i(s(J_{(0)}-(N-2)j)+\tau (L_{(1)}-jJ_{(0)}+j(j-1)(N-2)/2))}.
\endaligned
\eqno{(5.2.20)}
$$
Then (5.2.18) can be rewritten as follows
$$
\text{Eu}(\text{LG}_{f}^{(j)})(s,\tau)=
e^{\pi
i\left(N-2\right)\left\{-2js+\left(j^{2}-j\right)\tau+j^{2}\right\}}
E(\tau,s-j\tau).
\eqno{(5.2.21)}
$$
Finally, in order to extract  $g$-invariants, one applies
the averaging operator $\frac{1}{N}\sum_{l}e^{2\pi i l J_{(0)}}$
which results in shifts $s\mapsto s +l$ and gives
(4.11.4).

{\bf 5.2.17. Vertex algebra structure.} Let us focus on the case
where $\text{deg}f=p=N$. {\it A priori} the orbifold
$\text{LG}_{f,\text{orb}}$ carries only uninteresting vertex
algebra structure: by definition it can be regarded as an abelian
extension of $\text{LG}_{f}^{g}$ by the sum of vertex modules
$\oplus_{i}\left(\text{LG}^{(i)}_{f}\right)^{g}$. Cut-off
definition (4.6.13), which applies to the present situation, see
5.2.13, fares somewhat better: it endows
$\text{LG}_{f,\text{orb}}$ with a filtered vertex algebra
structure such that the corresponding graded object is the
mentioned above uninteresting vertex algebra. And yet this vertex
algebra structure does not seem to be quite right either.
Formally, one would expect  an orbifold multiplicative structure
to be $\BZ_{N}$-graded. Instead, our vertex algebra is graded by
the semigroup $\BZ_{+}/(N+\BZ_{+})$. A similar problem has already
manifested itself: by virtue of (4.9.3) the chiral ring of
$\text{LG}_{f, orb}$ is isomorphic to the cohomology algebra
$H_{\fF}(\fF,\Lambda^{*}\CT_{\fF})$ as a vector space but not as a
ring, see (4.13.12) and the sentence that follows. Let us propose
a  way out in the case where
 $f=\sum_{i}x_{i}^{N}$ and $\text{LG}_{f, orb}$ is replaced with
the Landau-Ginzburg cohomology  $H_{df(z)_{(0)}}(\text{LG}_{f,
orb})$.

In this case, Corollary 5.1.4 combined
with Lemma 5.2.15 provides us with a very explicit description
of the cohomology:
$$
H_{df(z)_{(0)}}(\text{LG}_{f,orb})=\bigoplus_{n=0}^{N-1}
\bigoplus_{j_{0},...,j_{N-1}\leq N-2}
\bigotimes_{t=0}^{N-1}S_{n}\left(U_{j_{t},N}\right)^{g}.
\eqno{(5.2.22)}
$$
The point is, the category of unitary highest weight $N2$-modules
has a remarkable periodicity property: there are ismorphisms [FS]
$$
T_{m}: S_{mN+n}(U_{j,N})\iso S_{n}(U_{j,N}),\; m\in\BZ,\;0\leq
n\leq N-1.
\eqno{(5.2.23)}
$$
Here is what this ``Bott periodicity'' suggests: since
 $\text{LG}_{f}^{(n)}$ is naturally defined for any
$n\in\BZ$,  instead of taking $\oplus_{n\geq
0}H_{df(z)_{(0)}}(\text{LG}_{f}^{(n)})^{g}$ and then quotienting
out by the ideal $\oplus_{n\geq
N}H_{df(z)_{(0)}}(\text{LG}_{f}^{(n)})^{g}$, consider
$\oplus_{n\in \BZ}H_{df(z)_{(0)}}(\text{LG}_{f}^{(n)})^{g}$ and
identify the components
$H_{df(z)_{(0)}}(\text{LG}_{f}^{(\alpha)})^{g}$ and
$H_{df(z)_{(0)}}(\text{LG}_{f}^{(\beta)})^{g}$ with
$\alpha-\beta\in N\BZ$  by using (5.2.23). Note that (5.2.23) is
determined up to proportionality and herein lies the difficulty.
Nevertheless, the following holds true.

{\bf 5.2.18. Theorem.} {\it Collection of isomorphisms (5.2.23)
can be normalized so that the linear span
$$
\CJ=<x-T_{m}^{\otimes
N}x;\;x\in\bigotimes_{t=0}^{N-1}S_{mN+n}\left(U_{j_{t},N}\right)^{g},\;0\leq
j_{t}\leq N-2,0\leq n\leq N-1>
$$
is an ideal of the vertex algebra $\oplus_{n\in
\BZ}H_{df(z)_{(0)}}(\text{LG}_{f}^{(n)})^{g}$.}

\bigskip

Proof of this theorem requires a little plunge in the genus zero
$N2$ correlation functions and will be skipped. Note that there is
an obvious vector space isomorphism
$$
H_{df(z)_{(0)}}(\text{LG}_{f, orb})\iso \bigoplus_{n\in
\BZ}H_{df(z)_{(0)}}(\text{LG}_{f}^{(n)})^{g}/\CJ,
 \eqno{(5.2.24)}
$$
and according to Theorem 5.2.18 the vector space
$H_{df(z)_{(0)}}(\text{LG}_{f, orb})$ becomes a vertex algebra.
The following corollary suggests that this vertex algebra
structure is the ``right'' one.

\begin{sloppypar}
{\bf 5.2.19. Corollary.} {\it The chiral ring of
$H_{df(z)_{(0)}}(\text{LG}_{f, orb})$ is isomorphic to the ring
$H_{\fF}(\fF,\Lambda^{*}\CT_{\fF})$, where $\fF$ is the Calabi-Yau
hypersurface (4.1.6).}
\end{sloppypar}

\bigskip

{\bf 5.2.20.} {\it Proof.} At the vector space level, Corollary
5.2.19 is nothing but Theorem 4.9 (ii) or, rather, its proof, see
sect. 4.13.1, formula (4.13.2), which shows that the differential
$D_{orb}$ is zero on the chiral ring. The multiplicative structure
of $H_{\fF}(\fF,\Lambda^{*}\CT_{\fF})$ has already been detected
inside $\widetilde{\text{LG}}_{orb}$ in sect. 4.13.2. It is clear
that the part of that argument which has to do with the untwisted
sector (i.e. the Milnor ring) or twisted sectors of ``low'' degree
of twisting, see (4.13.10-4.13.11), carries over to the present
situation without any changes. What remains to be done is to check
that (4.13.13) holds in $H_{df(z)_{(0)}}(\text{LG}_{f, orb})$.

In order to verify (4.13.13), we have to rely on (4.13.12).
 For the reader's convenience,
let us copy the latter:
$$
\aligned
&\left(e^{sX^{*}_{orb}-\sum_{j}(X_{i}+\chi_{i})}\right)_{(-1)}
e^{(N-s)X^{*}_{orb}-\sum_{j}(X_{i}+\chi_{i})}= \\
&(-1)^{s}
e^{NX^{*}_{orb}-2\sum_{j}(X_{i}+\chi_{i})}
\in \text{LG}_{f,orb}^{(N)}.
\endaligned
\eqno{(5.2.25)}
$$
We assert that up to proportionality
$$
T_{1}^{\otimes N} \left(\text{class of }
e^{\sum_{j}(X^{*}_{j}-2(X_{j}+\chi_{j}))}\right)= \text{class of
}e^{(N-2)\sum_{j}X_{j}} \in
H_{df(z)_{(0)}}(\text{LG}_{f,orb}^{(N)}). \eqno{(5.2.26)}
$$
If proven, (5.2.26) will recover the desired (4.13.13).

In order to prove (5.2.26), we shall obtain a representation
theoretic interpretation of both  the elements that determines
them up to proportionality.

Since in the case at hand everything is a tensor product of 1-dimensional
Landau-Ginzburg models and their cohomology groups,
we can drop the subindex $j$, and our task is
to show that
$$
T_{1}(\text{ class of }e^{X^{*} -2(X +\chi )})= \text{class of
}e^{(N-2)X}\in H_{df(z)_{(0)}}(\text{LG}_{f}).
\eqno{(5.2.27)}
$$
According to Theorem 5.1.3,  $e^{(N-2)X}$ generates in
the cohomology the unitary $N2$-module $U_{N-2,N}$. Actually,
$e^{(N-2)X}$ is a {\it highest weight vector} of $U_{N-2,N}$
meaning that it is annihilated by the Lie algebra span of
$$
\aligned
&G_{(2)},G_{(3)},...,\\
Q_{(1)}, &Q_{(2)}, Q_{(3)}.
\endaligned
\eqno{(5.2.28)}
$$
It follows from Lemma 5.2.15 (i) that
$S_{N}(U_{N-2,N})\subset H_{df(z)_{(0)}}(\text{LG}^{(N)}_{f})$
is generated by $e^{X^{*}  -2X-N\chi }$. The same lemma implies that
the latter is rather a {\it twisted highest weight vector} meaning
that
it is annihilated by the Lie algebra span of
$$
\aligned
&G_{(N+2)},G_{(N+3)},...,\\
Q_{(-N+1)}, Q_{(-N+2)},..., &Q_{(N+2)}, Q_{(N+3)},...
\endaligned
\eqno{(5.2.29)}
$$
Now we would like to  find a highest weight vector of
$S_{N}(U_{N-2,N})\subset H_{df(z)_{(0)}}(\text{LG}^{(N)}_{f})$. A
direct computation shows that $e^{X^{*}  -2X-N\chi }$ is
annihilated by $G_{(N+1)}$ in addition to  subalgebra (5.2.29).
Commutation relations (1.10.1) imply then
 that the vector
$$
G_{(3)}G_{(4)}\cdots G_{(N-1)}G_{(N)}e^{X^{*}  -2X-N\chi }
\eqno{(5.2.30)}
$$
is annihilated by the untwisted annihilating subalgebra, (5.2.28),
except perhaps the element $G_{(2)}$. It also follows from [FS,
(A.5)] that the class of (5.2.30) is non-zero in
$S_{N}(U_{N-2,N})$. Now we compute
$$
\aligned
&G_{(N)}e^{X^{*}  -2X-N\chi }=
\left(:(X^{*}(z)-\chi(z))e^{\chi}(z):\right)_{(N)}e^{X^{*}  -2X-N\chi }=\\
&e^{\chi}(X^{*}_{(0)}-\chi_{(0)})e^{X^{*}  -2X-N\chi }=
(N-2)e^{X^{*}  -2X-(N-1)\chi }.
\endaligned
$$
Likewise,
$$
G_{(N-1)}G_{(N)}e^{X^{*}  -2X-N\chi }=
(N-3)(N-2)e^{X^{*}  -2X-(N-2)\chi }.
$$
Continuing in the same vein we obtain
$$
G_{(3)}G_{(4)}\cdots G_{(N-1)}G_{(N)}e^{X^{*}  -2X-N\chi } =(N-2)!
e^{X^{*}  -2X-2\chi }
\eqno{(5.2.31)}
$$
The same computation shows that an application of $G_{(2)}$
annihilates  vector (5.2.31). Hence  (5.2.31) is a highest weight
vector and, since it is
 determined by this condition up to proportionality, (5.2.27) follows.
The $N$th tensor power of the latter gives the desired (5.2.26).
$\qed$

\bigskip

\bigskip

\centerline{{\bf References}}

[BD] A. ~Beilinson, V. ~Drinfeld, Chiral algebras, Preprint.

[B] L.Borisov, Vertex algebras and mirror symmetry, {\it Comm. Math. Phys.}
{\bf 215} (2001), no. 3, 517-557

[BL] L.Borisov, A.Libgober, Elliptic genera of toric varieties and applications
to mirror symmetry, {\it Inv. Math.} {\bf 140} (2000) no. 2, 453-485,

[D1] C.Dong, Vertex algebras associated with
even lattices, {\it J. Algebra} {\bf 161} (1993),  245-265

[D2] C.Dong, Twisted modules for vertex algebras associated with
even lattices, {\it J. Algebra} {\bf 165} (1994),  91-112

[Don] R.Donagi, Generic Torelli for projective hypersurfaces,
{\it Comp. Math.} {\bf 50} (1983) 325-353

[FS] B.L.Feigin, A.M.Semikhatov, Free-field resolutions of the unitary
N=2 super-Virasoro representations, preprint, hep-th/9810059

[FFR] A.Feingold, I.Frenkel, J.Ries, Spinor constructions of vertex
operator algebras, tiality, and $E^{(1)}_{8}$, {\it Contemp. Math.},
{\bf 121}, AMS, 1991

[FB-Z] E.Frenkel, D.Ben-Zvi, Vertex algebras and algebraic curves,
{\it Mathematical Surveys and Monographs} {\bf 88}, 2001,

%[FS] E.Frenkel, M.Szczesny, Chiral de Rham complex
%and orbifolds, preprint, posted on the net:
%AG/0307181

[FLM] I.Frenkel, J.Lepowski, A.Meurman, Vertex operator algebras
and the Monster, {\it Academic Press}, 1988

[G] D.Gepner,  Exactly solvable string compactifications on manifolds
of $SU(N)$ holonomy, {\it Phys.Lett.} {\bf B199} (1987)380-388,

[GMS] V.Gorbounov, F.Malikov, V.Schechtman,
Gerbes of chiral differential operators. II, to appear in
{\it Inv. Math.},  AG/0003170

[Gr] P.Griffiths, On the periods of certain rational integrals,
{\it Ann. Math.} {\bf 90} (1969) 460-541,

[K] V.Kac, Vertex algebras for beginners, 2nd edition, {\it AMS}, 1998

[KR] V.Kac, A.Radul, Representation theory of the vertex algebra
$W_{1+\infty}$, {\it Transform. groups}. {\bf 1}
(1996), no. 1-2, 41-70

[KV1] M.Kapranov, E.Vasserot,
Vertex algebras and the formal loop space, preprint, math.AG/0107143,

[KV2] M.Kapranov, E.Vasserot, private communication

[KO] A.Kapustin, D.Orlov, Vertex algebras, mirror
symmetry, and D-branes: the case of complex tori,
{\it Comm. Math. Phys.} {\bf 233} (2003) 79-163,

[KYY] T.Kawai, Y.Yamada, S-K. Yang, Elliptic genera
 and N=2 superconformal field theory, {\it Nucl.Phys.}
 {\bf B414} (1994) 191-212

[LVW] W.Lerche, C.Vafa, N.Warner, Chiral rings in N=2 superconformal theory,
{\it Nuclear Physics} {\bf B324} (1989), 427,

[MS] F.Malikov, V.Schechtman, Deformations of vertex algebras, quantum
cohomology of toric varieties, and
elliptic genus, {\it Comm. Math. Phys.}
{\bf 234} (2003), no.1 77-100,

[MSV] F. ~Malikov, V. ~Schechtman, A. ~Vaintrob, Chiral de Rham complex,
{\it Comm. Math. Phys.}, {\bf 204} (1999), 439 - 473,

[OR] K. Ono, S-S. Roan,  Vafa's formula and equivariant
$K$-theory, {\it J. Geom. Phys. } {\bf 10}  (1993), no. 3,
287--294.

 [R] S-S. Roan,  On Calabi-Yau orbifolds in weighted projective
spaces, {\it Internat. J. Math.} {\bf 1} (1990), no. 2, 211--232.

[S] A.Schwarz, Sigma-models having supermanifolds as target
spaces, {\it Lett. Math. Phys.} {\bf 38} (1996), no.4, 349-353

[V] C.Vafa, String vacua and orbifoldized Landau-Ginzburg model,
 {\it Mod.Phys.Lett.} {\bf A4} (1989) 1169

[VW] C.Vafa, N.P.Warner, Catastrophes and the classification of conformal
theories, {\it  Phys.Lett.} {\bf B218} (1989) 51

[W1] E. Witten, Phases of N=2 theories in two-dimensions,
{\it Nucl.Phys.} {\bf B403}, 159-222,1993,
 hep-th/9301042

[W2] E. Witten, On the Landau-Ginzburg description of N=2 minimal
models, {\it Int.J.Mod.Phys.} {\bf A9} (1994), 4783-4800, hep-th/9304026

\bigskip\bigskip

V.G.: Department of Mathematics, University of Kentucky,
Lexington, KY 40506, USA;\ vgorb\@ms.uky.edu

F.M.: Department of Mathematics, University of Southern California,
Los Angeles, CA 90089, USA;\ fmalikov\@mathj.usc.edu

\end{document}